\documentclass[amssymb,12pt]{amsart}
\usepackage{latexsym}
\usepackage{amsfonts}
\usepackage{amssymb}
\usepackage{graphics}
\usepackage{lscape}
\usepackage{rotating}
\usepackage{multirow}

\newtheorem{propo}{Proposition}[section]
\newtheorem{proposition}[propo]{Proposition}

\newtheorem{lemma}[propo]{Lemma}

\newtheorem{corollary}[propo]{Corollary}

\newtheorem{theorem}[propo]{Theorem}

\newcommand{\Ker}{\operatorname{Ker}}

\newcommand{\Irr}{{\mathrm {Irr}}}

\newcommand{\Spec}{{\mathrm {Spec}}}

\newcommand{\ZZ}{{\mathbb Z}}

\newcommand{\FF}{{\mathbb F}}

\newcommand{\tn}{\hspace{0.5mm}^{t}\hspace*{-1mm}}

\newcommand{\p}{^{\prime}}

\begin{document}
\title[low-dimensional complex characters of finite classical groups]
{Low-dimensional complex characters of the symplectic and orthogonal
groups}

\author{Hung Ngoc Nguyen}
\address{Department of Mathematics, Michigan State University, East Lansing,
MI 48824, USA} \email{hungnguyen@math.msu.edu}

\thanks{Part of this work was done during my PhD study at the University of Florida. I am grateful
to my advisor, Professor P. H. Tiep, for the devoted guidance and
fruitful discussion. I would like to thank Professor M. W. Liebeck
for motivating this work. Also, I want to thank Professor J. I. Hall
and Professor U. Meierfrankenfeld for their generous hospitality and
support during my time at the Michigan State University. Finally, I
thank the referee for his careful reading of the manuscript and
helpful suggestions. }

\date{\today}

\begin{abstract}
We classify the irreducible complex characters of the symplectic
groups $Sp_{2n}(q)$ and the orthogonal groups $Spin_{2n}^\pm(q)$,
$Spin_{2n+1}(q)$ of degrees up to the bound $D$, where
$D=(q^n-1)q^{4n-10}/2$ for symplectic groups, $D=q^{4n-8}$ for
orthogonal groups in odd dimension, and $D=q^{4n-10}$ for orthogonal
groups in even dimension.
\end{abstract}

\maketitle

\section{Introduction}

Lower bounds for the degrees of non-trivial irreducible
representations of finite groups of Lie type in cross characteristic
were found by Landazuri, Seitz, and Zalesskii in \cite{LS},
\cite{SZ}. These bounds have proved to be very useful in various
applications. Thanks to the effort of many people, the exact value
of the smallest degree of non-trivial irreducible representations
was determined for many groups of Lie type (see \cite{T}).

We are interested in not only the smallest representation, but more
importantly, the low-dimensional representations. In the case of
complex representations, this problem was first studied by Tiep and
Zalesskii \cite{TZ1} for finite classical groups and then by
L\"{u}beck \cite{Lu} for exceptional groups. For representations
over fields of cross characteristics, this problem has been studied
recently in \cite{BK}, \cite{GT1} for $SL_n(q)$; \cite{GMST},
\cite{HM} for $SU_n(q)$; \cite{GMST} for $Sp_{2n}(q)$ with $q$ odd;
and \cite{GT2} for $Sp_n(q)$ with $q$ even. Let $H$ be one of these
groups and denote by $\mathfrak{d}(H)$ the smallest degree of
nontrivial irreducible representations in cross characteristic. The
purpose of these papers is to classify the irreducible
representations of $H$ of degrees close to $\mathfrak{d}(H)$ and to
prove that there is a relatively ``big" gap between the degrees of
these representations and the next degree.

From now on, we denote by $G$ either the symplectic groups
$Sp_{2n}(q)$ or the orthogonal groups $Spin_{n}^\pm(q)$, where $q$
is a power of a prime $p$. The smallest non-trivial complex
characters of $G$ was determined in \cite{TZ1}. It was also proved
in \cite{TZ1} that, up to the bound $(q^{2n}-1)/2(q+1)$, the odd
characteristic symplectic group $Sp_{2n}(q)$ has four irreducible
characters of degrees $(q^n\pm1)/2$, which are so-called \emph{Weil
characters}, and the smallest unipotent character of degree
$(q^n-1)(q^n-q)/2(q+1)$. We want to extend these results to a larger
bound. More precisely, we classify the irreducible complex
characters of $G$ of degrees up to the bound $D$, where
$D=(q^n-1)q^{4n-10}/2$ for the odd characteristic symplectic groups,
$D=q^{4n-8}$ for the orthogonal groups in odd dimension, and
$D=q^{4n-10}$ for the orthogonal groups in even dimension. When $q$
is even, the irreducible complex characters of $Sp_{2n}(q)\simeq
Spin_{2n+1}(q)$ were classified up to an already good enough bound
$(q^{2n}-1)(q^{n-1}-1)(q^{n-1}-q^2)/2(q^4-1)$ in \cite{GT2} and
therefore it is not in our consideration. Theorems
\ref{orthogonal-odddimension}, \ref{orthogonal-qeven},
\ref{orthogonal-qodd}, and Proposition \ref{smallcase} in this paper
have been used in \cite{LOST} to establish the Ore's conjecture for
the orthogonal groups. Our main results are the following theorems.

\begin{theorem}\label{symplectic}
Let $\chi$ be an irreducible complex character of $G=Sp_{2n}(q)$,
where $n\geq 6$ and $q$ is an odd prime power. Then either
$\chi(1)>(q^n-1)q^{4n-10}/2$ or $\chi$ belongs to a list of
$q^2+12q+36$ characters in Tables 1, 2, 3 (at the end of $\S4$).
\end{theorem}

\begin{theorem}\label{orthogonal-odddimension}
Let $\chi$ be an irreducible complex character of
$G=Spin_{2n+1}(q)$, where $n\geq 5$ and $q$ is an odd prime power.
Then $\chi(1)=1$ ($1$ character), $\chi(1)=(q^{2n}-1)/(q^2-1)$ ($1$
character), $\chi(1)=q(q^{2n}-1)/(q^2-1)$ ($1$ character),
$(q^n+\alpha_1)(q^n+\alpha_2q)/2(q+\alpha_1\alpha_2)$ ($4$
characters, $\alpha_{1,2}=\pm1$), $\chi(1)=(q^{2n}-1)/(q+\alpha)$
($(q+\alpha-2)/2$ characters for each $\alpha=\pm1$), or
$\chi(1)>q^{4n-8}$.
\end{theorem}

\begin{theorem}\label{orthogonal-qeven}
Let $\chi$ be an irreducible complex character of
$G=Spin_{2n}^\alpha(q)\simeq\Omega^{\alpha}_{2n}(q)$, where $n\geq
5, \alpha=\pm$, $(n,q,\alpha)\neq(5,2,+)$, and $q$ is a power of
$2$. Let
$$D(n,q,\alpha)=\left\{\begin {array}{ll}
(q-1)(q^2+1)(q^3-1)(q^4+1), & n=5, \alpha=-,\\
q^{4n-10}+1, & otherwise.
 \end {array} \right.$$ Then $\chi(1)=1$ ($1$ character),
$\chi(1)=(q^{n}-\alpha)(q^{n-1}+\alpha q)/(q^2-1)$ ($1$ character),
$\chi(1)=(q^{2n}-q^2)/(q^2-1)$ ($1$ character),
$\chi(1)=(q^n-\alpha)(q^{n-1}+\alpha\beta)/(q-\beta)$
($(q-\beta-1)/2$ characters for each $\beta=\pm$), or $\chi(1)\geq
D(n,q,\alpha)$.

Furthermore, when $(n,\alpha)=(5,-)$, $G$ has exactly $q$ characters
of degree $D(5,q,-)$, $1$ character of degree
$q^2(q^4+1)(q^5+1)/(q+1)$, and no more characters of degrees up to
$q^{10}$.
\end{theorem}

\begin{theorem}\label{orthogonal-qodd}
Let $\chi$ be an irreducible complex character of
$G=Spin_{2n}^\alpha(q)$, where $\alpha=\pm$, $n\geq 5$, and $q$ is
an odd prime power. Let
$$D(n,q,\alpha)=\left\{\begin {array}{ll}
(q-1)(q^2+1)(q^3-1)(q^4+1), & n=5, \alpha=-,\\
q^{4n-10}+1, & otherwise.
 \end {array} \right.$$
Then $\chi(1)=1$ ($1$ character),
$\chi(1)=(q^n-\alpha)(q^{n-1}+\alpha q)/(q^2-1)$ ($1$ character),
$\chi(1)=(q^{2n}-q^2)/(q^2-1)$ ($1$ character),
$\chi(1)=(q^{n}-\alpha)(q^{n-1}-\alpha)/2(q+1)$ ($2$ characters),
$\chi(1)=(q^{n}-\alpha)(q^{n-1}+\alpha)/2(q-1)$ ($2$ characters),
$\chi(1)=(q^{n}-\alpha)(q^{n-1}+\alpha\beta)/(q-\beta)$
($(q-\beta-2)/2$ characters for each $\beta=\pm$), or $\chi(1)\geq
D(n,q,\alpha)$.

Furthermore, when $(n,\alpha)=(5,-)$, $G$ has exactly $q$ characters
of degree $D(5,q,-)$, $1$ character of degree
$q^2(q^4+1)(q^5+1)/(q+1)$, and no more characters of degrees up to
$q^{10}$.
\end{theorem}

This paper is organized as follows. In $\S2$, we recall the
Lusztig's classification of the irreducible complex characters of
finite groups of Lie type. Also, the structures of centralizers of
semi-simple elements in symplectic and orthogonal groups are
collected in $\S2$. The low-dimensional unipotent characters of $G$
are classified in $\S3$. Each of the groups $Sp_{2n}(q)$ ($q$ odd),
$Spin_{2n+1}(q)$ ($q$ odd), $Spin_{2n}^\pm(q)$ ($q$ even), and
$Spin_{2n}^\pm(q)$ ($q$ odd) is treated individually in $\S4$,
$\S5$, $\S6$, and $\S7$, respectively. Finally, in $\S8$, we count
the number of irreducible complex characters of $G=Spin_{12}^\pm(3)$
of degrees up to $4\cdot 3^{15}$.

%%% -----------------------------------------------------------------------------------

\section{Preliminaries}

\subsection{Notation} Throughout this paper, $GL^+_n(q)$ stands for
$GL_n(q)$ and $GL^-_n(q)$ stands for $GU_n(q)$. Furthermore, $\pm$
could be understood as $\pm$ itself or $\pm1$, depending on the
context. In order to label the eigenvalues of semi-simple elements,
as well as conjugacy classes of finite classical groups, we follow
the notation of \cite{S} and \cite{W2}. Let $\kappa$ denote a
generator of the field of $q^4$ elements, $\zeta=\kappa^{q^2-1}$,
$\theta=\kappa^{q^2+1}$, $\eta=\theta^{q-1}$, $\gamma=\theta^{q+1}$.
Let $T_1=\{1,...,(q-3)/2\}$, $T_2=\{1,...,(q-1)/2\}$ if $q$ is odd
and $T_1=\{1,...,(q-2)/2\}$, $T_2=\{1,...,q/2\}$ if $q$ is even.
Furthermore, when $q$ is odd, let $R_1^\ast=\{j\in \mathbb{Z}: 1\leq
j< q^2+1, j\neq (q^2+1)/2\}$ and define $R_1$ to be a complete set
of class representatives of the equivalence relation $\sim$ on
$R_1^\ast$: $$i\sim j \Leftrightarrow i\equiv \pm j ~{\text{or}} \pm
qj~(\bmod~ q^2+1).$$ Similarly, let $R_2^\ast=\{j\in \mathbb{Z}:
1\leq j\leq q^2-1, q-1\nmid j, q+1\nmid j\}$ and define $R_2$ to be
a complete set of class representatives of the equivalence relation
$\sim$ on $R_2^\ast$: $$i\sim j \Leftrightarrow i\equiv \pm j
~{\text{or}} \pm qj~(\bmod~ q^2-1).$$ Note that $|R_1|=(q^2-1)/4$
and $|R_2|=(q-1)^2/4$.

\subsection{Lusztig's classification} The main ingredient to prove
our results is Lusztig's classification of irreducible complex
characters of finite groups of Lie type. To make the paper
self-contained, we recall it here.

Let $G$ be either the symplectic groups $Sp_{2n}(q)$ or the
orthogonal groups $Spin_{2n}^\pm(q)$, $Spin_{2n+1}(q)$ where $q$ is
a power of a prime $p$. Let $\mathcal{G}$ be the algebraic group and
$F$ the Frobenius endomorphism on $\mathcal{G}$ such that
$G=\mathcal{G}^F$. Let $\mathcal{G}^*$ be the dual group of
$\mathcal{G}$ and $F^\ast$ the dual Frobenius endomorphism and
denote $G^\ast={\mathcal{G}^\ast}^{F^\ast}$. Lusztig's
classification (see chapter 13 of \cite{DM}) says that the set of
irreducible complex characters of $G$ is partitioned into Lusztig
series $\mathcal{E}(G,(s))$ associated to various geometric
conjugacy classes $(s)$ of semi-simple elements of $G^\ast$. In
fact, $\mathcal{E}(G,(s))$ is the set of irreducible constituents of
a Deligne-Lusztig character $R_\mathcal{T}^\mathcal{G}(\theta)$,
where $(\mathcal{T},\theta)$ is of the geometric conjugacy class
associated to $(s)$. The elements of $\mathcal{E}(G,(1))$ are called
\emph{unipotent characters} of $G$. When $\mathcal{G}$ is a
connected reductive group, for any semi-simple element $s\in
G^\ast$, there is a bijection $\chi\mapsto\psi$ from
$\mathcal{E}(G,(s))$ to $\mathcal{E}(C_{G^\ast}(s),(1))$ such that
\begin{equation}\label{lusztig}
\chi(1)=\frac{|G|_{p\p}}{|C_{G^\ast}(s)|_{p\p}}\psi(1).
\end{equation}
In this situation, we will say that the irreducible complex
character $\chi$ of $G$ is parameterized by the pair $((s),\psi)$.

Suppose that we are determining irreducible complex characters of
$G$ of degrees up to a certain bound $D$. Unipotent characters of
$G$ as well as $G^\ast$ are classified by Lusztig (see \S13.8 of
\cite{C1}) and we will use that to find the low-dimensional
unipotent characters of $G$. When the character $\chi$ is not
unipotent, i.e. $(s)\neq (1)$, based on the structure of
$C_{G^\ast}(s)$, we estimate $(G^\ast:C_{G^\ast}(s))_{p\p}$ and come
up with certain cases when $C_{G^\ast}(s)$ is large enough. More
specifically, from formula \ref{lusztig}, we see that if $\chi(1)<D$
then $|C_{G^\ast}(s)|>|G|_{p\p}/D$. The following Proposition is
used frequently to determine unipotent characters of
$C_{G^\ast}(s)$.

\begin{proposition}[Proposition 13.20 of \cite{DM}]\label{DM}
Let $\mathcal{G}$ and $\mathcal{G}_1$ be two reductive groups
defined over $\mathbb{F}_q$, and let $f:\mathcal{G}\rightarrow
\mathcal{G}_1$ be a morphism of algebraic groups with a central
kernel, defined over $\mathbb{F}_q$ and such that $f(\mathcal{G})$
contains the derived group $\mathcal{G}_1\p$; then the unipotent
characters of $\mathcal{G}^F$ are the $\theta\circ f$, where
$\theta$ runs over the unipotent characters of $\mathcal{G}_1^F$.
\end{proposition}

\subsection{Centralizers of semi-simple elements} We collect here
some well-known results about the structures of centralizers of
semi-simple elements in finite classical groups (cf. \cite{C2},
\cite{FS}, and \cite{TZ2}). Since the following lemmas are similar,
we omit most of their proofs except the last one's.

\begin{lemma}\label{lemma1}
Let $q$ be odd and $s\in SO_{2n+1}(q)$ be semi-simple. Then
$$C_{SO_{2n+1}(q)}(s)\simeq SO_{2k+1}(q)\times O_{2(m-k)}^\pm(q)
\times\prod_{i=1}^tGL_{a_i}^{\alpha_i}(q^{k_i})$$ and
$$C_{O_{2n+1}(q)}(s)\simeq O_{2k+1}(q)\times O_{2(m-k)}^\pm(q)
\times\prod_{i=1}^tGL_{a_i}^{\alpha_i}(q^{k_i}),$$ where $0\leq
k\leq m\leq n$, $\alpha_i=\pm$, and $\sum_{i=1}^tk_ia_i=n-m$.
\end{lemma}

\begin{lemma}\label{lemma2}
Let $q$ be even and $s\in O^\pm_{2n}(q)$ be semi-simple. Then
$$C_{O^\pm_{2n}(q)}(s)\simeq O^\pm_{2m}(q)
\times\prod_{i=1}^tGL_{a_i}^{\alpha_i}(q^{k_i}),$$ where
$\alpha_i=\pm$ and $\sum_{i=1}^tk_ia_i=n-m$.
\end{lemma}

Let $(.,.)$ be a non-degenerate symplectic form on the space
$V=\mathbb{F}_q^{2n}$, then the \emph{conformal symplectic group}
$CSp_{2n}(q)$ is defined to be
$$\{g\in GL(V)\mid \exists \hspace{3pt}\tau(g)\in \FF_q^\ast, \forall u,v\in V,
(gu,gv)=\tau(g)(u,v)\}.$$

\begin{lemma}\label{lemma3}
Let $q$ be odd and $s\in CSp_{2n}(q)$ be semi-simple. Then
$$C_{CSp_{2n}(q)}(s)\simeq C_{Sp_{2n}(q)}(s)\cdot
\mathbb{Z}_{q-1}.$$ Moreover,
\begin{enumerate}
\item[(i)] If $\tau(s)$ is not a square in $\mathbb{F}_q$, then
$$C_{Sp_{2n}(q)}(s)\simeq Sp_m(q^2)\times\prod_{i=1}^tGL_{a_i}^{\alpha_i}(q^{k_i}),$$
where $m$ is even, $\alpha_i=\pm$, and $\sum_{i=1}^tk_ia_i=n-m$.
\item[(ii)] If $\tau(s)$ is a square in $\mathbb{F}_q$, then
$$C_{Sp_{2n}(q)}(s)\simeq Sp_{2k}(q)\times Sp_{2(m-k)}(q)
\times\prod_{i=1}^tGL_{a_i}^{\alpha_i}(q^{k_i}),$$ where $0\leq
k\leq m\leq n$, $\alpha_i=\pm$, and $\sum_{i=1}^tk_ia_i=n-m$.
\end{enumerate}
\end{lemma}

Let $Q(.)$ be a non-degenerate quadratic form on the space
$V=\mathbb{F}_q^{2n}$. Then the \emph{conformal orthogonal group}
$CO^\pm_{2n}(q)$ is defined to be $$\{g\in GL(V)\mid \exists
\hspace{3pt}\tau(g)\in \FF_q^\ast, \forall v\in V,
Q(g(v))=\tau(g)Q(v)\}.$$ We use the description of
$CO^\pm_{2n}(q)^0$ in \cite[Remark 7.3]{TZ1}. When $q$ is odd, it is
also shown in \cite[Lemma 7.4]{TZ1} that $CO^\pm_{2n}(q)^0=\{g\in
CO^\pm_{2n}(q)\mid \det(g)=\tau(g)^n\}$ and $CO^\pm_{2n}(q)^0$ is
actually a subgroup of index $2$ of $CO^\pm_{2n}(q)$ such that
$CO^\pm_{2n}(q)^0/SO^\pm_{2n}(q)\simeq
CO^\pm_{2n}(q)/O^\pm_{2n}(q)\simeq \mathbb{Z}_{q-1}$.

\begin{lemma}\label{lemma4}
Let $q$ be odd and $s\in CO^\pm_{2n}(q)^0$ be semi-simple. Then
$$C_{CO^\pm_{2n}(q)}(s)\simeq C_{O^\pm_{2n}(q)}(s)\cdot
\mathbb{Z}_{q-1},$$ and $$C_{CO^\pm_{2n}(q)^0}(s)\simeq
C_{SO^\pm_{2n}(q)}(s)\cdot \mathbb{Z}_{q-1}.$$ Moreover,
\begin{enumerate}
\item[(i)] If $\tau(s)$ is not a square in $\mathbb{F}_q$, then
$$C_{O^\pm_{2n}(q)}(s)\simeq O^\pm_m(q^2)\times\prod_{i=1}^tGL_{a_i}^{\alpha_i}(q^{k_i}),$$
where $m$ is even, $\alpha_i=\pm$, and $\sum_{i=1}^tk_ia_i=n-m$.
\item[(ii)] If $\tau(s)$ is a square in $\mathbb{F}_q$, then
$$C_{O^\pm_{2n}(q)}(s)\simeq O^\pm_{2k}(q)\times O^\pm_{2(m-k)}(q)
\times\prod_{i=1}^tGL_{a_i}^{\alpha_i}(q^{k_i}),$$ where $0\leq
k\leq m\leq n$, $\alpha_i=\pm$, and $\sum_{i=1}^tk_ia_i=n-m$.
\end{enumerate}
\end{lemma}

\begin{proof}
Let $(.,.)$ be the non-degenerate bilinear form associated with
$Q(.)$. Fix a basis of $V$ and let $J$ be the Gram matrix of $(.,.)$
corresponding to this basis. Then $\tn sJs=\tau J$, where
$\tau:=\tau(s)$. Hence $\Spec(s)=\Spec(\tn
s)=\tau\Spec(Js^{-1}J^{-1})=\tau\Spec(s^{-1})$. Denote the
characteristic polynomial of $s$ acting on $V$ by $P(x)\in
\FF_q[x]$. There are two cases:

\textbf{Case} $\mathbf{1}$: If $\tau$ is not a square in
$\mathbb{F}_q$, suppose $\tau=\lambda_1^2$ where $\lambda_1\in
\FF_{q^2}\setminus \FF_q$. Decompose $P(x)$ into irreducible
polynomials over $\mathbb{F}_q$:
$$P(x)=(x^2-\tau)^m\prod_{i=1}^l f_i^{m_i}(x)
\prod_{j=1}^{l'} g_j^{n_j}(x)\widehat{g_j}^{n_j}(x),$$ where
\begin{itemize}
  \item If $\lambda$ is a root of $f_i(x)$ then $\tau\lambda^{-1}$
  is also a root of $f_i(x)$, $\deg(f_i)$ is even, $\pm \lambda_1$ are not the roots of $f_i$,
  \item If $\lambda$ is a root of $g_j(x)$ then $\tau\lambda^{-1}$
  is a root of $\widehat{g_j}(x)$, $\deg(g_j)=\deg(\widehat{g_j})$,
  $g_j\neq \widehat{g_j}$,
  \item $n=m+\Sigma_{i=1}^l
  m_i\deg(f_i)/2+\Sigma_{j=1}^{l'}n_j\deg(g_j)$.
\end{itemize}
Since $s\in CO^\pm_{2n}(q)^0$, we have
$\tau^n=\det(s)=(-\tau)^m\tau^{n-m}$. Therefore, $m$ is even. Set
$$V_\tau:=(\frac{P(x)}{(x^2-\tau)^m})_{x=s}(V), V_i:=
(\frac{P(x)}{f^{m_i}_i(x)})_{x=s}(V),
U_j:=(\frac{P(x)}{g^{n_j}_j(x)\widehat{g}^{n_j}_j(x)})_{x=s}(V).$$
Then $V=V_\tau\oplus V_1\oplus\cdot\cdot\cdot V_l\oplus
U_1\oplus\cdot\cdot\cdot \oplus U_{l'}$ is an orthogonal
decomposition of $V$. Let $W\in
\{V,V_\tau,V_1,...,V_l,U_1,...,U_{l'}\}$. Then $s|_W$ is semi-simple
and $(.,.)$ is non-degenerate on $W$. For short notation, we denote
$C_{O(W)}(s|_W)$ and $C_{CO(W)}(s|_W)$ by $C_{O(W)}(s)$ and
$C_{CO(W)}(s)$, respectively.

1a) First, we show that $C_{O(W)}(s)\simeq GU_{m_0}(q^{k})$ and
$C_{CO(W)}(s)\simeq C_{O(W)}(s)\cdot \mathbb{Z}_{q-1}$, where
$W=V_i$ for $i=1,...,l$, $m_0=m_i$, and $k=\deg(f_i)/2$. Note that
the characteristic polynomial of $s$ acting on $W$ is
$f^{m_0}_i(x)$. Let $\lambda\in \mathbb{F}_{q^{2k}}$ be an
eigenvalue of the action of $s$ on $W$. Then all eigenvalues of the
action of $s$ on $W$ are $\lambda,\lambda^q,...,\lambda^{q^{2k-1}}$.
Since $\tau\lambda^{-1}$ is also a root of $f_i(x)$ and
$\pm\lambda_1$ are not, $\tau\lambda^{-1}=\lambda^{q^k}$. Consider
$\widetilde{W}=W\otimes_{\mathbb{F}_q}\mathbb{F}_{q^{2k}}$. Fix a
basis $(f_i)$ in $W$, and define a Frobenius endomorphism $\sigma:
\Sigma_ix_if_i\mapsto \Sigma_ix_i^qf_i$ on $\widetilde{W}$, where
$x_i\in \mathbb{F}_{q^{2k}}$. The simplicity of $s$ implies that
$\widetilde{W}=\widetilde{W}_1\oplus\cdot\cdot\cdot\oplus
\widetilde{W}_{2k}$, where
$\widetilde{W_i}=\Ker(s-\lambda^{q^{i-1}})$. We see that $\sigma$
permutes the $\widetilde{W}_i$'s cyclically:
$\sigma(\widetilde{W}_i)=\widetilde{W}_{i+1}$, where
$\widetilde{W}_{2k+1}=\widetilde{W}_1$. Let $g\in CO(W)$ be
commuting with $s|_{W}$. Then $g$ preserves each $\widetilde{W}_i$.
Moreover, it is easy to see that $g$ also commutes with $\sigma$.
This implies that the action of $g$ on $\widetilde{W}$ is completely
determined by its action on $\widetilde{W}_1$:
$g(\sigma^iw)=\sigma^i(gw)$ for $w\in \widetilde{W}_1$. So
$C_{CO(W)}(s)\hookrightarrow GL_{m_0}(q^{2k})$. If $u\in
\widetilde{W}_i$ and $v\in\widetilde{W}_j$ then
$\tau(u,v)=(su,sv)=\lambda^{q^{i-1}+q^{j-1}}(u,v)$. Therefore,
\begin{equation}\label{equation1}
\widetilde{W}^\perp_1=\bigoplus_{i\neq1+k}\widetilde{W}_i \text {
and } \widetilde{W}_{1+k}\cap \widetilde{W}^\perp_1=0.
\end{equation}
Choose a basis $(u_1,...,u_{m_0})$ in $\widetilde{W}_1$. Then
$(v_i)$ is a basis for $\widetilde{W}_{1+k}$, where
$v_i=\sigma^k(u_i)$. We see that $(\sigma u_i, \sigma
v_j)=(u_i,v_j)^q$. Hence, $(u_i,v_j)^{q^k}=(\sigma^ku_i,
\sigma^kv_j)=(v_i, u_j)=(u_j,v_i)$. In other words, $\tn U=U^{q^k}$.
Thus, together with \ref{equation1}, $U$ determines a non-degenerate
Hermitian form on an $m_0$-dimensional $\mathbb{F}_{q^{2k}}$-space.

If $g$ acts on $\widetilde{W}_1$ with matrix $A=(a_{ij})$ (with
respect to the basis $(u_i)$), then $g$ acts on
$\widetilde{W}_{1+k}$ with matrix $A^{q^k}=(a_{ij}^{q^k})$ (with
respect to the basis $(v_i)$). From \ref{equation1}, $g\in
C_{O(W)}(s)$ if and only if
$$(u_k,v_l)=(gu_k,gv_l)=(\sum_ia_{ik}u_i,\sum_ja_{jl}^{q^k}v_j)=
\sum_{i,j}a_{ik}a_{jl}^{q^k}(u_i,v_j),$$ i.e. $\tn AUA^{q^k}=U$.
Therefore, $C_{O(W)}(s)\simeq GU_{m_0}(q^k)$. Similarly, $g\in
C_{CO(W)}(s)$ if and only if $\tn AUA^{q^k}=\tau(g)U$. Note that if
$\tn AUA^{q^k}=U$ then $\tn (\varepsilon A)U(\varepsilon
A)^{q^k}=\tau(g)U$, where $\varepsilon$ is a scalar in
$\mathbb{F}_{q^{2k}}$ such that $\varepsilon^{q^k+1}=\tau(g)$ (there
always exists such an $\varepsilon$ for any $\tau(g)\in
\mathbb{F}_q^\ast$). That means $\tau: C_{CO(W)}(s)\rightarrow
\mathbb{F}_q^\ast$ is an epimorphism and hence $C_{CO(W)}(s)\simeq
C_{O(W)}(s)\cdot\mathbb{Z}_{q-1}$.

1b) Next, we show that $C_{O(W)}(s)\simeq GL_{n_0}(q^{k})$ and
$C_{CO(W)}(s)\simeq GL_{n_0}(q^k)\cdot \mathbb{Z}_{q-1}$, where
$W=U_j$ for $j=1,...,l'$, $n_0=n_j$, and $k=\deg(g_j)$. Note that
the characteristic polynomial of $s$ acting on $W$ is
$g^{n_0}_j(x)\widehat{g}^{n_0}_j(x)$. Let $\lambda\in
\mathbb{F}_{q^{k}}$ be a root of $g_j$. Then all the roots of $g_j$
are $\lambda,\lambda^q,...,\lambda^{q^{k-1}}$ and all the roots of
$\widehat{g}_j$ are $\tau\lambda^{-1},...,\tau\lambda^{-q^{k-1}}$.
Let $(f_i)$ be a basis of $W$, and define a Frobenius endomorphism
$\sigma: \Sigma_ix_if_i\mapsto \Sigma_ix_i^qf_i$ on
$\widetilde{W}:=W\otimes_{\mathbb{F}_q}\mathbb{F}_{q^k}$, where
$x_i\in \mathbb{F}_{q^{k}}$. The simplicity of $s$ implies that
$\widetilde{W}=\widetilde{W}_1\oplus\cdot\cdot\cdot\oplus
\widetilde{W}_{k}\oplus
\widetilde{W}_1'\oplus\cdot\cdot\cdot\oplus\widetilde{W}_k'$, where
$\widetilde{W}_i=\Ker(s-\lambda^{q^{i-1}})$ and
$\widetilde{W}_j'=\Ker(s-\tau\lambda^{-q^{j-1}})$. We see that
$\sigma$ permutes the $\widetilde{W}_i$'s and $\widetilde{W}_j'$'s
cyclically: $\sigma(\widetilde{W}_i)=\widetilde{W}_{i+1}$,
$\sigma(\widetilde{W}_j')=\widetilde{W}'_{j+1}$, where
$\widetilde{W}_{2k+1}=\widetilde{W}_1$ and
$\widetilde{W}'_{2k+1}=\widetilde{W}'_1$. Let $g\in CO(W)$ commuting
with $s|_{W}$. Then $g$ preserves each $\widetilde{W}_i$,
$\widetilde{W}'_j$. Again, the action of $g$ on $\widetilde{W}$ is
completely determined by its action on $\widetilde{W}_1$ and
$\widetilde{W}'_1$. So $C_{CO(W)}(s)\hookrightarrow
GL_{2n_0}(q^{k})$. We also have
\begin{equation}\label{equation2}
\widetilde{W}^\perp_1=\bigoplus_{i}\widetilde{W}_i \bigoplus _{j\geq
2}\widetilde{W}'_j\text { and }
\widetilde{W}'_{1}\cap\widetilde{W}^\perp_1=0.
\end{equation}
Choose basis $(u_i)$ and $(v_i)$ of $\widetilde{W}_1$ and
$\widetilde{W}'_1$, respectively. Suppose that $g$ acts on
$\widetilde{W}_1$ with matrix $A=(a_{ij})$ (with respect to the
basis $(u_i)$) and on $\widetilde{W}'_1$ with matrix $(b_{ij})$
(with respect to the basis $(v_i)$). From \ref{equation2}$, g\in
C_{O(W)}(s)$ if and only if $(gu,gv)=(u,v)$ for any $u\in
\widetilde{W}_1$, $v\in \widetilde{W}'_1$, i.e.
$$(u_k,v_l)=(gu_k,gv_l)=(\sum_ia_{ik}u_i,
\sum_jb_{jl}v_j)=\sum_{i,j}a_{ik}b_{jl}(u_i,v_j).$$ In other words,
$g\in C_{O(W)}(s)$ if and only if $\tn AUB=U$. Since $U$ is
non-degenerate, $\tn AUB=U$ means that $B$ is uniquely determined by
$A$. Therefore, $C_{O(W)}(s)\simeq GL_{n_0}(q^k)$. Similarly, $g\in
C_{CO(W)}(s)$ if and only if $(gu,gv)=\tau(g)(u,v)$ for any $u\in
\widetilde{W}_1$, $v\in \widetilde{W}'_1$ or equivalently $\tn
AUB=\tau(g)U$. Therefore, $C_{CO(W)}(s)\simeq C_{O(W)}(s)\cdot
\mathbb{Z}_{q-1}$.

1c) Lastly, we show that $C_{O(V_\tau)}(s)\simeq O^\pm_m(q^2)$ and
$C_{CO(V_\tau)}(s)\simeq C_{O(V_\tau)}(s)\cdot \mathbb{Z}_{q-1}$.
Again, we have $\widetilde{V}_\tau=\widetilde{V}_1\oplus
\widetilde{V}_2$, where
$\widetilde{V}_\tau=V_\tau\otimes_{\mathbb{F}_q}\mathbb{F}_{q^2}$,
$\widetilde{V}_1=\Ker(s-\lambda_1)$, and
$\widetilde{V}_2=\Ker(s+\lambda_1)$. Furthermore,
$\widetilde{V}_1\perp\widetilde{V}_2$ and therefore
$\widetilde{V}_1$, $\widetilde{V}_2$ are non-degenerate. Similar
arguments as in $1a)$, we see that the action of an element $g\in
C_{CO(V_\tau)}(s)$ on $\widetilde{V}_\tau$ is completely determined
by its action on $\widetilde{V}_1$. Hence, $C_{O(V_\tau)}(s)\simeq
O(\widetilde{V}_1)=O^\pm_m(q^2)$ and $C_{CO(V_\tau)}(s)\simeq
C_{O(V_\tau)}(s)\cdot \mathbb{Z}_{q-1}$, since $m$ is even.

Combining what we have proved in $1a), 1b), 1c)$, we have that
$$C_{O(V)}(s)\simeq C_{O(V_\tau)}(s)\times
\prod_iC_{O(V_i)}(s)\times\prod_jC_{O(U_j)}(s)\simeq$$$$\simeq
O^\pm_m(q^2)
\times\prod_iGU_{m_i}(q^{\deg(f_i)/2})\prod_jGL_{n_j}(q^{\deg(g_j)}).$$
We have also proved that $C_{CO(W)}(s)\simeq C_{O(W)}(s)\cdot
\mathbb{Z}_{q-1}$ for any $W\in\{V_\tau,V_1,...,V_l,U_1,$
$...,U_{l'}\}$. More precisely, given any $\alpha\in
\mathbb{F}_q^\ast$, there exists $g_W\in C_{CO(W)}(s)$ such that
$\tau(g_W)=\alpha$ for $W\in\{V_\tau,V_1,...,V_l,U_1,...,U_{l'}\}$.
Set $g=(g_{V_\tau}, g_{V_1},...,g_{V_l},g_{U_1},...,g_{U_{l'}})$.
Then $g\in C_{CO(V)}(s)$ and $\tau(g)=\alpha$. Hence,
$C_{CO(V)}(s)\simeq C_{O(V)}(s)\cdot\mathbb{Z}_{q-1}$. The Lemma is
proved in this case.

\textbf{Case} $\mathbf{2}$: If $\tau$ is a square in $\mathbb{F}_q$,
suppose that $\tau=\lambda_2^2$ with $\lambda_2\in \mathbb{F}_q$.
Replacing $s$ by $\lambda_2^{-1}s$ if necessary, we can assume that
$\tau=1$. Decompose $P(x)$ into irreducible polynomials over
$\mathbb{F}_q$: $$P(x)=(x-1)^{k'}(x+1)^{2m-k'}\prod_{i=1}^l
f_i^{m_i}(x) \prod_{j=1}^{l'} g_j^{n_j}(x)\widehat{g_j}^{n_j}(x),$$
where
\begin{itemize}
  \item If $\lambda$ is a root of $f_i(x)$ then $\lambda^{-1}$
  is also a root of $f_i(x)$, $\deg(f_i)$ is even, $\pm 1$ are not the roots of $f_i$,
  \item If $\lambda$ is a root of $g_j(x)$ then $\lambda^{-1}$
  is a root of $\widehat{g_j}(x)$,
  $\deg(g_j)=\deg(\widehat{g_j})$, $g_j\neq \widehat{g_j}$,
  \item $n=m+\Sigma_{i=1}^l
  m_i\deg(f_i)/2+\Sigma_{j=1}^{l'}n_j\deg(g_j)$.
\end{itemize}

Since $1=\tau^n=\det(s)=(-1)^{2m-k'}$, $k'$ is even. Set $k'=2k$ and
$$V':=\Ker(s-1), V'':=\Ker(s+1), V_i:=
(\frac{P(x)}{f^{m_i}_i(x)})_{x=s}(V),
U_j:=(\frac{P(x)}{g^{n_j}_j(x)\widehat{g}^{n_j}_j(x)})_{x=s}(V).$$
Then $V=V'\oplus V"\oplus V_1\oplus\cdot\cdot\cdot V_l\oplus
U_1\oplus\cdot\cdot\cdot \oplus U_{l'}$ is an orthogonal
decomposition of $V$. It is obvious that
$C_{O(V')}(s)=O(V')=O^\pm_{2k}(q)$ and
$C_{O(V")}(s)=O(V")=O^\pm_{2(m-k)}(q)$. Furthermore,
$C_{CO(V')}(s)=CO(V')=O^\pm_{2k}(q)\cdot \mathbb{Z}_{q-1}$ and
$C_{CO(V")}(s)=CO(V")=O^\pm_{2(m-k)}(q)\cdot \mathbb{Z}_{q-1}$. Now
we can repeat the arguments in Case 1 and complete the proof.
\end{proof}

\textbf{Remark}. Lemma \ref{lemma4} is wrong when $s\notin
CO_{2n}^\pm(q)^0$. Here is one example. Suppose that $q$ is odd. Let
$V=\mathbb{F}_q^2=\{(x,y)\mid x,y\in \mathbb{F}_q\}$ be the vector
space of dimension $2$. The quadratic form $Q(x,y)=xy$ is
non-degenerate of Witt index $1$. The Gram matrix of $Q$
corresponding to the basis $\{(1,0), (0,1)\}$ is
$J={0\hspace{6pt}1\choose 1\hspace{6pt}0}$. By definition,
$O_2^+(q)$ is the group of matrices $A\in GL_2(q)$ such that $\tn
AJA=J$.
Therefore, $$O_2^+(q)=\left\{\left(\hspace{-2mm}\begin{array}{cc} a&0\\
0& a^{-1}\end{array}\hspace{-2mm}\right), \left(\hspace{-2mm}\begin{array}{cc} 0&a\\
a^{-1}& 0\end{array}\hspace{-2mm}\right)\mid a\in
\mathbb{F}_q^\ast\right\}.$$ Similarly, $$CO_2^+(q)=\{A\in
GL_2(q)\mid \exists \tau\in \mathbb{F}_q^\ast, \tn AJA=\tau J\}=
\left\{\left(\hspace{-2mm}\begin{array}{cc} a&0\\
0& b\end{array}\hspace{-2mm}\right), \left(\hspace{-2mm}\begin{array}{cc} 0&a\\
b& 0\end{array}\hspace{-2mm}\right)\mid a,b\in
\mathbb{F}_q^\ast\right\}.$$ Take the element $s=J$, which is
semi-simple. Since $\tau(s)=1$ and $\det(s)=-1$, $s\notin
CO_2^+(q)^0$. Direct
computation shows that $$C_{O_2^+(q)}(s)=\left\{\left(\hspace{-2mm}\begin{array}{cc} 1&0\\
0& 1\end{array}\hspace{-2mm}\right), \left(\hspace{-2mm}\begin{array}{cc} -1&0\\
0&-1\end{array}\hspace{-2mm}\right), \left(\hspace{-2mm}\begin{array}{cc} 0&1\\
1& 0\end{array}\hspace{-2mm}\right), \left(\hspace{-2mm}\begin{array}{cc} 0&-1\\
-1& 0\end{array}\hspace{-2mm}\right)\right\}.$$ On the other hand,
$$C_{CO_2^+(q)}(s)=
\left\{\left(\hspace{-2mm}\begin{array}{cc} a&0\\
0& a\end{array}\hspace{-2mm}\right), \left(\hspace{-2mm}\begin{array}{cc} 0&a\\
a& 0\end{array}\hspace{-2mm}\right)\mid a\in
\mathbb{F}_q^\ast\right\}.$$ Therefore, $C_{CO_2^+(q)}(s)\ncong
C_{O_2^+(q)}(s)\cdot \mathbb{Z}_{q-1}$ by order comparison.

%%% -----------------------------------------------------------------------------------

\section{Unipotent characters}
\subsection{Unipotent Characters of $SO_{2n+1}(q)$ and $PCSp_{2n}(q)$}
Proposition $5.1$ of \cite{TZ1} shows that $SO_{2n+1}(q)$ as well as
$PCSp_{2n}(q)$ have a unique unipotent character of minimal degree
$(q^n-1)(q^n-q)/2(q+1)$ and any other non-trivial unipotent
character has degree greater than $q^{2n}/2(q+1)$. We mimic its
proof and get Proposition \ref{unipotentpro}, which classifies
unipotent characters of degrees up to $\approx {q^{6n-15}}/{2}$.

From \cite[p. 466, 467]{C1} , we know that the unipotent characters
of $SO_{2n+1}(q)$ and $PCSp_{2n}(q)$ are labeled by symbols of the
form
$${\lambda \choose \mu}={\lambda_1\hspace{6pt}\lambda_2\hspace{6pt}\lambda_3\hspace{6pt}
...\hspace{6pt}\lambda_a\choose\mu_1\hspace{6pt}\mu_2\hspace{6pt}...\hspace{6pt}\mu_b}$$
where $0\leq\lambda_1<\lambda_2<... < \lambda_a$,
$0\leq\mu_1<\mu_2<...<\mu_b$, $a-b$ is odd and positive,
$(\lambda_1, \mu_1)\neq (0,0)$, and
$$\sum_i\lambda_i+\sum_j\mu_j-\left(\frac{a+b-1}{2}\right)^2=n.$$
The integer $n$ is called the rank of the symbol $\lambda
\choose\mu$. The degree of the unipotent character
$\chi^{\lambda,\mu}$ corresponding to the symbol $\lambda\choose\mu$
is equal to
$$\frac{(q^2-1)(q^4-1)\cdot\cdot\cdot(q^{2n}-1)\prod_{i\p<i}(q^{\lambda_i}-q^{\lambda_{i\p}})
\prod_{j\p<j}(q^{\mu_j}-q^{\mu_{j\p}})\prod_{i,j}(q^{\lambda_i}+q^{\mu_j})}
{2^{\frac{a+b-1}{2}}q^{{a+b-2\choose2}+
{a+b-4\choose2}+\cdot\cdot\cdot}\prod_i\prod_{k=1}^{\lambda_i}(q^{2k}-1)
\prod_j\prod_{k=1}^{\mu_j}(q^{2k}-1)}.$$

\begin{proposition}\label{unipotentpro}
Let $G^\ast$ be either $(B_n)_{ad}(q)=SO_{2n+1}(q)$ or
$(C_n)_{ad}(q)=PCSp_{2n}(q)$. Suppose that $n\geq6$ and $\chi\in
\Irr(G^\ast)$ is unipotent. Then either $\chi$ is one of characters
labeled by ${n\choose-}$, $0\hspace{3pt} 1 \hspace{3pt}n\choose-$,
$0\hspace{3pt} 1 \choose n$, $1\hspace{3pt} n \choose0$,
$0\hspace{3pt} n \choose1$, $0\hspace{3pt} 2
 \hspace{3pt}n-1\choose-$,
$0\hspace{3pt} 2 \choose n-1$, $2\hspace{3pt} n-1 \choose0$,
$0\hspace{3pt} n-1 \choose2$, $1\hspace{3pt} n-1 \choose1$,
$0\hspace{3pt} 1 \hspace{3pt} 2 \hspace{3pt}n\choose1$,
$0\hspace{3pt} 1 \hspace{3pt}2\choose1\hspace{3pt}n$, $1\hspace{3pt}
2 \hspace{3pt}n\choose0 \hspace{3pt}1$, $0\hspace{3pt} 1
\hspace{3pt}n\choose1\hspace{3pt} 2$, or $\chi(1)\geq
(q^{2n-2}-1)(q^{2n}-1)(q^{n-2}-1) (q^{n-2}-q^3)/2$
$(q^2-1)(q^4-1)(q^3+1)$ $(\approx {q^{6n-15}}/{2}),$ which is the
degree of the unipotent character labeled by $0\hspace{3pt} 3
\hspace{3pt}n-2\choose-$.
\end{proposition}

\textbf{Note}: For reader's convenience, we put these
low-dimensional unipotent characters in Table 1, where the degree of
each character is calculated. To keep the paper not too lengthy, we
will skip detailed arguments of some inequalities in this and
ongoing proofs. We frequently use the following obvious estimates
without notice:

(i) $q^{a-b}\leq\dfrac{q^a-1}{q^b-1}\leq
q^{a-b}(1+\dfrac{1}{q^b-1})$ for $a\geq b\geq 1$.

(ii)
$\prod_{i=2}^n\dfrac{q^{a_i}-q^{a_1}}{q^{a_i}-q^{a_1-1}}>\dfrac{1}{2}$
for integers $1\leq a_1<a_2<...<a_n$.

\begin{proof} $1)$ Define
$$D(n)=(q^{2n-2}-1)(q^{2n}-1)(q^{n-2}-1)(q^{n-2}-q^3)/2(q^2-1)(q^4-1)(q^3+1).$$
We have that
$$\frac{D(n)}{D(n-1)}=\frac{(q^{2n}-1)(q^{n-2}-1)(q^{n-5}-1)}
{(q^{2n-4}-1)(q^{n-3}-1)(q^{n-6}-1)}$$ and $q^6<D(n)/D(n-1)<q^7$ for
every $n\geq 7$. Also, $D(n)<q^{6n-14}/2$ for every $n\geq 6$. Let
$\mathcal{L}_n$ be the set of $14$ symbols in this Proposition. Note
that if $(n,q)\neq (6,2)$ then the degrees of characters
corresponding to the symbols in $\mathcal{L}_n$ is smaller than
$D(n)$. We will prove by induction on $n\geq 6$ that if
$\chi=\chi^{\lambda,\mu}$ with $\lambda\choose\mu$ is not in
$\mathcal{L}_n$, then $\chi(1)\geq D(n)$. The induction base $n=6$
can be checked easily by using \cite[Table 27]{CH}. The rest of the
proof establishes the induction step for $n\geq7$, by means of
induction on $b\geq0$.

$2)$ Now we consider the case $b=0$. Then $a\geq 1$ is odd. We may
assume $a\geq 3$, and $\lambda\neq (0, 1,n)$, $(0,2,n-1)$,
$(0,3,n-2)$. First we assume $a=3$. If $\lambda_1=0$ then
$\lambda=(0,k,n+1-k)$, where $4\leq k<(n+1)/2$; in particular,
$n\geq8$. In this case,
$$\chi(1)=\frac{1}{2}\cdot \frac{(q^{2(k+1)}-1)\cdot\cdot\cdot
(q^{2n}-1)}{(q^4-1)\cdot\cdot\cdot
(q^{2(n-k+1)}-1)}\cdot\frac{(q^k-1)(q^{n-k+1}-1)(q^{n-k+1}-q^k)}{q^2-1}>$$
$$>\frac{1}{2}q^{2(k-1)(n-k)}\cdot\frac{(q^k-1)(q^{n-k+1}-1)(q^{n-k+1}-q^k)}{q^2-1}.$$
Since $3\leq k-1<n-k$, $2(k-1)(n-k)\geq6(n-4)=6n-24$. Moreover,
$n-k+1>k$. Therefore
$(q^k-1)(q^{n-k+1}-1)(q^{n-k+1}-q^k)/(q^2-1)>q^{k-2}q^{2k}\geq
q^{10}$. It follows that
$\chi(1)>\frac{1}{2}q^{6n-24}q^{10}=q^{6n-14}/2>D(n)$ as required.
If $\lambda_1\geq1$, then $\lambda=(k,l,n+1-k-l)$ with $1\leq
k<l<n+1-k-l$. Then
$$\chi(1)=\frac{1}{2}\cdot\frac{(q^{2(k+1)}-1)\cdot\cdot\cdot
(q^{2(k+l)}-1)}{(q^2-1)\cdot\cdot\cdot (q^{2l}-1)}\cdot
\frac{(q^{2(k+l+1)}-1)\cdot\cdot\cdot
(q^{2n}-1)}{(q^4-1)\cdot\cdot\cdot (q^{2(n+1-k-l)}-1)}\cdot$$
$$\cdot \frac{(q^l-q^k)(q^{n+1-k-l}-q^k)(q^{n+1-k-l}-q^l)}{q^2-1}>
\frac{1}{2}q^{2kl}q^{2(k+l-1)(n-k-l)}q^{2k+l}.$$ If $(k,l)=(1,2)$
then $\lambda=(1,2,n-2)$ and
$$\chi(1)=\frac{(q^2-q)(q^{2n-2}-1)(q^{2n}-1)(q^{n-2}-q)(q^{n-2}-q^2)}{2(q^2-1)^2(q^4-1)}>D(n)$$
for every $n\geq 6$. It remains to consider $(k,l)\neq(1,2)$. Then
$k+l\geq 4$ and $l\geq 3$. It follows that $3\leq k+l-1$ and $3\leq
l\leq n-k-l$. Therefore, $2(k+l-1)(n-k-l)\geq 6(n-4)$. So, in this
case,
$\chi(1)>\frac{1}{2}q^{2kl}q^{6n-24}q^{2k+l}>\frac{1}{2}q^{6n-14}>D(n)$
and we are done. Now we may assume $a\geq 5$.

Suppose that $\lambda_1\geq 1$. Consider the unipotent character
$\chi\p$ labeled by the symbol $\lambda\p\choose-$ of rank $n-1$,
where $\lambda\p=(\lambda_1-1, \lambda_2,...,\lambda_a)$. We have
$$\frac{\chi(1)}{\chi\p(1)}=\frac{(q^{\lambda_2}-q^{\lambda_1})
\cdot\cdot\cdot(q^{\lambda_a}-q^{\lambda_1})}{(q^{\lambda_2}-q^{\lambda_1-1})
\cdot\cdot\cdot(q^{\lambda_a}-q^{\lambda_1-1})}\cdot\frac{q^{2n}-1}{q^{2\lambda_1}-1}
\geq \frac{q^{2(n-\lambda_1)}}{2}.$$ Note that
$\sum_{i=1}^a\lambda_i-\left(\frac{a-1}{2}\right)^2=n$ and hence
$n-\lambda_1=\lambda_2+\cdot\cdot\cdot+\lambda_a-\left(\frac{a-1}{2}\right)^2
\geq(a-1)\lambda_1+\frac{a^2-1}{4}\geq10$. Therefore,
$\chi(1)/\chi\p(1)\geq q^{20}/2>q^7$. It follows that $\chi(1)\geq
q^7D(n-1)>D(n)$. Now we may assume that $\lambda_1=0$.

Next we assume that $\lambda_i-\lambda_{i-1}\geq 2$ for some
$i\geq2$. Consider the unipotent character $\chi\p$ labeled by the
symbol $\lambda\p\choose-$ of rank $n-1$, where
$\lambda\p=(\lambda_1,...,\lambda_{i-1},\lambda_i-1,\lambda_{i+1},...,\lambda_a)$.
Similarly, we have $\chi(1)/\chi\p(1)>q^{2(n-\lambda_i)}/2$. By
induction hypotheses, $\chi\p(1)\geq D(n-1)$. If $n-\lambda_i\geq 4$
then $\chi(1)\geq q^7D(n-1)>D(n)$. Now suppose $n-\lambda_i\leq 3$.
Note that $\sum_i\lambda_i=\left(\frac{a-1}{2}\right)^2+n$. If
$i\leq a-1$ then $\lambda_{a-1}\geq n-3$ and $\lambda_{a}\geq n-2$.
Then $\left(\frac{a-1}{2}\right)^2+n\geq
0+1+\cdot\cdot\cdot+(a-3)+(n-3)+(n-2)$ or $a^2-8a-9+4n\leq 0$. This
is a contradiction because $a\geq 5$ and $n\geq 7$. Therefore $i=a$
and we have
$0+1+\cdot\cdot\cdot+(a-2)+\lambda_a=\left(\frac{a-1}{2}\right)^2+n$.
Since $\lambda_a\geq (n-3)$, $a^2-4a-9\leq 0$ and therefore $a=5$.
Then $\lambda=(0,1,2,3,n-2)$ and it is easy to check that
$\chi(1)>D(n)$. Now we may assume $\lambda=(0,1,2,...,a-2,a-1)$.
Consider $\chi\p=\chi^{\lambda\p,\mu\p}$, where
$\lambda\p=(0,1,2,...,a-4,a-3)$, and $\mu\p=\mu$. Note that
$a^2=4n+1$ and the rank of $\chi\p$ is $n-a+1$. We have
$$\frac{\chi(1)}{\chi\p(1)}=\frac{(q^{2(n-a+2)}-1)...(q^{2n}-1)\prod_{i\geq a-2, i>i\p}
(q^i-q^{i\p})}{2q^{{a-2\choose2}}
\prod_{k=1}^{a-2}(q^{2k}-1)\prod_{k=1}^{a-1}(q^{2k}-1)}=$$
$$=\frac{(q^{2(n-a+2)}-1)...(q^{2n}-1)}{(q^2-1)...(q^{2(a-1)}-1)}
\cdot \frac{q^{a-1}-q^{a-2}}{2} \cdot
\prod_{i=0}^{a-3}\frac{(q^{a-1}-q^i)(q^{a-2}-q^i)}{q^i(q^{2(i+1)}-1)}.$$
Remark that
$$\prod_{i=0}^{a-3}\frac{(q^{a-1}-q^i)(q^{a-2}-q^i)}{q^i(q^{2(i+1)}-1)}=
\prod_{i=0}^{a-3}\frac{q^i(q^{i+2}-1)}{q^{i+1}+1}\geq\prod_{i=0}^{2}
\frac{q^i(q^{i+2}-1)}{q^{i+1}+1}>q^4.$$ Therefore,
$$\frac{\chi(1)}{\chi\p(1)}\geq \frac{q^{2(n-a+1)(a-1)}(q^{a-1}-q^{a-2})q^4}{2}.$$
Since $n\geq 7$ and $a^2=4n+1$, $n\geq 12$ and $a\geq 7$. Then
$(n-a+1)(a-1)\geq 6(n-6)$. Hence,
$\chi(1)/\chi\p(1)>q^{12n-63}/2>q^7$ and we are done.

\medskip

$3)$ From now on we assume that $b\geq 1$. At this point we suppose
that $(\lambda_1,\mu_1)\neq (1,0)$ and $\lambda_1\geq 1$. Consider
the character $\chi\p$ labeled by the symbol $\lambda\p\choose\mu\p$
of rank $n-1$, where
$\lambda\p=(\lambda_1-1,\lambda_2,...,\lambda_a)$ and $\mu\p=\mu$.
If ${\lambda\p\choose\mu\p}\in \mathcal{L}_{n-1}$ then
${\lambda\choose\mu}\in \left\{{1\hspace{3pt}2\choose n-2}, {3
\hspace{3pt}n-2\choose0}, {1\hspace{3pt}n-2\choose2}, {2
\hspace{3pt}n-2\choose1}\right\}$. It is easy to check that the
degrees of characters corresponding to these symbols are greater
than $D(n)$. Now we can assume $\lambda\p\choose\mu\p$ is not in
$\mathcal{L}_{n-1}$. We have
$$\frac{\chi(1)}{\chi\p(1)}=\frac{q^{2n}-1}{q^{2\lambda_1}-1}\cdot
\prod_{i=2}^a\frac{q^{\lambda_i}-q^{\lambda_1}}{q^{\lambda_i}-q^{\lambda_1-1}}
\cdot
\prod_{j=1}^b\frac{q^{\mu_j}+q^{\lambda_1}}{q^{\mu_j}+q^{\lambda_1-1}}>
\frac{q^{2n}-1}{2(q^{2\lambda_1}-1)}>\frac{q^{2(n-\lambda_1)}}{2}.$$
If $n-\lambda_1\leq3$ then $\lambda_1\geq n-3$ and we have
$n-3+n-2+...+n+a-4+\frac{(b-1)b}{2}\leq\frac{(a+b-1)^2}{4}+n$. This
implies $4a(n-3)+(a-b)^2\leq 4n+1$ and therefore $a(n-3)\leq n$.
Since $n\geq 7$, $a=1$ which leads to a contradiction since
$b\geq1$. So $n-\lambda_1\geq 4$. Then we have
$\chi(1)/\chi\p(1)>q^7$ and $\chi(1)>q^7D(n-1)>D(n)$ as desired.

Similarly, for $(\lambda_1,\mu_1)\neq (0,1)$ and $\mu_1\geq 1$, we
consider the character $\chi\p$ labeled by the symbol
$\lambda\p\choose\mu\p$ of rank $n-1$, where $\lambda\p=\lambda$ and
$\mu\p=(\mu_1-1,\mu_2,...,\mu_b)$. If ${\lambda\p\choose\mu\p}\in
\mathcal{L}_{n-1}$ then ${\lambda\choose\mu}\in \left\{{2
\hspace{3pt}n-2\choose1}, {0 \hspace{3pt}n-2\choose3}, {1
\hspace{3pt}n-2\choose2}, {0
\hspace{3pt}1\hspace{3pt}2\hspace{3pt}n-1\choose2}, {0
\hspace{3pt}1\hspace{3pt}2\choose2\hspace{3pt}n-1}\right\}$. Again,
it is easy to check that the degrees of characters corresponding to
these symbols are greater than $D(n)$. Now we can assume
$\lambda\p\choose\mu\p$ is not in $\mathcal{L}_{n-1}$. We have
$$\frac{\chi(1)}{\chi\p(1)}=\frac{q^{2n}-1}{q^{2\mu_1}-1}\cdot
\prod_{j=2}^b\frac{q^{\mu_j}-q^{\mu_1}}{q^{\mu_j}-q^{\mu_1-1}} \cdot
\prod_{i=1}^a\frac{q^{\lambda_i}+q^{\mu_1}}{q^{\lambda_i}+q^{\mu_1-1}}>
\frac{q^{2n}-1}{2(q^{2\mu_1}-1)}>\frac{q^{2(n-\mu_1)}}{2}.$$ If
$n-\mu_1\geq 4$ then $\chi(1)/\chi\p(1)>q^7$ and
$\chi(1)>q^7D(n-1)>D(n)$ as desired. If $n-\mu_1\leq3$ then
$\mu_1\geq n-3$ and we have
$\frac{(a-1)a}{2}+n-3+n-2+...+n+b-4\leq\frac{(a+b-1)^2}{4}+n$. This
implies $4b(n-3)+(a-b)^2\leq 4n+1$. So $b(n-3)\leq n$. Since $n\geq
7$, $b=1$. Then we have $(a-1)^2\leq 13$ and therefore $a=2$ or
$a=4$.

\textbf{Case} $\mathbf{a=4,b=1:}$ Then
$\lambda_1+\lambda_2+\lambda_3+\lambda_4+\mu_1=n+4$. This implies
$\mu_1\leq n-2$. Hence ${\lambda\choose\mu} ={0
\hspace{3pt}1\hspace{3pt}2\hspace{3pt}3\choose n-2}$ or $0
\hspace{3pt}1\hspace{3pt}2\hspace{3pt}4\choose n-3$. The degrees of
characters corresponding to both these symbols are greater than
$D(n)$.

\textbf{Case} $\mathbf{a=2,b=1:}$ Since $n-3\leq\mu_1\leq n$ and
$\lambda\choose\mu$ is not in $\mathcal{L}_{n}$, we have
${\lambda\choose\mu} ={0\hspace{3pt}4\choose n-3}$, $1
\hspace{3pt}3\choose n-3$, $0 \hspace{3pt}3\choose n-2$ or $1
\hspace{3pt}2\choose n-2$. Again, the degrees of characters
corresponding to these symbols are greater than $D(n)$.

It remains to consider the case where $(\lambda_1, \mu_1)=(0,1)$ or
$(1,0)$.

\medskip

$4)$ Here we suppose that $(\lambda_1,\mu_1)=(1,0)$ but $\lambda\neq
(1,2,...,a)$. Then there exists an $i\geq2$ such that
$\lambda_i\geq\lambda_{i-1}+2$. We choose $i$ to be smallest
possible. If $a=2$ then $b=1$ and ${\lambda\choose\mu}
={1\hspace{3pt}n\choose0}\in \mathcal{L}_n$. Hence $a\geq 3$.
Consider the character $\chi\p$ labeled by the symbol
 $\lambda\p\choose\mu\p$ of
rank $n-1$, where
$\lambda\p=(\lambda_1,...,\lambda_{i-1},\lambda_i-1,\lambda_{i+1},...,\lambda_a)$
and $\mu\p=\mu$. If ${\lambda\p\choose\mu\p}\in \mathcal{L}_{n-1}$
then $\lambda\choose\mu$ can only be $1
\hspace{3pt}3\hspace{3pt}n-1\choose0\hspace{3pt}1$. The degree of
the corresponding character is greater than $D(n)$. Now we can
assume that $\lambda\p\choose\mu\p$ is not in $\mathcal{L}_{n-1}$.
By induction hypothesis, $\chi\p(1)\geq D_{n-1}$. Set
$$T_1=\prod_{i\p<i}\frac{q^{\lambda_i}-q^{\lambda_{i\p}}}{q^{\lambda_i-1}-q^{\lambda_{i\p}}},
 T_2=\prod_{i\p>i}\frac{q^{\lambda_{i\p}}-q^{\lambda_{i}}}{q^{\lambda_{i\p}}-q^{\lambda_{i}-1}},
 T_3=\prod_j\frac{q^{\lambda_i}+q^{\mu_j}}{q^{\lambda_i-1}+q^{\mu_j}}.$$
Then
$\chi(1)/\chi\p(1)=[(q^{2n}-1)/(q^{2\lambda_i}-1)]T_1T_2T_3>q^{2(n-\lambda_i+1)}/2$
(see \cite[p. 2121]{TZ1}). If $n-\lambda_i\geq 3$ then
$\chi(1)/\chi\p(1)>q^7$ and therefore $\chi(1)\geq D(n)$ as desired.
Now we assume $n-\lambda_i\leq 2$. There are two cases:

\textbf{Case} $\mathbf{n-\lambda_i=2:}$ Then $i\geq a-1$. If $i=a-1$
then $a-b=1$, $a=n-2$ and
${\lambda\choose\mu}={1\hspace{3pt}\cdot\cdot\cdot\hspace{3pt}n-4\hspace{3pt}n-2\hspace{3pt}
n-1\choose0\hspace{3pt}\cdot\cdot\cdot\hspace{3pt}n-4}$. Since
$n\geq 7$, $a\geq5$ and $i=b\geq 4$. It follows that $T_1\geq q^3$,
$T_3\geq q^{b-1}\geq q^3$. Hence $\chi(1)/\chi\p(1)>q^{10}/2>q^7$
and we are done. If $i=a\geq3$ then $T_1\geq
q^{a-2}\frac{q^{\lambda_i}-q}{q^{\lambda_i-1}-q}$, $T_2=1$ and
$T_3\geq \frac{q^{\lambda_a}+1}{q^{\lambda_a-1}+1}$. Then
$\chi(1)/\chi\p(1)>q^7$ and we are done.

\textbf{Case} $\mathbf{n-\lambda_i\leq 1:}$ Then $i=a$ and
$\lambda_a\geq n-1$. If $\lambda_a\geq n$ then $\lambda_a=n$,
$a-b=1$ and
${\lambda\choose\mu}={1\hspace{3pt}2\hspace{3pt}...\hspace{3pt}a-1\hspace{3pt}n
\choose0\hspace{3pt}1\hspace{3pt}...\hspace{3pt}a-2}$. Since
$\lambda\choose\mu$ is not in $\mathcal{L}_n$, $a\geq 4$. We have
$$\frac{\chi(1)}{\chi\p(1)}=\prod_{j=1}^{a-1}\left(\frac{q^n-q^j}{q^{n-1}-q^j}
\cdot\frac{q^n+q^{j-1}}{q^{n-1}+q^{j-1}}\right)=\frac{(q^n+1)(q^{n-1}-1)}
{(q^{n-a}-1)(q^{n-a+1}+1)}>q^{2a-2}.$$ If $a\geq5$ then
$\chi(1)/\chi\p(1)>q^8$ and therefore we are done. If $a=4$ then
${\lambda\choose\mu}={1\hspace{3pt}2\hspace{3pt}3\hspace{3pt}n\choose
0\hspace{3pt}1\hspace{3pt}2}$ and $\chi(1)>D(n)$. If $\lambda_a=n-1$
then we also have $a-b=1$ and
${\lambda\choose\mu}={1\hspace{3pt}2\hspace{3pt}...\hspace{3pt}a-1\hspace{3pt}n-1
\choose0\hspace{3pt}1\hspace{3pt}...\hspace{3pt}a-3\hspace{3pt}a-1}$.
If $a=3$ then
${\lambda\choose\mu}={1\hspace{3pt}2\hspace{3pt}n-1\choose0\hspace{3pt}2}$.
Check directly we see that $\chi(1)>D(n)$. If $a\geq 4$ then $T_1>
q^{a-1}\geq q^3$, $T_2=1$ and $T_3\geq q^{b-1}\geq q^2$. So
$\chi(1)/\chi\p(1)>q^7$ and we are done again.

Similarly, suppose that $(\lambda_1,\mu_1)=(0,1)$ but $\mu\neq
(1,2,...,b)$; in particular, $b\geq2$. Then there exists an index
$j\geq2$ such that $\mu_j\geq\mu_{j-1}+2$. We choose $j$ to be
smallest possible. Consider the character $\chi\p$ labeled by the
symbol
 $\lambda\p\choose\mu\p$ of rank $n-1$, where $\lambda\p=\lambda$ and
$\mu\p=(\mu_1,...,\mu_{j-1},\mu_j-1,\mu_{j+1},...,\mu_b)$. If
${\lambda\p\choose\mu\p}\in \mathcal{L}_{n-1}$ then
$\lambda\choose\mu$ can only be
$0\hspace{3pt}1\hspace{3pt}n-1\choose1\hspace{3pt}3$. The degree of
the corresponding character is greater than $D(n)$ and hence we are
done in this case. Now we can assume $\lambda\p\choose\mu\p$ is not
in $\mathcal{L}_{n-1}$. Set
$$U_1=\prod_{j\p<j}\frac{q^{\mu_j}-q^{\mu_{j\p}}}{q^{\mu_j-1}-q^{\mu_{j\p}}},
 U_2=\prod_{j\p>j}\frac{q^{\mu_{j\p}}-q^{\mu_{j}}}{q^{\mu_{j\p}}-q^{\mu_{j}-1}},
 U_3=\prod_i\frac{q^{\lambda_i}+q^{\mu_j}}{q^{\lambda_i}+q^{\mu_j-1}}.$$
Then
$\chi(1)/\chi\p(1)=[(q^{2n}-1)/(q^{2\mu_j}-1)]U_1U_2U_3>q^{2(n-\mu_j+1)}/2$
(see \cite[p. 2122]{TZ1}). If $n-\mu_j\geq 3$ then
$\chi(1)/\chi\p(1)>q^7$ and therefore $\chi(1)\geq D_n$ as desired.
Now we assume $n-\mu_j\leq 2$. There are two cases:

\textbf{Case} $\mathbf{n-\mu_j=2:}$ Then $j\geq b-1$. If $j=b-1$
then $a-b=1$, $b=n-2$ and
${\lambda\choose\mu}={0\hspace{3pt}1\hspace{3pt}...\hspace{3pt}n-2
\choose1\hspace{3pt}...\hspace{3pt}n-4\hspace{3pt}n-2\hspace{3pt}n-1}$.
Since $n\geq 7$, $b\geq5$ and $j\geq 4$. It follows that $U_1\geq
q^2\frac{q^{\mu_j}-q}{q^{\mu_j-1}-q}$, $T_3\geq
\frac{q^{\mu_j}+1}{q^{\mu_j-1}+1}$. Hence
$\chi(1)/\chi\p(1)>q^8/2>q^7$ and we are done. If $j=b\geq3$ then
$V_1\geq q\frac{q^{\mu_j}-q}{q^{\mu_j-1}-q}$, $V_2=1$ and $V_3\geq
\frac{q^{\lambda_a}+1}{q^{\lambda_a-1}+1}$. Therefore,
$\chi(1)/\chi\p(1)>q^7$ and we are done again. If $j=b=2$ then $a=3$
or $a=5$. In any case, $\lambda_a\leq4<n-2$. Therefore, in this
case, $U_1>q$, $U_2=1$ and $U_3>q^2$. Then $\chi(1)/\chi\p(1)>q^7$
and $\chi(1)\geq D_n$ as required.

\textbf{Case} $\mathbf{n-\mu_j=1:}$ Then $j=b$ and $\mu_b\geq n-1$.
If $\mu_b\geq n$ then $\mu_b=n$, $a-b=1$ and ${\lambda\choose\mu}
={0\hspace{3pt}1\hspace{3pt}...\hspace{3pt}b
\choose1\hspace{3pt}2\hspace{3pt}...\hspace{3pt}b-1\hspace{3pt}n}$.
Since $\lambda\choose\mu$ is not in $\mathcal{L}_n$, $b\geq 3$. We
have
$$\frac{\chi(1)}{\chi\p(1)}=\prod_{i=0}^{b}\frac{q^n+q^i}{q^{n-1}+q^i}
\cdot\prod_{i=1}^{b-1}\frac{q^n-q^{i}}{q^{n-1}-q^{i}}=\frac{(q^n+1)(q^{n-1}-1)}
{(q^{n-b-1}+1)(q^{n-b}-1)}>q^{2b-1}.$$ If $b\geq4$ then
$\chi(1)/\chi\p(1)>q^7$ and therefore we are done. If $b=3$ then
${\lambda\choose\mu}={0\hspace{3pt}1\hspace{3pt}2\hspace{3pt}3\choose
1\hspace{3pt}2\hspace{3pt}n}$. In this case, $\chi(1)>D(n)$. If
$\mu_b=n-1$ then we also have $a-b=1$ and ${\lambda\choose\mu}
={0\hspace{3pt}1\hspace{3pt}...\hspace{3pt}b-1\hspace{3pt}b+1
\choose1\hspace{3pt}2\hspace{3pt}...\hspace{3pt}b-1\hspace{3pt}n-1}$
where $n\geq b+2$. If $b=2$ then ${\lambda\choose\mu}
={0\hspace{3pt}1\hspace{3pt}3\choose1\hspace{3pt}n-1}$. Check
directly we see that $\chi(1)>D(n)$. If $b\geq 3$ then $U_1>
q^{b-1}\geq q^2$, $U_2=1$ and $U_3\geq q^3$. So
$\chi(1)/\chi\p(1)>q^7$ and we are done.

\medskip

$5)$ Here we suppose that $\mu_1=0$ and $\lambda=(1,2,...,a)$.
Consider the character $\chi\p$ labeled by the symbol
$${\lambda\p\choose\mu\p}=
{0\hspace{6pt}1\hspace{6pt}2\hspace{6pt}
...\hspace{6pt}a\choose\mu_2\hspace{6pt}\mu_3\hspace{6pt}
...\hspace{6pt}\mu_b}$$ of the same rank $n$, but with the
parameters $a\p=a+1$ and $b\p=b-1$. It is easy to see that
$\lambda\p\choose\mu\p$ is not in $\mathcal{L}_{n}$. So
$\chi\p(1)>D(n)$ by the induction hypothesis (on $b$). But
$$\frac{\chi(1)}{\chi\p(1)}=\frac{\prod_{i=1}^a(1+\frac{2}{q^i-1})}
{\prod_{j=2}^b(1+\frac{2}{q^{\mu_j}-1})}\geq
\frac{\prod_{i=1}^a(1+\frac{2}{q^i-1})}
{\prod_{j=1}^{b-1}(1+\frac{2}{q^{j}-1})}>1,$$ so we are done.

Similarly, suppose that $\lambda_1=0$ and $\mu=(1,2,...,b)$ but
$\lambda\neq (0,1,...,a-1)$. Then there exists an $i\geq2$ such that
$\lambda_i\geq\lambda_{i-1}+2$. We choose $i$ to be smallest
possible. Consider the character $\chi\p$ labeled by the symbol
$${\lambda\p\choose\mu\p}={\lambda_1\hspace{6pt}...\hspace{6pt}
\lambda_{i-1}\hspace{6pt}\lambda_i-1\hspace{6pt}\lambda_{i+1}
\hspace{6pt}...\hspace{6pt}\lambda_a\choose\mu_1\hspace{6pt}...\hspace{6pt}\mu_b}$$
of rank $n\p=n-1$. If ${\lambda\p\choose\mu\p}\in\mathcal{L}_{n-1}$,
then $\lambda\choose\mu$ is one of the symbols $0
\hspace{3pt}1\hspace{3pt}3\hspace{3pt} n-1\choose1$, $0
\hspace{3pt}2\hspace{3pt}n-1\choose1\hspace{3pt}2$. In this case, we
can check that $\chi(1)>D(n)$. Now we can suppose that
$\lambda\p\choose\mu\p$ is not in $\mathcal{L}_{n-1}$. We have
$\frac{\chi(1)}{\chi\p(1)}=\frac{q^{2n}-1}{q^{2\lambda_i}-1}\cdot
V_1V_2V_3$, where
$$V_1=\prod_{i\p<i}\frac{q^{\lambda_i}-q^{\lambda_{i\p}}}{q^{\lambda_i-1}
-q^{\lambda_{i\p}}},
V_2=\prod_{i\p>i}\frac{q^{\lambda_{i\p}}-q^{\lambda_{i}}}{q^{\lambda_{i\p}}
-q^{\lambda_{{i}}-1}},
V_3=\prod_{j}\frac{q^{\lambda_i}+q^{\mu_j}}{q^
{\lambda_i-1}+q^{\mu_j}}.$$ Note that $V_1\geq
(q^{\lambda_i}-1)/(q^{\lambda_i-1}-1)>q$, $V_2>1/2$ and $V_3>1$. So
$\chi(1)/\chi\p(1)> q^{2(n-\lambda_i)+1}/2$. If $n-\lambda_i\geq4$,
then $\chi(1)/\chi\p(1)>q^9/2$ and we are done. Now we assume that
$n-\lambda_i\leq3$. Then $i\geq a-1$. There are two cases:

\textbf{Case} $\mathbf{i=a-1:}$ We have
$(a-b)^2+4b+4\lambda_{a-1}+4\lambda_a=4n+8a-11$. Since
$(a-b)^2+4b\geq 4a-3$ and $\lambda_a\geq a$, $\lambda_{a-1}\leq
n-2$. If $\lambda_{a-1}=n-2$ then $\lambda_a=a$. It follows that
$a\geq n-1\geq 6$. Hence, $V_1\geq q^{a-2}\geq q^4$ and we have
$$\frac{\chi(1)}{\chi\p(1)}>\frac{1}{2}q^4\frac{q^{2n}-1}{q^{2n-4}-1}\geq q^7,$$
so $\chi(1)>D_n$. If $\lambda_{a-1}=n-3$ then $\lambda_a\leq a+1$.
Then $a+1\geq n-2$ and therefore $a\geq 4$. Hence $V_1\geq
q^{a-2}\geq q^2$. Again, we have
$$\frac{\chi(1)}{\chi\p(1)}>\frac{1}{2}q^2\frac{q^{2n}-1}{q^{2n-6}-1}\geq
q^7$$ and we are done.

\textbf{Case} $\mathbf{i=a:}$ Then $(a-b)^2+4b+4\lambda_a=4n+4a-3$.
Therefore, $(a-b-2)^2=4(n-\lambda_a)+1$. Since $n-\lambda_a\leq 3$,
$n-\lambda_a=0$ or $2$. If $n-\lambda_a=2$ then $(a-b-2)^2=9$ and we
have $a-b=5$. In particular, $a\geq 6$. Then $V_1> q^{a-1}\geq q^5$.
Moreover, in this case, $V_2=1$ and $V_3> 1$. So
$\chi(1)/\chi\p(1)>q^9$ and we are done. If $n-\lambda_a=0$ then
$a-b=1$ or $3$. In this case, $V_1>q^{a-1}$, $V_2=1$ and
$V_3>q^{b-1}$. So $\chi(1)/\chi\p(1)>q^{a+b-2}$. If $(a,b)=(4,3)$ or
$(5,2)$ then
${\lambda\choose\mu}={0\hspace{3pt}1\hspace{3pt}2\hspace{3pt}n
\choose1\hspace{3pt}2\hspace{3pt}3}$ or
$0\hspace{3pt}1\hspace{3pt}2\hspace{3pt}3\hspace{3pt}
n\choose1\hspace{3pt}2\hspace{3pt}3$. Checking directly we have
$\chi(1)>D(n)$. Otherwise, $a+b\geq 9$ and therefore
$\chi(1)/\chi\p(1)\geq q^7$ and we are done again.

\medskip

$6)$ Finally, we consider the case where ${\lambda\choose\mu}=
{0\hspace{3pt}1\hspace{3pt}2\hspace{3pt}...\hspace{3pt}
a-1\choose1\hspace{3pt}2\hspace{3pt}...\hspace{3pt} b}$. Consider
the character $\chi\p$ labeled by the symbol
${\lambda\p\choose\mu\p}=
{0\hspace{3pt}1\hspace{3pt}2\hspace{3pt}...\hspace{3pt}
a-2\choose1\hspace{3pt}2\hspace{3pt}...\hspace{3pt}b-1}$ of rank
$n-1$. Then we have
$$\frac{\chi(1)}{\chi\p(1)}=\frac{q^{a+b-2}(q^{2n}-1)}{(q^b-1)(q^{a-1}+1)}>
\frac{q^{2n}-1}{q}>q^7.$$ Note that $\lambda\p\choose\mu\p$ is not
in $\mathcal{L}_{n-1}$. Therefore, $\chi(1)>q^7\chi\p(1)\geq
q^7D(n-1)>D(n)$ as desired.
\end{proof}

\begin{corollary}\label{corollary}
Let $G^\ast$ be either $(B_n)_{ad}(q)=SO_{2n+1}(q)$ or
$(C_n)_{ad}(q)=PCSp_{2n}(q)$. Suppose that $n\geq5$ and $\chi\in
\Irr(G^\ast)$ is unipotent. Then either $\chi$ is one of characters
labeled by $n\choose-$, $0\hspace{3pt} 1 \hspace{3pt}n\choose-$,
$0\hspace{3pt} 1 \choose n$, $1\hspace{3pt} n \choose0$,
$0\hspace{3pt} n \choose1$, with degrees $1$,
$(q^n-1)(q^n-q)/2(q+1)$, $(q^n+1)(q^n+q)/2(q+1)$,
$(q^n+1)(q^n-q)/2(q-1)$, $(q^n-1)(q^n+q)/2(q-1)$, respectively, or
$\chi(1)>q^{4n-8}$.
\end{corollary}

\begin{proof}
If $n\geq6$, Corollary comes from Proposition \ref{unipotentpro}.
The case $n=5$ can be verified easily by using Table $26$ of
\cite{CH}.
\end{proof}

%%% -----------------------------------------------------------------------------------

\subsection{Unipotent Characters of $P(CO_{2n}^-(q)^0)$}
Proposition $7.1$ in \cite{TZ1} shows that the projective conformal
orthogonal group of type $-$, $P(CO^-_{2n}(q)^0)$, has a unique
unipotent character of minimal degree $(q^n+1)(q^{n-1}-q)/(q^2-1)$
and any other non-trivial unipotent character has degree greater
than $q^{2n-2}$. We mimic its proof and get Proposition \ref{prop1},
which classifies unipotent characters of degrees up to $q^{4n-10}$.

From \cite[p. 475, 476]{C1}, we know that the unipotent characters
of $P(CO_{2n}^-(q)^0)$ are labeled by symbols of the form
$${\lambda\choose\mu}=
{\lambda_1\hspace{6pt}\lambda_2\hspace{6pt}\lambda_3\hspace{6pt} ...
\hspace{6pt}\lambda_a\choose\mu_1\hspace{6pt}
\mu_2\hspace{6pt}...\hspace{6pt} \mu_b},$$ where
$0\leq\lambda_1<\lambda_2<... < \lambda_a$,
$0\leq\mu_1<\mu_2<...<\mu_b$, $a>b$, $a-b\equiv2(\bmod$ $4)$,
$(\lambda_1, \mu_1)\neq (0,0)$, and
$$\sum_i\lambda_i+\sum_j\mu_j-\left[\left(\frac{a+b-1}{2}\right)^2\right]=n.$$
The integer $n$ is called the rank of the symbol
$\lambda\choose\mu$. The degree of the unipotent character
$\chi^{\lambda,\mu}$ corresponding to the symbol $\lambda\choose\mu$
is equal to
$$\frac{(q^n+1)\prod_{i=1}^{n-1}(q^{2i}-1)\prod_{i\p<i}(q^{\lambda_i}-q^{\lambda_{i\p}})
\prod_{j\p<j}(q^{\mu_j}-q^{\mu_{j\p}})\prod_{i,j}(q^{\lambda_i}+q^{\mu_j})}
{2^{\frac{a+b-2}{2}}q^{{a+b-2\choose2}+{a+b-4\choose2}
+\cdot\cdot\cdot}\prod_i\prod_{k=1}^{\lambda_i}(q^{2k}-1)
\prod_j\prod_{k=1}^{\mu_j}(q^{2k}-1)}.$$

\begin{proposition}\label{prop1}
Let $G^\ast=(^2D_n)_{ad}(q)=P(CO_{2n}^-(q)^0)$. Suppose that
$n\geq6$ and $\chi\in \Irr(G^\ast)$ is unipotent. Then either $\chi$
is one of the characters labeled by $0\hspace{3pt} n\choose-$, $1
\hspace{3pt} n-1\choose-$, $0 \hspace{3pt} 1 \hspace{3pt} n\choose1$
with degrees $1$, $(q^n+1)(q^{n-1}-q)/(q^2-1)$,
$(q^{2n}-q^2)/(q^2-1)$, respectively, or $\chi(1)> q^{4n-10}$.
Furthermore, when $n=5$, $G^\ast$ has one more character of degree
$q^2(q^4+1)(q^5+1)/(q+1)$ (corresponding to the symbol $2
\hspace{3pt} 3\choose-$), which is less than $q^{4n-10}$.
\end{proposition}

\begin{proof}
Denote $\mathcal{L}_n={0 \hspace{3pt} n\choose-}$, $1 \hspace{3pt}
n-1\choose-$, $0\hspace{3pt}1\hspace{3pt} n\choose1$ if $n\geq 6$
and $\mathcal{L}_n={0\hspace{3pt} n\choose-}$, $1 \hspace{3pt}
n-1\choose-$, $0 \hspace{3pt} 1 \hspace{3pt} n\choose1$, $2
\hspace{3pt} n-2\choose-$ if $n=5$. We will prove by induction on
$n\geq 5$ that if $\chi=\chi^{\lambda,\mu}$ and
${\lambda\choose\mu}$ $\notin\mathcal{L}_n$, then
$\chi(1)>q^{4n-10}$. The induction base $n=5$ can be checked easily
by using \cite[Table 31]{CH}. The rest of the proof establishes the
induction step for $n\geq 6$.

\medskip

$1)$ At this point we suppose that $(\lambda_1,\mu_1)\neq (1,0)$ and
$\lambda_1\geq 1$ (eventually $b$ may be zero). Consider the
unipotent character $\chi\p$ labeled by the symbol
$\lambda\p\choose\mu\p$ of rank $n-1$, where
$\lambda\p=(\lambda_1-1, \lambda_2,...,\lambda_a)$ and $\mu\p=\mu$.
If ${\lambda\p\choose\mu\p}\in \mathcal{L}_{n-1}$ then
${\lambda\choose\mu}={2 \hspace{3pt} n-2\choose-}$ and therefore
$\chi(1)=(q^{2n-2}-1)(q^n+1)(q^{n-2}-q^2)/(q^2-1)(q^4-1)>q^{4n-10}$.
So we can suppose that $\lambda\p\choose\mu\p$ is not in
$\mathcal{L}_{n-1}$. By induction hypothesis,
$\chi\p(1)>q^{4(n-1)-10}$. Observe that
$n=\sum_i\lambda_i+\sum_j\mu_j-\frac{(a+b)^2-2(a+b)}{4}\geq
a\lambda_1+\frac{(a-b)^2}{4}$. Since $a\geq2$ and $a-b\geq2$,
$2\lambda_1<n$. Note that $n\geq6$. It follows that $\lambda_1\leq
n-4$. We have
$$\frac{\chi(1)}{\chi\p(1)}=\frac{(q^n+1)(q^{n-1}-1)}{q^{2\lambda_1}-1}\cdot
\prod_{i=2}^a\frac{q^{\lambda_i}-q^{\lambda_1}}{q^{\lambda_i}-q^{\lambda_1-1}}
\cdot\prod_{j=1}^b\frac{q^{\lambda_1}+q^{\mu_j}}{q^{\lambda_1-1}+q^{\mu_j}}>$$
$$>\frac{(q^n+1)(q^{n-1}-1)}{2(q^{2\lambda_1}-1)}
\geq\frac{(q^n+1)(q^{n-1}-1)}{2(q^{2n-8}-1)}>q^4,$$ and therefore
$\chi(1)>q^{4n-10}$.

\medskip

$2)$ Similarly, supposing $(\lambda_1,\mu_1)\neq (0,1)$ and
$\mu_1\geq 1$, consider the unipotent character $\chi\p$ labeled by
the symbol $\lambda\p\choose\mu\p$ of rank $n-1$, where
$\lambda\p=\lambda$ and $\mu\p=(\mu_1-1, \mu_2,..., \lambda_b)$. If
${\lambda\p\choose\mu\p}\in \mathcal{L}_{n-1}$ then
${\lambda\choose\mu}={0 \hspace{3pt} 1\hspace{3pt} n-1\choose2}$ and
therefore
$\chi(1)=(q^{n}+1)(q^{n-1}-1)(q^{n-1}-q)(q^{n-1}+q^2)/2(q^2-1)^2>q^{4n-10}$.
So we can suppose that $\lambda\p\choose\mu\p$ is not in
$\mathcal{L}_{n-1}$. By induction hypothesis,
$\chi\p(1)>q^{4(n-1)-10}$. Set
$$T_1=\prod_{j=2}^{b}\frac{q^{\mu_j}-q^{\mu_1}}{q^{\mu_j}-q^{\mu_1-1}},
T_2=\prod_{i=1}^{a}\frac{q^{\lambda_i}+q^{\mu_1}}{q^{\lambda_i}+q^{\mu_1-1}}.$$
Then, if $n-\mu_1\geq 4$,
$$\frac{\chi(1)}{\chi\p(1)}=\frac{(q^n+1)(q^{n-1}-1)}{q^{2\mu_1}-1}T_1T_2>
\frac{(q^n+1)(q^{n-1}-1)}{2(q^{2\mu_1}-1)}>q^4.$$ If $n-\mu_1=3$
then $b=1$ since $b\mu_1<n$ and $n\geq6$. Then $T_1=1$ and we still
have
$\frac{\chi(1)}{\chi\p(1)}>\frac{(q^n+1)(q^{n-1}-1)}{q^{2\mu_1}-1}>q^4$.
It follows that $\chi(1)>q^{4n-10}$ when $n-\mu_1\geq3$. If
$n-\mu_1\leq 2$, then ${\lambda\choose\mu}={0 \hspace{3pt}
1\hspace{3pt} 2\choose n-1}$, $0 \hspace{3pt} 1\hspace{3pt} 3\choose
n-2$. Checking directly, we again have $\chi^{\lambda,
\mu}(1)>q^{4n-10}$.

\medskip

$3)$ Now we consider the case where $b=0$ and $\lambda_1=0$. If
$a=2$ then we can suppose $\lambda=(k,n-k)$ with $3\leq k<n-k$
($k=2$ is considered already in $1)$). Then
$$\chi(1)=\frac{(q^n+1)(q^{n-k}-q^k)\prod_{i=k+1}^{n-1}(q^{2i}-1)}
{\prod_{i=1}^{n-k}(q^{2i}-1)}>$$$$>q^{n+k-2}q^{2(k-1)(n-k-1)}\geq
q^{n+1}q^{4n-16}>q^{4n-10}.$$ Now we may assume $a\geq6$. First we
suppose that $\lambda\neq (0, 1, ..., a-1)$. Then there exists
$i\geq 2$ such that $\lambda_i-\lambda_{i-1}\geq 2$. Consider the
character $\chi\p$ corresponding to the symbol $\lambda\p\choose-$,
where $\lambda\p=(..., \lambda_{i-1}, \lambda_i-1, \lambda_{i+1},
...)$. Note that $n-\lambda_i\geq(a-2)^2/4\geq4$. Therefore,
$\chi(1)/\chi\p(1)> (q^n+1)(q^{n-1}-1)/2(q^{2\lambda_i}-1)>q^4$. The
induction hypothesis $\chi\p(1)>q^{4n-14}$ implies that
$\chi(1)>q^{4n-10}$.

Next we suppose that $\lambda=(0, 1, ..., a-1)$. Then
$\chi(1)=(q^n+1)f(a)$ and $n=a^2/4$, where
\begin{equation}\label{functionf}f(a)=\frac{(q^2-1)(q^4-1)\cdot\cdot\cdot(q^{(a^2-4)/2}-1)
\prod_{0\leq i\p\leq i\leq
a-1}(q^i-q^{i\p})}{2^{\frac{a-2}{2}}q^{{a-2\choose2}
+{a-4\choose2}+\cdot\cdot\cdot}
\prod_{i=0}^{a-1}\prod_{k=1}^i(q^{2k}-1)}\end{equation} for any even
$a$. It is easy to show that
$f(a)/f(a-2)>q^{a(a-1)(a-3)/2}/2^6>q^{3a^2/4}=q^{3n}$ for $a\geq6$.
Furthermore, $f(4)>1$. It follows that $\chi(1)>q^{4n-10}$ and we
are done.

It remains to consider the case where $b\geq 1$ and
$(\lambda_1,\mu_1)=(0,1)$ or $(1,0)$.

\medskip

$4)$ Here we suppose that $(\lambda_1,\mu_1)=(1,0)$ but $\lambda\neq
(1,2,...,a)$. Then there exists $i\geq2$ such that
$\lambda_i\geq\lambda_{i-1}+2$. Choose $i$ to be smallest possible.
Consider the unipotent character $\chi\p$ labeled by
$\lambda\p\choose\mu\p$ of rank $n-1$, with
$\lambda\p=(\lambda_1,...,\lambda_{i-1},\lambda_i-1,\lambda_{i+1},...,\lambda_a)$
and $\mu\p=\mu$. By the induction hypothesis, $\chi\p(1)>q^{4n-14}$.
Set
$$U_1=\prod_{i\p<i}\frac{q^{\lambda_i}-q^{\lambda_{i\p}}}{q^{\lambda_i-1}-q^{\lambda_{i\p}}},
U_2=\prod_{i\p>i}\frac{q^{\lambda_{i\p}}-q^{\lambda_{i}}}{q^{\lambda_{i\p}}-q^{\lambda_i-1}},
U_3=\prod_{j}\frac{q^{\lambda_i}+q^{\mu_{j}}}{q^{\lambda_i-1}+q^{\mu_{j}}}.$$
If $n-\lambda_i\geq3$ then
$$\frac{\chi(1)}{\chi\p(1)}=\frac{(q^n+1)(q^{n-1}-1)}{q^{2\lambda_i}-1}U_1U_2U_3>
\frac{(q^n+1)(q^{n-1}-1)}{q^{2\lambda_i}-1}\cdot\frac{q^{\lambda_i}-q}{q^{\lambda_i-1}-q}\cdot\frac{1}{2}
\cdot\frac{q^{\lambda_i}+1}{q^{\lambda_i-1}+1}>$$$$>\frac{(q^n+1)(q^{n-1}-1)q^2}{2(q^{2\lambda_i}-1)}
\geq\frac{(q^n+1)(q^{n-1}-1)q^2}{2(q^{2n-6}-1)}>q^4,$$ so
$\chi(1)>q^{4n-10}$. If $n-\lambda_i\leq2$ then $i=a$. Note that
$n-\lambda_i\geq1$. As $b\geq1$, we have $a\geq3$. Therefore
$$U_1\geq
\frac{(q^{\lambda_a}-q)(q^{\lambda_a}-q^2)}{(q^{\lambda_a-1}-q)(q^{\lambda_a-1}-q^2)},
U_2=1, U_3\geq \frac{q^{\lambda_a}+1}{q^{\lambda_a-1}+1}.$$ It
follows that
$$\frac{\chi(1)}{\chi\p(1)}\geq\frac{q^n+1}{q^{n-1}+1}\cdot
\frac{(q^{\lambda_a}-q)(q^{\lambda_a}-q^2)}{(q^{\lambda_a-1}-q)(q^{\lambda_a-1}-q^2)}\cdot
\frac{q^{\lambda_a}+1}{q^{\lambda_a-1}+1}>q^4,$$ so
$\chi(1)>q^{4n-10}$ as required.

\medskip

$5)$ Similarly, suppose that $(\lambda_1,\mu_1)=(0,1)$ but $\mu\neq
(1,2,...,b)$. Then there exists $j\geq2$ such that
$\mu_j\geq\mu_{j-1}+2$. Consider the unipotent character $\chi\p$
labeled by $\lambda\p\choose\mu\p$ of rank $n-1$, with
$\lambda\p=\lambda$ and
$\mu\p=(\mu_1,...,\mu_{j-1},\mu_j-1,\mu_{j+1},...,\mu_b)$. By the
induction hypothesis, $\chi\p(1)>q^{4n-14}$. Set
$$V_1=\prod_{j\p<j}\frac{q^{\mu_j}-q^{\mu_{j\p}}}{q^{\mu_j-1}-q^{\mu_{j\p}}},
V_2=\prod_{j\p>j}\frac{q^{\mu_{j\p}}-q^{\mu_{j}}}{q^{\mu_{j\p}}-q^{\mu_j-1}},
V_3=\prod_{i}\frac{q^{\lambda_i}+q^{\mu_{j}}}{q^{\lambda_i}+q^{\mu_{j}-1}}.$$
If $n-\mu_j\geq3$ then
$$\frac{\chi(1)}{\chi\p(1)}=\frac{(q^n+1)(q^{n-1}-1)}{q^{2\mu_j}-1}V_1V_2V_3>
\frac{(q^n+1)(q^{n-1}-1)}{q^{2\mu_j}-1}\cdot\frac{q^{\mu_j}-q}{q^{\mu_j-1}-q}\cdot\frac{1}{2}
\cdot\frac{q^{\mu_j}+1}{q^{\mu_j-1}+1}>$$$$>\frac{(q^n+1)(q^{n-1}-1)q^2}{2(q^{2\mu_j}-1)}>q^4,$$
so $\chi(1)>q^{4n-10}$ as desired. If $n-\mu_j=2$ then $j=b$ and
therefore $V_2=1$. We have
$\chi(1)/\chi\p(1)>(q^n+1)(q^{n-1}-1)q^2/(q^{2\mu_j}-1)>q^4$ and we
are done. Now suppose that $n-\mu_j\leq 1$. Then $n-\mu_j=1$,
$a-b=2$, and
$${\lambda\choose\mu}=
{0\hspace{6pt}1\hspace{6pt}...
\hspace{6pt}a-2\hspace{6pt}a-1\choose1\hspace{6pt}
2\hspace{6pt}...\hspace{6pt} b-1\hspace{6pt}n-1}.$$ As $b\geq2$, we
have $a\geq4$. It follows that $V_1>q$,
$V_3\geq(q^{n-1}+1)/(q^{n-5}+1)>q^3$. Hence $\chi(1)/\chi\p(1)>q^4$
and we are done again.

\medskip

$6)$ Here we suppose that $\mu_1=0$ and $\lambda=(1,2,...,a)$. First
we consider the case where $\mu\neq (0,1,...,b-1)$. Then there
exists an index $j\geq2$ such that $\mu_j\geq\mu_{j-1}+2$. Consider
the unipotent character $\chi\p$ labeled by $\lambda\p\choose\mu\p$
of rank $n-1$, with $\lambda\p=\lambda$ and
$\mu\p=(\mu_1,...,\mu_{j-1},\mu_j-1,\mu_{j+1},...,\mu_b)$. Since
$\lambda\p\choose\mu\p$ is not in $\mathcal{L}_{n-1}$, by the
induction hypothesis, we have $\chi\p(1)>q^{4n-14}$. Observe that
$n-\mu_j\geq (a-b+2)^2/4\geq4$. So we have
$$\frac{\chi(1)}{\chi\p(1)}=\frac{(q^n+1)(q^{n-1}-1)}{q^{2\mu_j}-1}
\cdot\prod_{j\p<j}\frac{q^{\mu_j}-q^{\mu_{j\p}}}{q^{\mu_j-1}-q^{\mu_{j\p}}}
\cdot\prod_{j\p>j}\frac{q^{\mu_{j\p}}-q^{\mu_{j}}}{q^{\mu_{j\p}}-q^{\mu_j-1}}
\cdot\prod_{i}\frac{q^{\lambda_i}+q^{\mu_{j}}}{q^{\lambda_i}+q^{\mu_{j}-1}}>$$
$$>\frac{(q^n+1)(q^{n-1}-1)}{2(q^{2\mu_j}-1)}>q^4,$$ and hence
$\chi(1)>q^{4n-10}$.

Next we consider the case where $\mu=(0,1,...,b-1)$. Observe that
$n=a+(a-b)^2/4$. If $b=1$ then $n=(a+1)^2/4$ and
$\chi(1)>(q^n+1)f(a+1)$ where $f$ is the function defined in formula
(\ref{functionf}). Since $n\geq6$, it follows that $a\geq7$. We have
already proved that $f(a+1)>q^{3n-10}$ for $a\geq5$. Therefore
$\chi(1)>q^{4n-10}$ as required. Now we can assume $b\geq2$.
Consider the unipotent character $\chi\p$ labeled by
$\lambda\p\choose\mu\p$ of rank $n-1$, with
$\lambda\p=(1,2,...,a-1)$ and $\mu\p=(0,1,...,b-2)$. Since
$\lambda\p\choose\mu\p$ is not in $\mathcal{L}_{n-1}$, again by the
induction hypothesis, $\chi\p(1)>q^{4n-14}$. Furthermore,
$$\frac{\chi(1)}{\chi\p(1)}=\frac{q^{a+b-2}(q^n+1)(q^{n-1}-1)}{(q^a-1)(q^{b-1}+1)}>q^{2n-3}>q^4.$$
Consequently, $\chi(1)>q^{4n-10}$.

\medskip

$7)$ Finally, we suppose that $\lambda_1=0$ and $\mu=(1,2,...,b)$.
First we consider the case where $\lambda\neq (0,1,...,a-1)$. Then
there exists $i\geq2$ such that $\lambda_i\geq\lambda_{i-1}+2$.
Choose $i$ to be smallest possible. Consider the unipotent character
$\chi\p$ labeled by $\lambda\p\choose\mu\p$ of rank $n-1$, with
$\lambda\p=(\lambda_1,...,\lambda_{i-1},\lambda_i-1,\lambda_{i+1},...,\lambda_a)$
and $\mu\p=\mu$. If ${\lambda\p\choose\mu\p}\in \mathcal{L}_{n-1}$
then ${\lambda\choose\mu}={0\hspace{3pt} 2\hspace{3pt} n-1\choose1}$
and by checking directly we have $\chi(1)>q^{4n-10}$. Therefore we
can suppose that $\lambda\p\choose\mu\p$ is not in
$\mathcal{L}_{n-1}$. In other words, by the induction hypothesis,
$\chi\p(1)>q^{4n-14}$. Remark that $n-\lambda_i\geq0$. Set
$$W_1=\prod_{i\p<i}\frac{q^{\lambda_i}-q^{\lambda_{i\p}}}{q^{\lambda_i-1}-q^{\lambda_{i\p}}},
W_2=\prod_{i\p>i}\frac{q^{\lambda_{i\p}}-q^{\lambda_{i}}}{q^{\lambda_{i\p}}-q^{\lambda_i-1}},
W_3=\prod_{j}\frac{q^{\lambda_i}+q^{\mu_{j}}}{q^{\lambda_i-1}+q^{\mu_{j}}}.$$
Note that $W_1\geq (q^{\lambda_i}-1)/(q^{\lambda_i-1}-1)>q$,
$W_2>1/2$ and $W_3>1$. Therefore, if $n-\lambda_i\geq3$,
$$\frac{\chi(1)}{\chi\p(1)}=\frac{(q^n+1)(q^{n-1}-1)}{q^{2\lambda_i}-1}W_1W_2W_3>
\frac{(q^n+1)(q^{n-1}-1)q}{2(q^{2\lambda_i}-1)}>q^4,$$ so
$\chi(1)>q^{4n-10}$. If $1\leq n-\lambda_i\leq2$ then either $i=a$
or
$${\lambda\choose\mu}={0\hspace{6pt}...\hspace{6pt}n-4\hspace{6pt}n-2\hspace{6pt}n-1\choose
1\hspace{6pt}...\hspace{6pt}n-3}.$$ In the former case, $W_2=1$ and
we have
$$\frac{\chi(1)}{\chi\p(1)}\geq\frac{(q^n+1)(q^{n-1}-1)}{q^{2\lambda_i}-1}\cdot
\frac{(q^{\lambda_i}-1)(q^{\lambda_i}-q)}{(q^{\lambda_i-1}-1)(q^{\lambda_i-1}-q)}
\cdot\frac{q^{\lambda_i}+q}{q^{\lambda_i-1}+q}>$$$$>\frac{(q^n+1)(q^{n-1}-1)}
{(q^{n-1}+1)(q^{n-2}-1)}\cdot\frac{q^{2\lambda_i}-q^2}{q^{2\lambda_i-2}-q^2}>q^4,$$
and therefore we are done. In the latter case, since $n\geq6$,
$W_1\geq q^4$ and we are done also. If $n=\lambda_i$, then $i=a$,
$a-b=2$ and $\lambda=(0,1,...,a-2,n)$. Since $\lambda\choose\mu$ is
not in $\mathcal{L}_n$, $a\geq 4$. Then we have
$$\frac{\chi(1)}{\chi\p(1)}=\prod_{i=1}^{a-2}\frac{q^n-q^i}{q^{n-1}-q^i}\cdot
\prod_{j=1}^{a-2}\frac{q^n+q^j}{q^{n-1}+q^j}>q^{2(a-2)}\geq q^4,$$
so $\chi(1)>q^{4n-10}$.

The last configuration we have to handle is that
$\lambda=(0,1,...,a-1)$ and $\mu=(1,2,...,b)$. Consider the
unipotent character $\chi\p$ labeled by $\lambda\p\choose\mu\p$ of
rank $n-1$, with $\lambda\p=(0,1,...,a-2)$ and $\mu\p=(1,2,...,b-1)$
(if $b=1$, $\mu\p$ is just empty). By the induction hypothesis,
$\chi\p(1)>q^{4n-14}$. Furthermore,
$$\frac{\chi(1)}{\chi\p(1)}=\frac{q^{a+b-2}(q^n+1)(q^{n-1}-1)}{(q^a-1)(q^{b-1}+1)}>q^{2n-3}>q^4.$$
Consequently, $\chi(1)>q^{4n-10}$.
\end{proof}

%========================================================================

\subsection{Unipotent Characters of $P(CO_{2n}^+(q)^0)$}
Proposition $7.2$ in \cite{TZ1} shows that the projective conformal
orthogonal group of type $+$, $P(CO^+_{2n}(q)^0)$, has a unique
unipotent character of minimal degree $(q^n-1)(q^{n-1}+q)/(q^2-1)$
and any other non-trivial unipotent character has degree greater
than $q^{2n-2}$. We mimic its proof and get Proposition \ref{prop2},
which classifies unipotent characters of degrees up to $q^{4n-10}$.

From \cite[p. 471, 472]{C1}, we know that the unipotent characters
of $G$ are labeled by symbols of the form
$${\lambda\choose\mu}=
{\lambda_1\hspace{6pt}\lambda_2\hspace{6pt}\lambda_3\hspace{6pt} ...
\hspace{6pt}\lambda_a\choose\mu_1\hspace{6pt}
\mu_2\hspace{6pt}...\hspace{6pt} \mu_b},$$ where
$0\leq\lambda_1<\lambda_2<... < \lambda_a$,
$0\leq\mu_1<\mu_2<...<\mu_b$, $a-b\equiv0(\bmod$ $4)$, $(\lambda_1,
\mu_1)\neq (0,0)$, and
$$\sum_i\lambda_i+\sum_j\mu_j-\left[\left(\frac{a+b-1}{2}\right)^2\right]=n.$$
The integer $n$ is called the rank of the symbol
$\lambda\choose\mu$. The degree of the unipotent character
$\chi^{\lambda,\mu}$ corresponding to the symbol $\lambda\choose\mu$
is equal to
$$\frac{(q^n-1)\prod_{i=1}^{n-1}(q^{2i}-1)\prod_{i\p<i}(q^{\lambda_i}-q^{\lambda_{i\p}})
\prod_{j\p<j}(q^{\mu_j}-q^{\mu_{j\p}})\prod_{i,j}(q^{\lambda_i}+q^{\mu_j})}
{2^{c}q^{{a+b-2\choose2}+ {a+b-4\choose2}
+\cdot\cdot\cdot}\prod_i\prod_{k=1}^{\lambda_i}(q^{2k}-1)
\prod_j\prod_{k=1}^{\mu_j}(q^{2k}-1)}.$$ Here $c=(a+b-2)/2$ if
$\lambda\neq\mu$, and $c=a$ if $\lambda=\mu$. In the former case,
$\lambda\choose\mu$ and $\mu\choose\lambda$ determine the same
unipotent character. In the latter case, there are two unipotent
characters with the same symbols.

\begin{proposition}\label{prop2}
Let $G^\ast=(D_n)_{ad}(q)=P(CO_{2n}^+(q)^0)$. Suppose that $n\geq5$,
$(n,q)\neq(5,2)$, and $\chi\in \Irr(G^\ast)$ is unipotent. Then
either $\chi$ is one of characters labeled by $n\choose0$,
$n-1\choose1$, $0\hspace{3pt}1\choose1\hspace{3pt}n$ with degrees
$1$, $(q^n-1)(q^{n-1}+q)/(q^2-1)$, $(q^{2n}-q^2)/(q^2-1)$,
respectively, or $\chi(1)> q^{4n-10}$. Furthermore, when
$(n,q)=(5,2)$, $G^\ast$ has one more character of degree $868$
(corresponding to the symbol
$0\hspace{3pt}1\hspace{3pt}2\hspace{3pt}4\choose-$), which is less
than $q^{4n-10}$.
\end{proposition}

\begin{proof}
When $n=5$, Proposition can be verified directly by using
\cite[Table 30]{CH}. Denote $\mathcal{L}_n=\{\{(n),(0)\},
\{(n-1),(1)\},\{(0,1),(1,n)\}\}$. We will prove by induction on
$n\geq 5$ that if $\chi=\chi^{\lambda,\mu}$ and
$\{\lambda,\mu\}\notin\mathcal{L}_n$, then $\chi(1)>q^{4n-10}$,
provided $(n,q)\neq(5,2)$. The rest of the proof establishes the
induction step for $n\geq 6$.

\medskip

$1)$ At this point we suppose that $(\lambda_1,\mu_1)\neq (1,0)$ and
$\lambda_1\geq 1$ (eventually $b$ may be zero). Consider the
unipotent character $\chi\p$ labeled by the symbol
$\lambda\p\choose\mu\p$ of rank $n-1$, where
$\lambda\p=(\lambda_1-1, \lambda_2,...,\lambda_a)$ and $\mu\p=\mu$.
If $\{\lambda,\mu\}\in \mathcal{L}_{n-1}$ then
${\lambda\choose\mu}={2\choose n-2}$ or
$2\hspace{3pt}n-1\choose0\hspace{3pt}1$ and therefore
$\chi(1)>q^{4n-10}$. So we can suppose that $\lambda\p\choose\mu\p$
is not in $\mathcal{L}_{n-1}$. By induction hypothesis,
$\chi\p(1)>q^{4(n-1)-10}$. The condition
$n=\sum_i\lambda_i+\sum_j\mu_j-\frac{(a+b)^2-2(a+b)}{4}\geq
a\lambda_1+\frac{(a-b)^2}{4}$ implies that $n-\lambda_1\geq0$. If
$\lambda_1\leq n-4$ then
$$\frac{\chi(1)}{\chi\p(1)}\geq \frac{(q^n-1)(q^{n-1}+1)}{2(q^{2\lambda_1}-1)}\cdot
\prod_{i=2}^a\frac{q^{\lambda_i}-q^{\lambda_1}}{q^{\lambda_i}-q^{\lambda_1-1}}
\cdot\prod_{j=1}^b\frac{q^{\lambda_1}+q^{\mu_j}}{q^{\lambda_1-1}+q^{\mu_j}}>$$
$$>\frac{(q^n-1)(q^{n-1}+1)}{4(q^{2\lambda_1}-1)}
\geq\frac{(q^n-1)(q^{n-1}+1)}{4(q^{2n-8}-1)}>q^4,$$ so
$\chi(1)>q^{4n-10}$. If $n-\lambda_1\leq3$ then $a=1$ and
$(a-b)^2/4\leq3$. Since $a-b\equiv0(\bmod$ $4)$, $a=b=1$. Note that
$\{\lambda,\mu\}$ is not in $\mathcal{L}_n$. So
${\lambda\choose\mu}={n-2\choose2}$ or $n-3\choose3$. Checking
directly, we have $\chi(1)>q^{4n-10}$ as required.

The case $(\lambda_1,\mu_1)\neq (0,1)$ and $\mu_1\geq1$ can be
reduced to $1)$ by interchanging $\lambda$ and $\mu$.

\medskip

$2)$ Now we consider the case where $b=0$ and $\lambda_1=0$. Then
$a\equiv0(\bmod$ $4)$. First we suppose that $\lambda\neq (0, 1,
..., a-1)$. Then there exists $i\geq 2$ such that
$\lambda_i-\lambda_{i-1}\geq 2$. We choose $i$ to be smallest
possible. Consider the character $\chi\p$ corresponding to the
symbol $\lambda\p\choose-$, where $\lambda\p=(..., \lambda_{i-1},
\lambda_i-1, \lambda_{i+1}, ...)$. Note that
$n-\lambda_i\geq(a-2)^2/4$. Therefore, $\chi(1)/\chi\p(1)>
(q^n+1)(q^{n-1}-1)/2(q^{2\lambda_i}-1)>q^4$ if $a\geq 8$. It is
obvious that ${\lambda\p\choose-}\notin\mathcal{L}_{n-1}$. If
$(n,q)=(6,2)$ and $\lambda\p=(0,1,2,5)$ or $(0,1,3,4)$, we can check
directly that $\chi(1)>q^{4n-10}$. Otherwise, the induction
hypothesis $\chi\p(1)>q^{4n-14}$ implies that $\chi(1)>q^{4n-10}$.
If $a=4$ and $n-\lambda_i\geq3$, then again $\chi(1)/\chi\p(1)>q^4$,
so $\chi(1)>q^{4n-10}$. If $a=4$ and $n-\lambda_i\leq2$, then $i=a$
and $\lambda=(0,1,2,n-1)$. Direct computation shows that
$\chi(1)>q^{4n-10}$ and we are done. Next we suppose that
$\lambda=(0, 1, ..., a-1)$. Then $\chi(1)=(q^n-1)f(a)$ and
$n=a^2/4$, where $f(a)$ is the function defined in
(\ref{functionf}). Note that $n\geq 6$, so $a\geq8$. We already
showed that $f(a)>q^{3n}$ for $a\geq6$. Hence, $\chi(1)>q^{4n-10}$
as desired.

Again, the case $a=0$ and $\mu_1=0$ can be reduced to $2)$ by
interchanging $\lambda$ and $\mu$. Therefore, it remains to consider
the case where $a,b\geq 1$ and $(\lambda_1,\mu_1)=(0,1)$ or $(1,0)$.
By a similar reason as above, it is enough to suppose that
$(\lambda_1,\mu_1)=(1,0)$.

\medskip

$3)$ Here we suppose that $\lambda\neq(1,2,...,a)$. Then there
exists $i\geq2$ such that $\lambda_i\geq\lambda_{i-1}+2$. Choose $i$
to be smallest possible. Consider $\chi\p$ labeled by
$\lambda\p\choose\mu\p$ of rank $n-1$, with
$\lambda\p=(\lambda_1,...,\lambda_{i-1},\lambda_i-1,\lambda_{i+1},...,\lambda_a)$
and $\mu\p=\mu$. By induction hypothesis, $\chi\p(1)>q^{4n-14}$. Set
$$U_1=\prod_{i\p<i}\frac{q^{\lambda_i}-q^{\lambda_{i\p}}}{q^{\lambda_i-1}-q^{\lambda_{i\p}}},
U_2=\prod_{i\p>i}\frac{q^{\lambda_{i\p}}-q^{\lambda_{i}}}{q^{\lambda_{i\p}}-q^{\lambda_i-1}},
U_3=\prod_{j}\frac{q^{\lambda_i}+q^{\mu_{j}}}{q^{\lambda_i-1}+q^{\mu_{j}}}.$$
If $n-\lambda_i\geq2$ then
$$\frac{\chi(1)}{\chi\p(1)}=\frac{(q^n-1)(q^{n-1}+1)}{q^{2\lambda_i}-1}U_1U_2U_3>
\frac{(q^n-1)(q^{n-1}+1)}{q^{2\lambda_i}-1}\cdot\frac{q^{\lambda_i}-q}{q^{\lambda_i-1}-q}\cdot\frac{1}{2}
\cdot\frac{q^{\lambda_i}+1}{q^{\lambda_i-1}+1}>$$$$>\frac{(q^n-1)(q^{n-1}+1)q^2}{2(q^{2\lambda_i}-1)}
\geq\frac{(q^n-1)(q^{n-1}+1)q^2}{2(q^{2n-4}-1)}>q^4,$$ so
$\chi(1)>q^{4n-10}$. If $n-\lambda_i\leq1$ then $i=a=b$. Note that
$n-\lambda_i\geq0$. If $n-\lambda_i=1$ then
$${\lambda\choose\mu}={1\hspace{6pt}2\hspace{6pt}...\hspace{6pt}a-1\hspace{6pt}n-1
\choose0\hspace{6pt}1\hspace{6pt}...\hspace{6pt}a-2\hspace{6pt}a}.$$
The case $a=2$ can be checked directly. So we can suppose $a\geq3$.
Then
$$U_1\geq
\frac{(q^{n-1}-q)(q^{n-1}-q^2)}{(q^{n-2}-q)(q^{n-2}-q^2)}, U_2=1,
U_3\geq \frac{q^{n-1}+1}{q^{n-2}+1}.$$ It follows that
$$\frac{\chi(1)}{\chi\p(1)}\geq\frac{q^n-1}{q^{n-1}-1}\cdot
\frac{(q^{n-1}-q)(q^{n-1}-q^2)}{(q^{n-2}-q)(q^{n-2}-q^2)}\cdot
\frac{q^{n-1}+1}{q^{n-2}+1}>q^4,$$ so $\chi(1)>q^{4n-10}$ as
required. If $n-\lambda_i=0$ then
$${\lambda\choose\mu}={1\hspace{6pt}2\hspace{6pt}...\hspace{6pt}a-1\hspace{6pt}n
\choose0\hspace{6pt}1\hspace{6pt}...\hspace{6pt}a-2\hspace{6pt}a-1}.$$
Since $\{\lambda,\mu\}$ is not in $\mathcal{L}_n$, $a\geq3$. Then
$$\frac{\chi(1)}{\chi\p(1)}\geq\frac{(q^n-1)(q^{n-1}+1)}{q^{2n}-1}\cdot
\frac{(q^{n}-q)(q^{n}-q^2)}{(q^{n-1}-q)(q^{n-1}-q^2)}\cdot
\frac{(q^{n}+1)(q^n+q)(q^n+q^2)}{(q^{n-1}+1)(q^{n-1}+q)(q^{n-1}+q^2)}$$$$=
\frac{(q^{n}-q)(q^{n}-q^2)}{(q^{n-1}-q)(q^{n-1}-q^2)}\cdot
\frac{(q^n+q)(q^n+q^2)}{(q^{n-1}+q)(q^{n-1}+q^2)}>q^4,$$ and
therefore $\chi(1)>q^{4n-10}$.

\medskip

$4)$ Finally, we suppose that $\mu_1=0$ and $\lambda=(1,2,...,a)$.
First we consider the case where $\mu\neq (0,1,...,b-1)$. Then there
exists an index $j\geq2$ such that $\mu_j\geq\mu_{j-1}+2$. Consider
the unipotent character $\chi\p$ labeled by $\lambda\p\choose\mu\p$
of rank $n-1$, with $\lambda\p=\lambda$ and
$\mu\p=(\mu_1,...,\mu_{j-1},\mu_j-1,\mu_{j+1},...,\mu_b)$. If
$\{\lambda\p,\mu\p\}\in\mathcal{L}_{n-1}$ then
${\lambda\p\choose\mu\p}={1\hspace{3pt}n-1 \choose0\hspace{3pt}2}$
and hence $\chi(1)>q^{4n-10}$ by direct computation. Therefore we
can assume $\{\lambda\p,\mu\p\}$ is not in $\mathcal{L}_{n-1}$. By
the induction hypothesis, $\chi\p(1)>q^{4n-14}$. Set
$$V_1=\prod_{j\p<j}\frac{q^{\mu_j}-q^{\mu_{j\p}}}{q^{\mu_j-1}-q^{\lambda_{j\p}}},
V_2=\prod_{j\p>j}\frac{q^{\mu_{j\p}}-q^{\mu_{j}}}{q^{\mu_{j\p}}-q^{\mu_j-1}},
V_3=\prod_{i}\frac{q^{\lambda_i}+q^{\mu_{j}}}{q^{\lambda_i}+q^{\mu_{j}-1}}.$$
Remark that $V_1>q, V_2>1/2, V_3>1$ and $n-\mu_j\geq1$. If
$n-\mu_j\geq3$, we have
$$\frac{\chi(1)}{\chi\p(1)}=\frac{(q^n-1)(q^{n-1}+1)}{q^{2\mu_j}-1}V_1V_2V_3
>\frac{(q^n-1)(q^{n-1}+1)q}{2(q^{2\mu_j}-1)}\geq \frac{(q^n-1)(q^{n-1}+1)q}{2(q^{2n-6}-1)}>q^4,$$ so
$\chi(1)>q^{4n-10}$. If $n-\mu_j=2$ then $j=b$ and $V_2=1$. We still
have $\chi(1)/\chi\p(1)>(q^n-1)(q^{n-1}+1)q/(q^{2n-4}-1)>q^4$ and we
are done. The last case we need to handle is $n-\mu_j=1$. Note that
$j=b$ and $V_2=1$ in this case. If $b\geq3$ then
$$V_1\geq\frac{(q^{n-1}-1)(q^{n-1}-q)}{(q^{n-2}-1)(q^{n-2}-q)}>
\frac{q(q^{n-1}-q)}{(q^{n-2}-q)}, V_3\geq
\frac{q^{n-1}+q}{q^{n-2}+q}.$$ It follows that
$$\frac{\chi(1)}{\chi\p(1)}>\frac{(q^n-1)(q^{n-1}+1)}{q^{2n}-1}\cdot
\frac{q(q^{n-1}-q)}{(q^{n-2}-q)}\cdot\frac{q^{n-1}+q}{q^{n-2}+q}>q^4,$$
so $\chi(1)>q^{4n-10}$. If $b\leq2$ then $a=b=2$ and
${\lambda\choose\mu}={1\hspace{3pt}2 \choose0\hspace{3pt}n-1}$,
hence $\chi(1)>q^{4n-10}$.

Next we consider the case where $\mu=(0,1,...,b-1)$. Consider the
unipotent character $\chi\p$ labeled by $\lambda\p\choose\mu\p$ of
rank $n-1$, with $\lambda\p=(1,2,...,a-1)$ and $\mu\p=(0,1,...,b-2)$
(if $a=1$, $\lambda\p$ is just empty; if $b=1$, $\mu\p$ is just
empty). It is obvious that $\{\lambda\p,\mu\p\}$ is not in
$\mathcal{L}_{n-1}$ and therefore $\chi\p(1)>q^{4n-14}$ by the
induction hypothesis. Furthermore,
$$\frac{\chi(1)}{\chi\p(1)}=\frac{q^{a+b-2}(q^n-1)(q^{n-1}+1)}{(q^a-1)(q^{b-1}+1)}>q^{2n-3}>q^4.$$
Consequently, $\chi(1)>q^{4n-10}$.
\end{proof}

%%% -----------------------------------------------------------------------------------

\section{Odd characteristic symplectic groups}
The low-dimensional unipotent characters of $G=Sp_{2n}(q)$, $q$ odd
(up to degree $(q^{2n-2}-1)(q^{2n}-1)(q^{n-2}-1) (q^{n-2}-q^3)/2$
$(q^2-1)(q^4-1)(q^3+1)$) are already classified in Proposition
\ref{unipotentpro}. Therefore, Theorem \ref{symplectic} will be
proved if we can classify non-unipotent characters of $G$ of degrees
up to $q^{4n-10}(q^n-1)/2$. We know that the dual group of $G$ is
$G^\ast=SO_{2n+1}(q)$ and every non-unipotent character $\chi\in
\Irr(G)$ is parameterized by a pair $((s),\psi)$ where $(s)$ is a
nontrivial geometric conjugacy class of a semi-simple element $s\in
G^\ast$ and $\psi$ is an irreducible unipotent character of the
centralizer $C:=C_{G^\ast}(s)$. Moreover,
$\chi(1)=(G^\ast:C)_{p'}\psi(1)$.

\begin{proposition} \label{non-unipro}Let $G=Sp_{2n}(q)$ where
$n\geq 6$ and $q$ is an odd prime power. Suppose that $\alpha,
\alpha_1, \alpha_2=\pm$ and $\chi\in \Irr(G)$ is not unipotent of
degree less than $q^{4n-10}(q^n-1)/2$. Then one of the following
holds:
\begin{enumerate}
   \item[1)] $\chi(1)=(q^{2n-2}-1)(q^{2n}-1)/|A|_{p\p}$, where $A=
GL^{\alpha}_2(q)$, $\mathbb{Z}_{q^2-\alpha}$, or
$\mathbb{Z}_{q-\alpha_1}\times \mathbb{Z}_{q-\alpha_2}$, and $\chi$
is parameterized by $((s), 1_C)$, where $C\simeq SO_{2n-3}(q)\times
A$ or $\chi(1)=q(q^{2n-2}-1)(q^{2n}-1)/(q-\alpha)(q^2-1)$, and
$\chi$ is parameterized by $((s), 1_{SO_{2n-3}(q)}\otimes \lambda)$,
where $C\simeq SO_{2n-3}(q)\times GL_2^{\alpha}(q)$ and $\lambda$ is
the unique nontrivial unipotent character of degree $q$ of
$GL_2^{\alpha}(q)$.
  \item[2)] $\chi(1)=(q^{n-1}+\alpha_1)(q^{2n}-1)/2(q- \alpha_2)$, and
$\chi$ is parameterized by $((s), \lambda\otimes
1_{\mathbb{Z}_{q-\alpha_2}})$, where $C\simeq
O^{\alpha_1}_{2n-2}(q)\times \mathbb{Z}_{q-\alpha_2}$ and $\lambda$
is one of two extensions of $1_{SO^{\alpha_1}_{2n-2}(q)}$ to
$O^{\alpha_1}_{2n-2}(q)$.
  \item[3)] $\chi(1)=(q^{2n}-1)\omega(1)/(q-\alpha)$, and $\chi$ is
parameterized by $((s),\omega\otimes 1_{\mathbb{Z}_{q-\alpha}})$,
where $C\simeq SO_{2n-1}(q)\times \mathbb{Z}_{q-\alpha}$ and
$\omega$ is one of unipotent characters of $SO_{2n-1}(q)$ of degrees
$1$, $(q^{n-1}-1)(q^{n-1}-q)/2(q+1)$,
$(q^{n-1}+1)(q^{n-1}+q)/2(q+1)$, $(q^{n-1}+1)(q^{n-1}-q)/2(q-1)$, or
$(q^{n-1}-1)(q^{n-1}+q)/2(q-1)$.
  \item[4)] $\chi(1)=(q^{2n-2}-1)(q^{2n}-1)/2(q-\alpha_1)(q-\alpha_2)$, and
$\chi$ is parameterized by $((s),\lambda\otimes
1_{SO_{2n-3}(q)\times\mathbb{Z}_{q-\alpha_2}})$, where $C\simeq
O_2^{\alpha_1}(q)\times SO_{2n-3}(q)\times \mathbb{Z}_{q-\alpha_2}$
and $\lambda$ is one of two extensions of $1_{SO_2^{\alpha_1}(q)}$
to $O_2^{\alpha_1}(q)$.
  \item[5)] $\chi(1)= (q^{2n}-1)\omega(1)/2(q-\alpha)$, and $\chi$ is
parameterized by $((s),\omega\otimes \lambda)$, where $C\simeq
SO_{2n-1}(q)\times O_2^\alpha(q)$, $\omega$ is one of unipotent
characters of $SO_{2n-1}(q)$ of degrees $1$,
$(q^{n-1}-1)(q^{n-1}-q)/2(q+1)$, $(q^{n-1}+1)(q^{n-1}+q)/2(q+1)$,
$(q^{n-1}+1)(q^{n-1}-q)/2(q-1)$, or $(q^{n-1}-1)(q^{n-1}+q)/2(q-1)$,
and $\lambda$ is one of two extensions of $1_{SO_2^\alpha(q)}$ to
$O_2^{\alpha}(q)$.
  \item[6)] $\chi(1)=(q^n+1)\lambda(1)/2$, and $\chi$ is parameterized by
$((s),\lambda)$, where $C\simeq O_{2n}^+(q)$ and $\lambda$ is an
extension of one of three unipotent characters of $SO_{2n}^+(q)$ of
degrees $1$, $(q^n-1)(q^{n-1}+q)/(q^2-1)$ or $(q^{2n}-q^2)/(q^2-1)$
to $O_{2n}^+(q)$.
  \item[7)] $\chi(1)=(q^n-1)\lambda(1)/2$, and $\chi$ is parameterized by
$((s),\lambda)$, where $C\simeq O_{2n}^-(q)$ and $\lambda$ is an
extension of one of three unipotent characters of $SO_{2n}^-(q)$ of
degrees $1$, $(q^n+1)(q^{n-1}-q)/(q^2-1)$ or $(q^{2n}-q^2)/(q^2-1)$
to $O_{2n}^-(q)$.
  \item[8)] $\chi(1)=(q^{n-1}+\alpha)(q^{2n}-1)\omega(1)/2(q^2-1)$, and
$\chi$ is parameterized by $((s),\omega\otimes\lambda)$, where
$C\simeq SO_3(q)\times O_{2n-2}^\alpha(q)$, $\omega=1_{SO_3(q)}$ or
the unique nontrivial unipotent character of degree $q$ of
$SO_3(q)$, and $\lambda$ is one of two extensions of
$1_{SO_{2n-2}^\alpha(q)}$ to $O_{2n-2}^\alpha(q)$.
  \item[9)] $\chi(1)=(q^{2n-2}-1)(q^{2n}-1)\lambda(1)/2(q^2-1)(q^2-\alpha)$, and
$\chi$ is parameterized by $((s),1_{SO_{2n-3}(q)}\otimes \lambda)$,
where $C\simeq SO_{2n-3}(q)\times O_4^\alpha(q)$ and $\lambda$ is
one of two extensions of $1_{SO_4^\alpha(q)}$ or the unique
unipotent character of degree $q^2$ of $SO_4^\alpha(q)$ to
$O_4^\alpha(q)$ or $\lambda$ is the unique unipotent character of
degree $2q$ of $O_4^+(q)$.
\end{enumerate}
\end{proposition}

\begin{proof} For short notation, we set $D(n):=(q^n-1)q^{4n-10}/2$.
By Lemma \ref{lemma1}, $C\simeq A\times B$ where $A\simeq
\prod_{i=1}^{t}GL_{a_i}^{\alpha_i}(q^{k_i})$, $a_i,k_i\in
\mathbb{N}$, $\alpha_i=\pm1$, $\sum_{i=1}^{t}k_ia_i=n-m$ and
$B=SO_{2k+1}(q)\times O^{\pm}_{2m-2k}(q)$, $0\leq k\leq m\leq n$. It
is easy to show that $|B|_{p\p}\leq \prod_{i=1}^{m}(q^{2i}-1)$ and
$|A|_{p\p}\leq \prod_{i=1}^{n-m}(q^i-(-1)^i)$. Therefore,
$$\frac{|G|_{p\p}}{|C_{G^\ast}(s)|_{p\p}}\geq
\frac{(q^{2(m+1)}-1)\cdot\cdot\cdot (q^{2n}-1)}
{(q+1)(q^2-1)\cdot\cdot\cdot (q^{n-m}-(-1)^{n-m})}:=f(m,n).$$

\medskip

$1)$ Case $m\leq n-3$. Since the function $f(x,n)$ is first
increasing and then decreasing in the interval $[0,n]$, we have
$f(m,n)\geq \mathrm{min}\{f(0,n), f(n-3,n)\}$. It is easy to check
that both $f(0,n)$ and $f(n-3,n)$ are greater than $D(n)$.

\medskip

$2)$ Case $m=n-2$. First we suppose that $k=n-2$. Then
$\chi(1)=\frac{(q^{2n-2}-1)(q^{2n}-1)}{|A|_{p\p}}\psi(1)$, where
$A=GL^{\alpha}_2(q)$, $GL^{\alpha}_1(q^2)$ or
$GL_1^{\alpha_1}(q)\times GL_1^{\alpha_2}(q)$, $\alpha, \alpha_1,
\alpha_2:=\pm1$ and $\psi$ is a unipotent character of
$C=SO_{2n-3}(q)\times A$. Suppose that $\psi=\omega\otimes \lambda$
where $\omega, \lambda$ are unipotent characters of $SO_{2n-3}(q),
A$, respectively. If $\omega$ is nontrivial then $\omega(1)\geq
(q^{n-2}-1)(q^{n-2}-q)/2(q+1)$ by \cite[Proposition 5.1]{TZ1}. Then
$$\chi(1)\geq
\frac{(q^{2n-2}-1)(q^{2n}-1)(q^{n-2}-1)(q^{n-2}-q)}{2(q+1)^2(q^2-1)}>D(n).$$
If $\omega$ is trivial then $1)$ holds since $GL^{\alpha}_1(q^2)$ or
$GL_1^{\alpha_1}(q)\times GL_1^{\alpha_2}(q)$ has only one unipotent
character, which is trivial, while $A=GL^{\alpha}_2(q)$ has two
unipotent characters, the trivial one and the one of degree $q$.

Next, we consider $k=0$ then
$$\chi(1)\geq\frac{(q^{n-2}-1)(q^{2n-2}-1)(q^{2n}-1)}{2(q+1)(q^2-1)}>D(n).$$

Finally, if $1\leq k\leq n-3$, we have
$$\chi(1)\geq \frac{(q^{2n-2}-1)(q^{2n}-1)}{(q+1)(q^2-1)}\cdot \frac{\prod^{n-2}
_{i=k+1}(q^{2i}-1)\cdot
(q^{n-2-k}-1)}{2\prod_{i=1}^{n-2-k}(q^{2i}-1)}
$$
$$> \frac{(q^{2n-2}-1)(q^{2n}-1)}{(q+1)(q^2-1)}\cdot\frac{1}{2}q^{2k(n-2-k)}
(q^{n-2-k}-1)$$
$$> \frac{(q^{2n-2}-1)(q^{2n}-1)}{(q+1)(q^2-1)}\cdot
\frac{1}{2}q^{2(n-3)}> D(n).$$

\medskip

$3)$ Case $m=n-1$. Assume $k=0$. Then
$\chi(1)=\frac{(q^{n-1}+\alpha_1)(q^{2n}-1)}{2(q-
\alpha_2)}\psi(1)$, where $\psi$ is a unipotent character of
$C\simeq O_{2n-2}^{\alpha_1}(q)\times
GL_1^{\alpha_2}(q)=(SO^{\alpha_1}_{2n-2}(q)\cdot 2)\times
GL_1^{\alpha_2}(q)$, $\alpha_1, \alpha_2=\pm1$. Note that
$GL_1^{\alpha_2}(q)\simeq \mathbb{Z}_{q-\alpha_2}$ has only one
unipotent character, which is the trivial one. Suppose that $\psi$
is an irreducible constituent of
$(\omega)^{SO^{\alpha_1}_{2n-2}(q)\cdot 2}\otimes
1_{GL_1^{\alpha_2}(q)}$ where $\omega$ is a unipotent character of
$SO^{\alpha_1}_{2n-2}(q)$. By Proposition \ref{DM}, the unipotent
characters of $SO^{\pm}_{2n-2}(q)$ are of the form $\theta\circ f$,
where $\theta$ runs over the unipotent characters of
$P(CO^{\pm}_{2n-2}(q)^0)$ and $f$ is the canonical homomorphism
$f:SO^{\pm}_{2n-2}\rightarrow P(CO^{\pm}_{2n-2}(q)^0)$. In
particular, degrees of the unipotent characters of
$SO^{\pm}_{2n-2}(q)$ are the same as those of
$P(CO^{\pm}_{2n-2}(q)^0)$. Therefore, by Propositions 7.1, 7.2 of
\cite{TZ1}, if $\omega$ is nontrivial then $\omega(1)\geq
(q^{n-1}+1)(q^{n-2}-q)/(q^2-1)$. Then
$$\chi(1)\geq
\frac{(q^{n-1}+\alpha_1)(q^{2n}-1)}{2(q- \alpha_2)}\cdot
\frac{(q^{n-1}+1)(q^{n-2}-q)}{q^2-1}> D(n).$$ Suppose that $\omega$
is trivial. Since
$(1_{SO^{\alpha_1}_{2n-2}(q)})^{O^{\alpha_1}_{2n-2}(q)}$ is the sum
of two different irreducible characters of degree $1$,
${O^{\alpha_1}_{2n-2}(q)}$ has two unipotent characters of degree
$1$. So, in this case, $\psi=\lambda\otimes 1_{GL_1^{\alpha_2}(q)}$
where $\lambda$ is one of two extensions of
$1_{SO^{\alpha_1}_{2n-2}(q)}$ to $O^{\alpha_1}_{2n-2}(q)$ and
therefore $2)$ holds.

Assume $k=n-1$. Then $\chi(1)=\frac{q^{2n}-1}{q-\alpha}\psi(1)$,
$\alpha=\pm1$, where $\psi$ is a unipotent character of
$SO_{2n-1}(q)\times GL_1^\alpha(q)$. Again, we have
$\psi=\omega\otimes 1_{GL_1^\alpha(q)}$ where $\omega$ is a
unipotent character of $SO_{2n-1}(q)$. By Corollary \ref{corollary},
$\omega$ is either one of characters of degrees $1$,
$(q^{n-1}-1)(q^{n-1}-q)/2(q+1)$, $(q^{n-1}+1)(q^{n-1}+q)/2(q+1)$,
$(q^{n-1}+1)(q^{n-1}-q)/2(q-1)$, $(q^{n-1}-1)(q^{n-1}+q)/2(q-1)$ or
$\omega(1)>q^{4n-12}$. In the former case, $3)$ holds. In the latter
case, $\chi(1)> \frac{q^{2n}-1}{q-\alpha}q^{4n-12}>D(n)$.

Assume $k=1$. Then $$\chi(1)\geq (G^\ast:C)_{p\p}\geq
\frac{(q^{n-2}-1)(q^{2n-2}-1)(q^{2n}-1)}{2(q+1)(q^2-1)}>D(n).$$

Assume $k=n-2$. Then
$$\chi(1)=\frac{(q^{2n-2}-1)(q^{2n}-1)}{2(q-\alpha_1)(q-\alpha_2)}\psi(1),$$
where $\psi$ is an irreducible unipotent character of $C\simeq
SO_{2n-3}(q)\times O_2^{\alpha_1}(q)\times GL_1^{\alpha_2}(q)$,
$\alpha_1,\alpha_2=\pm$. Note that both $SO_2^{\alpha_1}(q))$ and
$GL_1^{\alpha_2}(q)$ have only one unipotent character, which is the
trivial one. Therefore, $\psi=\omega\otimes \lambda\otimes
1_{GL_1^{\alpha_2}(q)}$, where $\omega$ is a unipotent character of
$SO_{2n-3}(q)$ and $\lambda$ is one of two extensions of
$1_{SO_2^{\alpha_1}(q)}$ to $O_2^{\alpha_1}(q)$. If $\omega$ is
trivial then $4)$ holds. If $\omega$ is nontrivial then
$\omega(1)\geq (q^{n-2}-1)(q^{n-2}-q)/2(q+1)$ by Proposition 5.1 of
\cite{TZ1}. Then,
$$\chi(1)\geq \frac{(q^{2n}-1)(q^{2n-2}-1)}{2(q+1)^2}\cdot
 \frac{(q^{n-2}-1)(q^{n-2}-q)}{2(q+1)}>D(n).$$

Finally, assume $2\leq k\leq n-3$. We have
$$\chi(1)\geq \frac{q^{2n}-1}{q+1}\cdot
\frac{(q^{2(k+1)}-1)\cdot\cdot\cdot
(q^{2(n-1)}-1)}{2(q^2-1)\cdot\cdot\cdot(q^{2(n-1-k)}-1)}\cdot
(q^{n-1-k}-1)$$
$$\geq \frac{q^{2n}-1}{q+1}\cdot \frac{1}{2}q^{2k(n-1-k)}\cdot (q^{n-1-k}-1)
\geq \frac{q^{2n}-1}{q+1}\cdot \frac{1}{2}q^{4(n-3)}\cdot
(q^2-1)>D(n).$$

\medskip

$4)$ Case $m=n$. Since $\chi$ is not unipotent, $k<n$. We have
$$\chi(1)=\frac{\prod_{i=k+1}^n(q^{2i}-1)\cdot (q^{n-k}+\alpha)}
{2\prod_{i=1}^{n-k}(q^{2i}-1)}\psi(1),$$ where $\psi$ is a unipotent
character of $C\simeq SO_{2k+1}(q)\times O_{2(n-k)}^\alpha(q),
\alpha=\pm$. So, $\chi(1)>\frac{1}{2}q^{2k(n-k)}(q^{n-k}- 1)$.

Assume $k=n-1$. Then $\chi(1)=\frac{q^{2n}-1}{2(q-\alpha)}\psi(1)$.
We have $\psi=\omega\otimes \lambda$, where $\omega$ is a unipotent
characters of $SO_{2n-1}(q)$ and $\lambda$ is one of two extensions
of $1_{SO_2^{\alpha_1}(q)}$ to $O_2^{\alpha_1}(q)$. By Corollary
\ref{corollary}, $\omega$ is either one of characters of degrees
$1$, $(q^{n-1}-1)(q^{n-1}-q)/2(q+1)$,
$(q^{n-1}+1)(q^{n-1}+q)/2(q+1)$, $(q^{n-1}+1)(q^{n-1}-q)/2(q-1)$,
$(q^{n-1}-1)(q^{n-1}+q)/2(q-1)$ or $\omega(1)>q^{4n-12}$. In the
former case, $5)$ holds. In the latter case, $\chi(1)>
\frac{q^{2n}-1}{2(q-\alpha)}q^{4n-12}>D(n)$.

Assume $k=0$. Then $\chi(1)=\frac{q^n+\alpha}{2}\psi(1)$, where
$\psi$ is a unipotent character of $C\simeq
O_{2n}^\alpha(q)=SO_{2n}^\alpha(q)\cdot 2$. Suppose that $\psi$ is
an irreducible constituent of $\omega^C$, where $\psi$ is a
unipotent character of $SO_{2n}^\alpha(q)$. We have
$\omega=\theta\circ f$ where $f$ is the canonical homomorphism $f:
SO_{2n}^\pm(q)\rightarrow P(CO_{2n}^\pm(q)^0)$ and $\theta$ is a
unipotent character of $P(CO_{2n}^\pm(q)^0)$. There are two cases:

\textbf{Case} $\mathbf{\alpha=+:}$ By Proposition \ref{prop2},
$\theta$ is either one of characters of degrees $1$,
$(q^n-1)(q^{n-1}+q)/(q^2-1)$, $(q^{2n}-q^2)/(q^2-1)$, or
$\theta(1)>q^{4n-10}$. In the latter case,
$$\chi(1)= \frac{q^n+1}{2}\psi(1)\geq
\frac{q^n+1}{2}\omega(1)>\frac{q^n+1}{2}q^{4n-10}> D(n).$$ In the
former case, by \cite{ST}, the character afforded by the rank $3$
permutation module of $O_{2n}^+(q)$ on the set of all singular
$1$-spaces of $\mathbb{F}_q^{2n}$ is the sum of three unipotent
characters of degrees $1$, $(q^n-1)(q^{n-1}+q)/(q^2-1)$, and
$(q^{2n}-q^2)/(q^2-1)$. Hence, $6)$ holds in this case.

\textbf{Case} $\mathbf{\alpha=-:}$ By Proposition \ref{prop1},
$\theta$ is either one of characters of degrees $1$,
$(q^n+1)(q^{n-1}-q)/(q^2-1)$, $(q^{2n}-q^2)/(q^2-1)$, or
$\theta(1)>q^{4n-10}$. In the latter case,
$$\chi(1)= \frac{q^n-1}{2}\psi(1)\geq
\frac{q^n-1}{2}\omega(1)>\frac{q^n-1}{2}q^{4n-10}= D(n).$$ In the
former case, again by \cite{ST}, the character afforded by the rank
$3$ permutation module of $O_{2n}^-(q)$ on the set of all singular
$1$-spaces of $\mathbb{F}_q^{2n}$ is the sum of three unipotent
characters of degrees $1$, $(q^n+1)(q^{n-1}-q)/(q^2-1)$, and
$(q^{2n}-q^2)/(q^2-1)$. Hence, $7)$ holds in this case.

Assume $k=1$. Then
$\chi(1)=\frac{(q^{2n}-1)(q^{n-1}+\alpha)}{2(q^2-1)}\psi(1)$, where
$\psi$ is a unipotent character of $C\simeq SO_{3}(q)\times
O_{2n-2}^\alpha(q)$. Suppose that $\psi$ is an irreducible
constituent of $\omega\otimes \varphi^C$ where $\omega, \varphi$ are
unipotent characters of $SO_{3}(q), SO_{2n-2}^\alpha(q)$,
respectively. Let us consider the first case where $\varphi$ is
nontrivial. Then by Propositions 7.1, 7.2 of \cite{TZ1}, we have
$\varphi(1)\geq (q^{n-1}+1)(q^{n-2}-q)/(q^2-1)$. Then
$$\chi(1)\geq
\frac{(q^{2n}-1)(q^{n-1}-1)(q^{n-1}+1)(q^{n-2}-q)}{2(q^2-1)^2}>D(n).$$
Now suppose $\varphi=1_{SO_{2n-2}^\alpha(q)}$. Since $SO_3(q)$ has
two unipotent characters, the trivial one and the unique nontrivial
unipotent character of degree $q$, $8)$ holds in this case.

Assume $k=2$. Then
$$\chi(1)=\frac{(q^{2n-2}-1)(q^{2n}-1)(q^{n-2}\pm1)}{2(q^2-1)(q^4-1)}\psi(1)>D(n).$$

Assume $k=n-2$. Then
$$\chi(1)=\frac{(q^{2n-2}-1)(q^{2n}-1)}{2(q^2-1)(q^2-\alpha)}\psi(1),$$
where $\psi$ is a unipotent character of $C\simeq SO_{2n-3}(q)\times
O_4^\alpha(q), \alpha=\pm$. Suppose that $\psi$ is an irreducible
constituent of $\omega\otimes \varphi^{O_4^\alpha(q)}$, where
$\omega, \varphi$ are unipotent characters of $SO_{2n-3}(q),
SO_{4}^\alpha(q)$, respectively. If $\omega$ is nontrivial, by
Proposition 5.1 of \cite{TZ1}, $\omega(1)\geq
(q^{n-2}-1)(q^{n-2}-q)/2(q+1)$ and therefore
$$\chi(1)\geq \frac{(q^{2n-2}-1)(q^{2n}-1)(q^{n-2}-1)(q^{n-2}-q)}{4(q+1)(q^4-1)}>D(n).$$
Now we assume $\omega=1_{SO_{2n-3}(q)}$. There are two cases:

\textbf{Case} $\mathbf{\alpha=+:}$ Let
$V:=M_{2\times2}(\mathbb{F}_q)$ be the space of $2\times 2$ matrices
over $\mathbb{F}_q$. Then $V$ is a vector space over $\mathbb{F}_q$
of rank $4$. The determinant function $Q(M)=\mathrm{det}(M)$, $M\in
V$ is a quadratic form on $V$ with Witt index $2$. So the group of
linear transformations $V\rightarrow V$ preserving $Q$ will be
$O_4^+(q)$. Consider an action of group $GL_2(q)\times GL_2(q)$ on
$V$ by $\tau(A,B)(M)=A^{-1}MB$. We see that $\tau(A,B)$ preserves
$Q$ if and only if $\mathrm{det}(A)=\mathrm{det}(B)$, and
$\tau(A,B)$ is the identity if and only if $A=B$ is a scalar matrix.
Moreover, $\tau(A,B)$ has determinant
$\mathrm{det}(A)^{-2}\mathrm{det}(B)^2$, which is $1$ when
$\mathrm{det}(A)=\mathrm{det}(B)$. Therefore, we have a homomorphism
$$\tau: (SL_2(q)\times SL_2(q))\cdot \mathbb{Z}_2\rightarrow SO_4^+(q),$$
where $\mathbb{Z}_2$ is generated by
$$U=\left(\left(\hspace{-2mm}\begin{array}{cc} -1&0\\
0& 1\end{array}\hspace{-2mm}\right),
\left(\hspace{-2mm}\begin{array}{cc} -1&0\\
0& 1\end{array}\hspace{-2mm}\right)\right).$$ We have
$\Ker(\rho)=\mathbb{Z}_2$, which is generated by $(-I,-I)$.
Therefore, $((SL_2(q)\times SL_2(q))\cdot
\mathbb{Z}_2)/\mathbb{Z}_2\simeq 2\cdot(PSL_2(q)\times
PSL_2(q))\cdot2$ is a subgroup of $SO_4^+(q)$. Since
$|2\cdot(PSL_2(q)\times PSL_2(q))\cdot2|=|SO_4^+(q)|=q^2(q^2-1)^2$,
$2\cdot(PSL_2(q)\times PSL_2(q))\cdot2\simeq SO_4^+(q)$. We denote
by $L_1, L_2$ the first and second terms respectively in
$PSL_2(q)\times PSL_2(q)$. Note that $P(CO_4^+(q)^0)$ as well as
$SO_4^+(q)$ have four unipotent characters of degrees $1, q, q,
q^2$. Since $PSL_2(q)$ has only one character of degree $q$ which is
denoted by $\mu_1$ for $L_1$ and $\mu_2$ for $L_2$ and $U$ fixes
$L_1$ and $L_2$, the unipotent characters of $SO_4^+(q)$ must be
extensions of characters $1_{L_1\times L_2}, \mu_1\otimes 1_{L_2},
1_{L_1}\otimes \mu_2$ and $\mu_1\otimes \mu_2$ which are considered
as characters of $2\cdot(L_1\times L_2)$.

Consider the transformation $T: V\rightarrow V$ defined by
$T(M)=M^{T}$, the transpose of $M$. It is clear that $T$ preserves
$Q$ and therefore $T\in O_4^+(q)$. Moreover $\mathrm{det}(T)=-1$. So
$O_4^+(q)=SO_4^+(q)\rtimes <T>$. We have
$T^{-1}\tau(A,B)T(M)=B^TM(A^T)^{-1}$ for every $A,B\in GL_2(q)$. So
$T$ fixes $\mu_1\otimes \mu_2$ and maps one of $\{\mu_1\otimes
1_{L_2}, 1_{L_1}\otimes \mu_2\}$ to the other. In other words, the
unipotent character of degree $q^2$ of $SO_4^+(q)$ has two
extensions to $O_4^+(q)$ and the inductions of two unipotent
characters of degree $q$ of $SO_4^+(q)$ to $O_4^+(q)$ are equal and
irreducible. So $9)$ holds in this case.

\textbf{Case} $\mathbf{\alpha=-:}$ Note that $SO_4^-(q)$ as well as
$P(CO_{4}^-(q)^0)$ have two unipotent characters of degrees $1,
q^2$. It is well known that $SO_4^-(q)\simeq PSL_2(q^2)\times 2_1$
and $(O_4^-(q)=(PSL_2(q^2)\times 2_1)\cdot 2_2\simeq 2_1\times
(PSL_2(q^2)\cdot 2_2)$, where $2_2$ acts trivially on $2_1$. Since
$PSL_2(q^2)$ has only one irreducible character of degree $q^2$
which is unipotent, the unipotent character of degree $q^2$ of
$SO_4^-(q)$ is invariant in $O_4^-(q)$. Therefore, $O_4^-(q)$ has
two unipotent characters of degree $q^2$. So $9)$ also holds in this
case.

Lastly, we assume $3\leq k\leq n-3$. Then $$\chi(1)\geq
\frac{\prod_{i=k+1}^n(q^{2i}-1)\cdot (q^{n-k}\pm
1)}{2\prod_{i=1}^{n-k}(q^{2i}-1)}>$$$$>\frac{1}{2}q^{2k(n-k)}(q^{n-k}\pm1)\geq
\frac{1}{2}q^{6(n-3)}(q^3-1)>D(n).$$
\end{proof}

\noindent \textbf{Counting semi-simple conjugacy classes in}
$\mathbf{G^\ast=SO_{2n+1}(q)}$: In Proposition \ref{non-unipro} and
its proof, we have not shown how to count the number of semi-simple
conjugacy classes $(s)$ in $SO_{2n+1}(q)$ for a certain
$C:=C_{G^\ast}(s)$, which will imply the number of irreducible
characters of $Sp_{2n}(q)$ at each degree. Actually, the way to
count them is pretty similar in all the cases from 1) to 9) in
Proposition \ref{non-unipro}. First, from the structure of
$C_{G^\ast}(s)$, we know the form of the characteristic polynomial
as well as the eigenvalues of $s$. In general, we have
$$\Spec(s)=\{\underbrace{1,...,1}_{2k+1},\underbrace{-1,...,-1}_{2(m-k)},
\underbrace{...}_{2(n-m)}\}.$$ By \cite[(2.6)]{W1}, if $-1\notin
\Spec(s)$, there is exactly one $O_{2n+1}(q)$-conjugacy class of
semi-simple elements $(s)$ for a given $\Spec(s)$. This class is
also an $SO_{2n+1}(q)$-conjugacy class since
$(O_{2n+1}(q):C_{O_{2n+1}(q)}(s))=(SO_{2n+1}(q):C_{SO_{2n+1}(q)}(s))$
(see Lemma \ref{lemma1}). The situation is a little bit different
when $-1\in \Spec(s)$. In that case, again by \cite[(2.6)]{W1},
there are exactly two $O_{2n+1}(q)$-conjugacy class of semi-simple
elements $(s)$ for a given $\Spec(s)$, in which
$C_{O_{2n+1}(q)}(s)\simeq O_{2k+1}(q)\times O^\pm_{2(m-k)}(q)\times
\prod_{i=1}^tGL_{a_i}^{\alpha_i}(q^{k_i})$ ($+$ for one class and
$-$ for the other). These classes are also $SO_{2n+1}(q)$-conjugacy
class by the same reason as before. So, in order to count the number
of semi-simple conjugacy classes $(s)$ with a given centralizer $C$,
we need to count the number of choices of $\Spec(s)$. This is
demonstrated in Tables 2, 3.

We finish this section by two following Corollaries of Theorem
\ref{symplectic}.

\begin{corollary}
Let $\chi$ be an irreducible complex character of $G=Sp_{2n}(q)$,
where $n\geq 6$ and $q$ is an odd prime power. Then $\chi(1)= 1$
($1$ character), $(q^n\pm 1)/2$ ($4$ characters),
$(q^n+\alpha_1)(q^n+\alpha_2q)/2(q+\alpha_1\alpha_2)$ ($4$
characters) $(\alpha_{1,2}=\pm1)$, $(q^{2n}-1)/2(q\pm1)$ ($4$
characters), $(q^{2n}-1)/(q\pm1)$ ($q-2$ characters), or
$\chi(1)\geq (q^{2n}-1)(q^{n-1}-q)/2(q^2-1)$.
\end{corollary}

\begin{corollary}
Let $\chi$ be an irreducible complex character of $G=Sp_{2n}(q)$,
where $n\geq 6$ and $q$ is an odd prime power. Then $\chi(1)= 1$
($1$ character), $(q^n\pm 1)/2$ ($4$ characters),
$(q^n+\alpha_1)(q^n+\alpha_2q)/2(q+\alpha_1\alpha_2)$ ($4$
characters) $(\alpha_{1,2}=\pm1)$, $(q^{2n}-1)/2(q\pm1)$ ($4$
characters), $(q^{2n}-1)/(q\pm1)$ ($q-2$ characters),
$(q^{2n}-1)(q^{n-1}\pm q)/2(q^2-1)$ ($4$ characters),
$(q^{2n}-1)(q^{n-1}\pm1)/2(q^2-1)$ ($4$ characters),
$(q^{2n}-1)(q^{n-1}\pm1)/2(q\pm1)$ ($4q-8$ characters),
$(q^{2n}-q^2)(q^{n}\pm1)/2(q^2-1)$ ($4$ characters),
$q(q^{2n}-1)(q^{n-1}\pm1)/2(q^2-1)$ ($4$ characters), or
$\chi(1)\geq (q^{2n}-1)(q^{n-1}-1)(q^{n-1}-q^2)/2(q^4-1)$.
\end{corollary}

\begin{table}
\caption{Low-dimensional unipotent characters of $Sp_{2n}(q)$,
$n\geq 6$, $q$ odd. See $\S 3.1$ for the explanation of notation.}
\begin{tabular}{lll}
\hline
 Characters & Symbols  & Degrees\\ \hline
 $1_G$&$\left(\hspace{-2.5mm}\begin{array}{c}n\\-\end{array} \hspace{-2.5mm}\right)$  & $1$\\
 $\chi_2$&$\left(\hspace{-2.5mm}\begin{array}{c}0\hspace{3pt} 1 \hspace{3pt}n\\-\end{array}
\hspace{-2.5mm}\right)$&$\dfrac{(q^n-1)(q^n-q)}{2(q+1)}$\\
$\chi_3$&$\left(\hspace{-2.5mm}\begin{array}{c}0\hspace{3pt} 1
\\n\end{array}\hspace{-2.5mm} \right)$&$\dfrac{(q^n+1)(q^n+q)}{2(q+1)}$\\
$\chi_4$&$\left(\hspace{-2.5mm}\begin{array}{c}1\hspace{3pt} n
\\0\end{array} \hspace{-2.5mm}\right)$&$\dfrac{(q^n+1)(q^n-q)}{2(q-1)}$\\
$\chi_5$& $\left(\hspace{-2.5mm}\begin{array}{c}0\hspace{3pt} n
\\1\end{array}\hspace{-2.5mm} \right)$&$\dfrac{(q^n-1)(q^n+q)}{2(q-1)}$\\
 $\chi_6$&$\left(\hspace{-2.5mm}\begin{array}{c}0\hspace{3pt} 2 \hspace{3pt}n-1\\-\end{array}
\hspace{-2.5mm}\right)$&$\dfrac{(q^{2n}-1)(q^{n-1}-1)(q^{n-1}-q^2)}{2(q^4-1)}$\\
$\chi_7$& $\left(\hspace{-2.5mm}\begin{array}{c}0\hspace{3pt}2
\\n-1\end{array} \hspace{-2.5mm}\right)$&$\dfrac{(q^{2n}-1)(q^{n-1}+1)(q^{n-1}+q^2)}{2(q^4-1)}$\\
$\chi_8$& $\left(\hspace{-2.5mm}\begin{array}{c}2\hspace{3pt} n-1
\\0\end{array} \hspace{-2.5mm}\right)$&$\dfrac{(q^{2n}-1)(q^{n-1}+1)(q^{n-1}-q^2)}{2(q^2-1)^2}$\\
$\chi_9$& $\left(\hspace{-2.5mm}\begin{array}{c}0\hspace{3pt} n-1
\\2\end{array} \hspace{-2.5mm}\right) $&$\dfrac{(q^{2n}-1)(q^{n-1}-1)(q^{n-1}+q^2)}{2(q^2-1)^2}$\\
$\chi_{10}$&$\left(\hspace{-2.5mm}\begin{array}{c}1\hspace{3pt} n-1
\\1\end{array} \hspace{-2.5mm}\right)$&$\dfrac{q(q^{2n}-1)(q^{2n-2}-q^2)}{(q^2-1)^2}$\\
$\chi_{11}$&$\left(\hspace{-2.5mm}\begin{array}{c}0\hspace{3pt} 1
\hspace{3pt} 2
\hspace{3pt}n\\1\end{array} \hspace{-2.5mm}\right)$&$\dfrac{(q^{2n}-q^2)(q^{n}-1)(q^n-q^2)}{2(q^4-1)}$\\
 $\chi_{12}$&$\left(\hspace{-2.5mm}\begin{array}{c}0\hspace{3pt} 1
\hspace{3pt}2\\1\hspace{3pt} n\end{array}
\hspace{-2.5mm}\right)$&$\dfrac{(q^{2n}-q^2)(q^{n}+1)(q^n+q^2)}{2(q^4-1)}$\\
$\chi_{13}$&$\left(\hspace{-2.5mm}\begin{array}{c}1\hspace{3pt} 2
\hspace{3pt}n\\0 \hspace{3pt}1\end{array}
\hspace{-2.5mm}\right)$&$\dfrac{(q^{2n}-q^2)(q^{n}+1)(q^n-q^2)}{2(q^2-1)^2}$\\
$\chi_{14}$&$\left(\hspace{-2.5mm}\begin{array}{c}0\hspace{3pt} 1
\hspace{3pt}n\\1\hspace{3pt} 2\end{array}
\hspace{-2.5mm}\right)$&$\dfrac{(q^{2n}-q^2)(q^{n}-1)(q^n+q^2)}{2(q^2-1)^2}$\\\hline
\end{tabular}
\end{table}

\begin{sidewaystable}
\caption{Low-dimensional irreducible characters of $Sp_{2n}(q)$,
$n\geq 6$, $q$ odd. See $\S 2.1$ and Proposition \ref{non-unipro}
for the explanation of notation.}
\begin{tabular}{lllll}\hline
$\Spec(s)$ & $C_{G^\star}(s)$ & $\psi(1)$ & $\chi(1)$ & $\sharp$ of
characters \\\hline
  \multirow{2}{*}{\begin{tabular}{c}$\{1,...,1,\gamma^k,\gamma^{-k},\gamma^k,\gamma^{-k}\}$\\
  $k\in T_1$\end{tabular}}&
 \multirow{2}{*}{$SO_{2n-3}(q)\times
 GL_{2}(q)$}&$1$&$\dfrac{(q^{2n-2}-1)(q^{2n}-1)}{(q-1)(q^2-1)}$&$(q-3)/2$\\\cline{3-5}
 &&$q$& $\dfrac{q(q^{2n-2}-1)(q^{2n}-1)}{(q-1)(q^2-1)}$&$(q-3)/2$\\\hline
 \multirow{2}{*}{\begin{tabular}{c}$\{1,...,1,\eta^k,\eta^{-k},\eta^k,\eta^{-k}\}$\\
  $k\in T_2$\end{tabular}}&
 \multirow{2}{*}{$SO_{2n-3}(q)\times
 GU_{2}(q)$}&$1$&$\dfrac{(q^{2n-2}-1)(q^{2n}-1)}{(q+1)(q^2-1)}$&$(q-1)/2$\\\cline{3-5}
 &&$q$& $\dfrac{q(q^{2n-2}-1)(q^{2n}-1)}{(q+1)(q^2-1)}$&$(q-1)/2$\\\hline
\begin{tabular}{c}$\{1,...,1,\theta^j,\theta^{qj},\theta^{-j},\theta^{-qj}\}$\\
$j\in R_2$\end{tabular}&$SO_{2n-3}(q)\times
 GL_{1}(q^2)$&$1$&$\dfrac{(q^{2n-2}-1)(q^{2n}-1)}{q^2-1}$&$(q-1)^2/4$\\\hline
 \begin{tabular}{c}$\{1,...,1,\zeta^j,\zeta^{qj},\zeta^{-j},\zeta^{-qj}\}$\\
  $j\in R_1$\end{tabular}&$SO_{2n-3}(q)\times
 GU_{1}(q^2)$&$1$&$\dfrac{(q^{2n-2}-1)(q^{2n}-1)}{q^2+1}$&$(q^2-1)/4$\\\hline
 $\begin{tabular}{c}$\{1,...,1,\gamma^k,\gamma^{-k},\gamma^l,\gamma^{-l}\}$\\
 $k,l\in T_1, k\neq l$\end{tabular}$&$SO_{2n-3}(q)\times
 GL_{1}(q)\times GL_1(q)$&$1$&$\dfrac{(q^{2n-2}-1)(q^{2n}-1)}{(q-1)^2}$&$(q-3)(q-5)/8$\\\hline
 $\begin{tabular}{c}$\{1,...,1,\eta^k,\eta^{-k},\eta^l,\eta^{-l}\}$\\
 $k,l\in T_2, k\neq l$\end{tabular}$&$SO_{2n-3}(q)\times
 GU_{1}(q)\times GU_1(q)$&$1$&$\dfrac{(q^{2n-2}-1)(q^{2n}-1)}{(q+1)^2}$&$(q-1)(q-3)/8$\\\hline
 $\begin{tabular}{c}$\{1,...,1,\gamma^k,\gamma^{-k},\eta^l,\eta^{-l}\}$\\
 $k\in T_1, l\in T_2$\end{tabular}$&$SO_{2n-3}(q)\times
 GL_{1}(q)\times GU_1(q)$&$1$&$\dfrac{(q^{2n-2}-1)(q^{2n}-1)}{q^2-1}$&$(q-1)(q-3)/4$\\\hline
 $\{1,-1,...,-1,\gamma^k,\gamma^{-k}\}$, $k\in T_1$&$O^{\alpha}_{2n-2}(q)\times
 GL_{1}(q)$&$1$&$\dfrac{(q^{2n}-1)(q^{n-1}+\alpha)}{2(q-1)}$&$2(q-3)$\\\hline
 $\{1,-1,...,-1,\eta^k,\eta^{-k}\}$, $k\in T_2$&$O^{\alpha}_{2n-2}(q)\times
 GU_{1}(q)$&$1$&$\dfrac{(q^{2n}-1)(q^{n-1}+\alpha)}{2(q+1)}$&$2(q-1)$\\\hline
$\{1,...,1,\gamma^k,\gamma^{-k}\}$, $k\in T_1$&$SO_{2n-1}(q)\times
 GL_{1}(q)$&$(\star)$&$\dfrac{q^{2n}-1}{q-1}\psi(1)$&$(q-3)/2$\\\hline
 $\{1,...,1,\eta^k,\eta^{-k}\}$, $k\in T_2$&$SO_{2n-1}(q)\times
 GU_{1}(q)$&$(\star)$&$\dfrac{q^{2n}-1}{q+1}\psi(1)$&$(q-1)/2$\\\hline
\end{tabular}
\end{sidewaystable}

\begin{sidewaystable}
\caption{Low-dimensional irreducible characters of $Sp_{2n}(q)$,
$n\geq 6$, $q$ odd. See $\S 2.1$ and Proposition \ref{non-unipro}
for the explanation of notation.}
\begin{tabular}{lllll}\hline
$\Spec(s)$ & $C_{G^\star}(s)$ & $\psi(1)$ & $\chi(1)$ & $\sharp$ of chars \\
\hline\begin{tabular}{c}$\{1,...,1,-1,-1,\gamma^k,\gamma^{-k}\}$\\
$k\in T_1$\end{tabular}&
\begin{tabular}{c}$SO_{2n-3}(q)\times O^{\alpha}_{2}(q)$\\$\times
 GL_{1}(q)$\end{tabular}&$1$&$\dfrac{(q^{2n-2}-1)(q^{2n}-1)}{2(q-\alpha)(q-1)}$
 &$2(q-3)$\\\hline \begin{tabular}{c}$\{1,...,1,-1,-1,\eta^k,\eta^{-k}\}$\\
$k\in T_2$\end{tabular}&
\begin{tabular}{c}$SO_{2n-3}(q)\times O^{\alpha}_{2}(q)$\\$\times
 GU_{1}(q)$\end{tabular}&$1$&$\dfrac{(q^{2n-2}-1)(q^{2n}-1)}{2(q-\alpha)(q+1)}$
 &$2(q-1)$\\\hline$\{1,...,1,-1,-1\}$&$SO_{2n-1}(q)\times
 O^{\alpha}_{2}(q)$&$(\star)$&$\dfrac{q^{2n}-1}{2(q-\alpha)}\psi(1)$&$4$\\\hline
\multirow{3}{*}{$\{1,-1,...,-1\}$}&
 \multirow{3}{*}{$O_{2n}^\alpha(q)$}&$1$&$(q^n+\alpha)/2$&$4$\\\cline{3-5}
 &&$\dfrac{(q^n-\alpha)(q^{n-1}+\alpha
 q)}{q^2-1}$&$\dfrac{(q^{2n}-1)(q^{n-1}+\alpha q)}{2(q^2-1)}$&$4$\\\cline{3-5}
 &&$\dfrac{q^{2n}-q^2}{q^2-1}$&$\dfrac{(q^{2n}-q^2)(q^n+\alpha)}{2(q^2-1)}$&$4$\\\hline
\multirow{2}{*}{$\{1,1,1,-1,...,-1\}$}&
 \multirow{2}{*}{$SO_3(q)\times O_{2n-2}^\alpha(q)$}&$1$&$\dfrac{(q^{2n}-1)(q^{n-1}+\alpha)}{2(q^2-1)}$
 &$4$\\\cline{3-5}
 &&$q$&$\dfrac{q(q^{2n}-1)(q^{n-1}+\alpha)}{2(q^2-1)}$&$4$\\\hline
\multirow{2}{*}{$\{1,...,1,-1,-1,-1,-1\}$}&
 \multirow{2}{*}{$SO_{2n-3}(q)\times O_{4}^\alpha(q)$}&$1$&
 $\dfrac{(q^{2n-2}-1)(q^{2n}-1)}{2(q^2-1)(q^2-\alpha)}$
 &$4$\\\cline{3-5}
 &&$q^2$&$\dfrac{q^2(q^{2n-2}-1)(q^{2n}-1)}{2(q^2-1)(q^2-\alpha)}$&$4$\\\hline
 $\{1,...,1,-1,-1,-1,-1\}$&$SO_{2n-3}(q)\times O_{4}^+(q)$&
 $2q$&$\dfrac{q(q^{2n-2}-1)(q^{2n}-1)}{(q^2-1)^2}$&$1$\\\hline
 \end{tabular}
\begin{tabular}{l} Where $(\star)$ is $1$,
$(q^{n-1}-1)(q^{n-1}-q)/2(q+1)$, $(q^{n-1}+1)(q^{n-1}+q)/2(q+1)$,
$(q^{n-1}+1)(q^{n-1}-q)/2(q-1)$,\\ or
$(q^{n-1}-1)(q^{n-1}+q)/2(q-1)$.
\end{tabular}
\end{sidewaystable}

%%% -----------------------------------------------------------------------------------

\clearpage
\section{Odd characteristic orthogonal groups in odd dimension}
The aim of this section is to classify low-dimensional complex
representations of groups $G=Spin_{2n+1}(q)$ where $q$ is a power of
an odd prime $p$. The dual group $G^\ast$ is the projective
conformal symplectic group $PCSp_{2n}(q)$, which is the quotient of
$\widetilde{G}=CSp_{2n}(q)$ by its center
$Z(\widetilde{G})\simeq\mathbb{Z}_{q-1}$.

\medskip

\textbf{Proof of Theorem \ref{orthogonal-odddimension}}. If $\chi$
is unipotent, Theorem is done by Corollary \ref{corollary}. So we
can assume $\chi$ is not unipotent such that $\chi(1)\leq q^{4n-8}$.
Suppose that $\chi$ is parameterized by $((s^\ast), \psi)$, where
$s^\ast$ is a non-trivial semi-simple element in $G^\ast$. Consider
an inverse image $s$ of $s^\ast$ in $\widetilde{G}$. Then $(s)$ is a
non-trivial conjugacy class of semi-simple elements in
$\widetilde{G}$. Set $C^\ast=C_{G^\ast}(s^\ast)$ and
$C=C_{\widetilde{G}}(s)$. Let $g$ be any element in $\widetilde{G}$
such that the image of $g$ in $G^\ast$ belongs to $C^\ast$. Then
$gsg^{-1}=\mu(g)s$ for some $\mu(g)\in \mathbb{Z}_{q-1}$. Therefore,
$\tau(s)=\tau(gsg^{-1})=\tau(\mu(g)s)=\mu(g)^2\tau(s)$. It follows
that $\mu(g)$ is $1$ or $-1$. Hence
$$(G^\ast:C^\ast)_{p'}\geq (\widetilde{G}:C)_{p'}/2.$$

\textbf{Case} $\mathbf{1}$: If $\tau(s)$ is not a square in
$\mathbb{F}_q$, by Lemma \ref{lemma3}, $C\simeq (Sp_{m}(q^2)\times
\prod_{i=1}^{t}GL_{a_i}^{\alpha_i}(q^{k_i}))\cdot \mathbb{Z}_{q-1}$,
where $\alpha_i=\pm$, $\Sigma_{i=1}^tk_ia_i=n-m$. Therefore,
$$(\widetilde{G}:C)_{p'}\geq\frac{(q^2-1)\cdot\cdot\cdot (q^{2n}-1)}{(q^4-1)
\cdot\cdot\cdot(q^{2m}-1)(q+1)(q^2-1)\cdot\cdot\cdot(q^{n-m}-(-1)^{n-m})}$$
$$ \geq \left\{\begin{array}{ll}\min\{\dfrac{(q^2-1)\cdot\cdot\cdot
(q^{2n}-1)}{(q^4-1) \cdot\cdot\cdot(q^{2n}-1)},
\dfrac{(q^2-1)\cdot\cdot\cdot
(q^{2n}-1)}{(q+1)(q^2-1)\cdot\cdot\cdot(q^{n}-(-1)^{n})}\}&{\text{
if }} 2\mid n\\
\min\{\dfrac{(q^2-1)\cdot\cdot\cdot (q^{2n}-1)}{(q^4-1)
\cdot\cdot\cdot(q^{2n-2}-1)(q+1)}, \dfrac{(q^2-1)\cdot\cdot\cdot
(q^{2n}-1)}{(q+1)(q^2-1)\cdot\cdot\cdot(q^{n}-(-1)^{n})}\}&{\text{
if }} 2\nmid n
 \end{array}\right.$$
 $$=\left\{\begin{array}{ll}\min\{(q^2-1)(q^6-1)\cdot\cdot\cdot(q^{2n-2}-1),
 (q-1)\cdot\cdot\cdot(q^n+(-1)^{n})\}&\text{if } 2\mid n\\
 \min\{(q^2-1)(q^6-1)\cdot\cdot\cdot(q^{2n}-1)/(q+1),
 (q-1)\cdot\cdot\cdot(q^n+(-1)^{n})\}&\text{if } 2\nmid n,
 \end{array}\right.$$
which is greater than $2q^{4n-8}$. Therefore,
$$\chi(1)\geq(G^\ast:C^\ast)_{p'}\geq (\widetilde{G}:C)_{p'}/2>q^{4n-8}.$$

\textbf{Case} $\mathbf{2}$: If $\tau(s)$ is a square in
$\mathbb{F}_q$, multiplying $s$ with a suitable scalar, we can
assume that $\tau(s)=1$ or in other words $s\in Sp_{2n}(q)$. By
Lemma \ref{lemma3}, $C\simeq (Sp_{2k}(q)\times Sp_{2(m-k)}(q) \times
\prod_{i=1}^{t}GL_{a_i}^{\alpha_i}(q^{k_i}))\cdot \mathbb{Z}_{q-1}$,
where $\alpha_i=\pm$, $\Sigma_{i=1}^tk_ia_i=n-m$. It is easy to see
that $|Sp_{2k}(q)\times Sp_{2(m-k)}(q)|_{p'}\leq
\prod_{i=1}^{m}(q^{2i}-1)$ and
$|\prod_{i=1}^{t}GL_{a_i}^{\alpha_i}(q^{k_i})|_{p'}\leq
\prod_{i=1}^{n-m}(q^i-(-1)^i)$. Therefore,
$$(\widetilde{G}:C)_{p'}\geq \frac{(q^{2(m+1)}-1)\cdot\cdot\cdot
(q^{2n}-1)} {(q+1)(q^2-1)\cdot\cdot\cdot
(q^{n-m}-(-1)^{n-m})}=:f(m,n).$$

$1)$ Case $m\leq n-2$. We have $(\widetilde{G}:C)_{p'}\geq
\min\{f(0,n), f(n-2,n)\}$. Since both
$f(0,n)=\prod_{i=1}^n(q^i+(-1)^i)$ and
$f(n-2,n)=(q^{2n-2}-1)(q^{2n}-1)/(q+1)(q^2-1)$ are greater than
$2q^{4n-8}$, $\chi(1)\geq (\widetilde{G}:C)_{p'}/2>q^{4n-8}$.

$2)$ Case $m=n-1$. If $1\leq k\leq n-2$ then
$(\widetilde{G}:C)_{p'}\geq (q^{2n-2}-1)(q^{2n}-1)/(q+1)(q^2-1)$.
Again, $\chi(1)\geq (\widetilde{G}:C)_{p'}/2>q^{4n-8}$. So $k=0$ or
$n-1$. Modulo $Z(\widetilde{G})$, we can assume that $k=n-1$. There
are two cases:
\begin{itemize}
\item $C=(Sp_{2n-2}(q)\times GL_1(q))\cdot \ZZ_{q-1}$, which is
happened when $\Spec(s)=\{1,...,1$, $\lambda, \lambda^{-1}\}$, where
$\pm1\neq \lambda\in \FF^\times_q$. Note that there are $(q-3)/2$
choices for $\lambda$, namely, $\lambda=\gamma^k, k\in T_1$. For
each such $\lambda$, there is exactly one semi-simple conjugacy
class of such elements $s$ in $Sp_{2n}(q)$ by \cite[(2.6)]{W1}. This
conjugacy class is also the conjugacy class of $s$ in $CSp_{2n}(q)$.
Therefore, there is exactly one conjugacy class of semi-simple
elements $(s^\ast)$ in $G^\ast$ such that $\Spec(s)=\{1,...,1$,
$\lambda, \lambda^{-1}\}$ for each $\lambda=\gamma^k, k\in T_1$.
\item $C=(Sp_{2n-2}(q)\times GU_1(q)) \cdot \ZZ_{q-1}$, which is happened when
$\Spec(s)=\{1,...,1$, $\lambda, \lambda^{-1}\}$, where
$\pm1\neq \lambda\in \FF^\times_{q^2}$ and
$\lambda^{-1}=\lambda^{q}$. Note that there are $(q-1)/2$ choices
for $\lambda$, namely, $\lambda=\eta^k, k\in T_2$. Similarly as
above, there is exactly one conjugacy class of semi-simple elements
$(s^\ast)$ in $G^\ast$ such that $\Spec(s)=\{1,...,1$, $\lambda,
\lambda^{-1}\}$ for each $\lambda=\eta^k, k\in T_2$.
\end{itemize}
In two above situations, if
$\{1,...,1,\lambda,\lambda^{-1}\}=\{\mu,...,\mu,\mu\lambda,\mu\lambda^{-1}\}$
for some $\mu\in \mathbb{F}_q^\times$, then $\mu=1$ since $n\geq5$.
Therefore, $C$ is the complete inverse image of $C^\ast$ in
$\widetilde{G}$. In other words, $C^\ast=C/Z(\widetilde{G})$ and
hence $(G^\ast:C^\ast)_{p'}=(\widetilde{G}:C)_{p'}$. Consider the
canonical homomorphism $f:Sp_{2n-2}(q)\times
GL_1^{\alpha}(q)\hookrightarrow C\rightarrow C^\ast$, $\alpha=\pm$,
whose kernel is contained in the center of $Sp_{2n-2}(q)\times
GL_1^{\alpha}(q)$ and image contains the commutator group of
$C^\ast$ since $Sp_{2n-2}(q)\times GL_1^{\alpha}(q)$ contains the
commutator group of $C\simeq (Sp_{2n-2}(q)\times
GL_1^{\alpha}(q))\cdot \mathbb{Z}_{q-1}$. By Proposition \ref{DM},
the unipotent characters of $Sp_{2n-2}(q)\times GL_1^{\alpha}(q)$
are of the form $\psi\circ f$, where $\psi$ runs over the unipotent
characters of $C^\ast$. In particular, $\psi$ is trivial or
$\psi(1)\geq (q^{n-1}-1)(q^{n-1}-q)/2(q+1)$ by \cite[Proposition
5.1]{TZ1}. In the latter case,
$$\chi(1)\geq \frac{q^{2n}-1}{q+1}\cdot
\frac{(q^{n-1}-1)(q^{n-1}-q)}{2(q+1)}>q^{4n-8}.$$ Therefore, in this
case, $\chi$ is one of $(q-3)/2$ characters of degree
$(q^{2n}-1)/(q-1)$ or $(q-1)/2$ characters of degree
$(q^{2n}-1)/(q+1)$.

$3)$ Case $m=n$. Since $(s^\ast)$ is non-trivial, we assume $1\leq
k\leq n-1$.

First, if $k=1$ or $n-1$, modulo $Z(\widetilde{G})$, we may assume
$\Spec(s)=\{-1,-1,1,...,1\}$. Again, $C^\ast=C/Z(\widetilde{G})$ and
$(G^\ast:C^\ast)_{p'}=(\widetilde{G}:C)_{p'}$. There is a unique
possibility for $(s^\ast)$ in this case. We have
$\chi(1)=\frac{q^{2n}-1}{q^2-1}\psi(1)$ where $\psi$ is a unipotent
character of $C^\ast$. Applying Proposition \ref{DM} again for the
canonical homomorphism $f: Sp_2(q)\times Sp_{2n-2}(q)\hookrightarrow
C\rightarrow C^\ast$, we see that either $\psi$ is trivial, or the
unipotent character of degree $q$, or $\psi(1)\geq
(q^{n-1}-1)(q^{n-1}-q)/2(q+1)$. In the last case,
$$\chi(1)\geq \frac{q^{2n}-1}{q^2-1}\cdot
\frac{(q^{n-1}-1)(q^{n-1}-q)}{2(q+1)}>q^{4n-8}.$$ The first two
cases explain why $G$ has two characters of degrees
$(q^{2n}-1)/(q^2-1)$ and $q(q^{2n}-1)/(q^2-1)$.

Next, if $k=2$ or $n-2$, since $n\geq 5$, again one can show that
$C^\ast=C/Z(\widetilde{G})$ and hence $$\chi(1)\geq
(\widetilde{G}:C)_{p'}\geq\dfrac{(q^{2n-2}-1)(q^{2n}-1)}{(q^2-1)(q^4-1)}>q^{4n-8}.$$

Finally, if $3\leq k\leq n-3$, in particular $n\geq 6$, then
$$\chi(1)\geq (\widetilde{G}:C)_{p'}/2\geq
\frac{(q^{2n-4})(q^{2n-2}-1)(q^{2n}-1)}{2(q^2-1)(q^4-1)(q^6-1)}>
q^{4n-8}.$$ \hfill$\Box$

%%% -----------------------------------------------------------------------------------

\section{Even characteristic orthogonal groups in even dimension}

In this section we classify low-dimensional complex characters of
the group $G=Spin_{2n}^\alpha(q)=\Omega^{\alpha}_{2n}(q)$, where
$\alpha=\pm$ and $q$ is a power of $2$. Since $q$ is even, we can
identify $G$ with its dual group. Note that,
$\Omega^{\alpha}_{2n}(q)=[O^{\alpha}_{2n}(q),O^{\alpha}_{2n}(q)]$
and $O^{\alpha}_{2n}(q)=\Omega^{\alpha}_{2n}(q)\cdot \mathbb{Z}_2$.

\medskip

{\bf Proof of Theorem \ref{orthogonal-qeven}.} If $\chi$ is
unipotent then we are done by Propositions \ref{prop1}, \ref{prop2}.
Now we assume $\chi$ is not unipotent. Suppose that $\chi$ is
parameterized by $((s),\psi)$, where $s$ is a non-trivial
semi-simple element in $G$ and $\psi$ is a unipotent character of
$C:=C_{G}(s)$ such that $\chi(1)\leq q^{4n-10}$. Set
$C'=C_{O^{\alpha}_{2n}(q)}(s)$. Since $C$ is a subgroup of $C'$ of
index $1$ or $2$, $(G:C)_{2'}=(O^{\alpha}_{2n}(q):C')_{2'}$.

By Lemma \ref{lemma2}, $C'\simeq O_{2m}^{\pm}(q) \times
\prod_{i=1}^{t}GL_{a_i}^{\alpha_i}(q^{k_i})$, where $\alpha_i=\pm$,
$\Sigma_{i=1}^tk_ia_i=n-m$. Since $(s)$ is non-trivial, $m<n$. It is
easy to see that
$|\prod_{i=1}^{t}GL_{a_i}^{\alpha_i}(q^{k_i})|_{2'}\leq
\prod_{i=1}^{n-m}(q^i-(-1)^i)$. Therefore, with convention that
$q^0-1=1$, we have $$(G:C)_{2'}\geq
\frac{(q^{2(m+1)}-1)\cdot\cdot\cdot (q^{2n}-1)(q^m-1)}
{(q+1)(q^2-1)\cdot\cdot\cdot (q^{n-m}-(-1)^{n-m})(q^n+1)}=:f(m,n).$$

1) When $1\leq m\leq n-2$. We have
$$f(1,n)=(q-1)(q^n-1)\prod_{i=1}^{n-1}(q^i+(-1)^i)>q^{4n-10}$$
and
$$f(n-2,n)=\frac{(q^{2n-2}-1)(q^n-1)(q^{n-2}-1)}{(q+1)(q^2-1)}>q^{4n-10}.$$
Therefore, $\chi(1)\geq \min\{f(1,n),f(n-2,n)\}>q^{4n-10}$.

2) When $m=n-1$. Then $C\simeq O_{2n-2}^{\pm}(q)\times
GL_1^{\beta}(q)$ with $\beta=\pm$. Since $GL_1^{\beta}(q)\simeq
\Omega_2^{\beta}(q)$, $C'\simeq O_{2n-2}^{\alpha\beta}(q)\times
GL_1^{\beta}(q)$. It is obvious that
$C=\Omega_{2n-2}^{\alpha\beta}(q)\times GL_1^\beta(q)$. There are
two cases:
\begin{itemize}
  \item $C'=O_{2n-2}^{\alpha}(q)\times GL_1(q)$ corresponding to the case
 $\Spec(s)=\{1,...,1,\lambda,\lambda^{-1}\}$, where $1\neq \lambda\in \FF^\times_q$.
 Note that there are $(q-2)/2$
choices for $\lambda$, namely, $\lambda=\gamma^k, k\in T_1$. For
each such $\lambda$, there is exactly one semi-simple conjugacy
class of such elements $s$ in $O^\alpha_{2n}(q)$ by
\cite[(3.7)]{W1}. Since $(G:C)=(O_{2n}^\alpha(q):C')$,
$O_{2n}^\alpha(q)$-conjugacy class of $s$ is also a $G$-conjugacy of
$s$. Therefore, there is also exactly one conjugacy class of
semi-simple elements $(s)$ in $G$ such that $\Spec(s)=\{1,...,1$,
$\lambda, \lambda^{-1}\}$ for each $\lambda=\gamma^k, k\in T_1$.
 \item $C'=O_{2n-2}^{-\alpha}(q)\times GU_1(q)$ corresponding to the case
$\Spec(s)=\{1,...,1, \lambda, \lambda^{-1}\}$, where $1\neq
\lambda\in \FF^\times_{q^2}$ and $\lambda^{-1}=\lambda^{q}$.
Similarly as above, there is exactly one conjugacy class of
semi-simple elements $s$ in $G$ such that $\Spec(s)=\{1,...,1,
\lambda, \lambda^{-1}\}$ for each $\lambda=\eta^k, k\in T_2$.
\end{itemize}
We have $\chi(1)=\frac{(q^n-\alpha)(q^{n-1}+\alpha\beta)}{q-
\beta}\psi(1)$, where $\psi$ is a unipotent character of $C$. If
$\psi$ is non-trivial and $(n,q,\alpha)\neq (5,2,-)$ then by
Propositions 7.1 and 7.2 of \cite{TZ1},
$$\psi(1)\geq \frac{(q^{n-1}-\alpha\beta)(q^{n-2}+\alpha\beta q)}{q^2-1}.$$ In that case,
$$\chi(1)\geq \frac{(q^n-\alpha)(q^{n-1}+\alpha\beta)}{q-\beta}\cdot
\frac{(q^{n-1}-\alpha\beta)(q^{n-2}+\alpha\beta
q)}{q^2-1}>q^{4n-10}.$$ If $\psi$ is non-trivial and
$(n,q,\alpha)=(5,2,-)$, then $\psi(1)\geq q^3(q-1)^3(q^3-1)/2=28$
and we still have $\chi(1)>q^{4n-10}$. So $\psi$ must be trivial. In
summary, in this case, $\chi$ is one of $(q-2)/2$ characters of
degree $(q^{n}-\alpha)(q^{n-1}+\alpha)/(q-1)$ or $q/2$ characters of
degree $(q^{n}-\alpha)(q^{n-1}-\alpha)/(q+1)$.

3) When $m=0$. If $n\geq6$, we still have
$f(0,n)=\prod_{i=1}^n(q^i+(-1)^i)/(q^n+1)>q^{4n-10}$. Now we
consider the case $n=5$. Since $\chi(1)\leq q^{4n-10}$,
$C'=GU_5(q)$. This forces $G=\Omega_{10}^-(q)$. Then
$(G:C)_{2'}=(q-1)(q^2+1)(q^3-1)(q^4+1)$ and therefore $\psi(1)=1$.
There are exactly $q/2$ conjugacy classes $(s)$ in $O_{10}^-(q)$ of
semi-simple elements such that $C'\simeq GU_5(q)$, which is happened
when
$\Spec(s)=\{\lambda,\lambda,\lambda,\lambda,\lambda,\lambda^{-1},
\lambda^{-1},\lambda^{-1},\lambda^{-1},\lambda^{-1}\}$,
$\lambda=\eta^k, k\in T_2$. Note that $GU_5(q)\leq
\Omega_{10}^-(q)$, so $C=GU_5(q)$ and
$(O_{10}^-(q):C')=2(\Omega_{10}^-(q):C)$. In other words,
$|s^{O_{10}^-(q)}|=2|s^{\Omega_{10}^-(q)}|$. Therefore there are
exactly $q$ conjugacy classes $(s)$ in $\Omega_{10}^-(q)$ such that
$C=GU_5(q)$. This gives $q$ characters of degree
$(q-1)(q^2+1)(q^3-1)(q^4+1)$ of $\Omega_{10}^-(q)$. \hfill $\Box$

%%% -----------------------------------------------------------------------------------

\section{Odd characteristic orthogonal groups in even dimension}

In this section, we classify low-dimensional complex representations
of the group $G=Spin^{\alpha}_{2n}(q)$ where $\alpha=\pm$ and $q$ is
an odd prime power. The dual group $G^\ast$ is the projective
conformal orthogonal group $P(CO^{\alpha}_{2n}(q)^0)$, which is
quotient of $CO^{\alpha}_{2n}(q)^0$ by its center
$\mathbb{Z}(CO^\alpha_{2n}(q)^0)\simeq\mathbb{Z}_{q-1}$.

\medskip

\textbf{Proof of Theorem \ref{orthogonal-qodd}}. If $\chi$ is
unipotent then we are done by Propositions \ref{prop1}, \ref{prop2}.
So we can assume that $\chi$ is not unipotent. We denote
$\widetilde{G}:=CO_{2n}^\alpha(q)$,
$\widetilde{G}^0:=CO_{2n}^\alpha(q)^0$, and
$Z:=Z(\widetilde{G})\simeq Z(\widetilde{G}^0)\simeq
\mathbb{Z}_{q-1}$. Suppose that $s^\ast$ is a non-trivial
semi-simple element in $G^\ast$. Consider an inverse image $s$ of
$s^\ast$ in $\widetilde{G}^0$. Then $(s)$ is a non-trivial conjugacy
class of semi-simple elements in $\widetilde{G}$. Set
$C^\ast:=C_{G^\ast}(s^\ast)$, $\widetilde{C}:=C_{\widetilde{G}}(s)$,
and $\widetilde{C}^0:=C_{\widetilde{G}^0}(s)$. Suppose that $\chi$
is parameterized by $((s^\ast),\psi), \psi\in \Irr(C^\ast)$, such
that $\chi(1)\leq q^{4n-10}$. Let $g$ be any element in
$\widetilde{G}^0$ such that the image of $g$ in $G^\ast$ belongs to
$C^\ast$. Then $gsg^{-1}=\mu(g)s$ for some $\mu(g)\in
\mathbb{Z}_{q-1}$. Therefore,
$\tau(s)=\tau(gsg^{-1})=\tau(\mu(g)s)=\mu(g)^2\tau(s)$. It follows
that $\mu(g)$ is $1$ or $-1$. Hence,
$$(G^\ast:C^\ast)_{p'}\geq (\widetilde{G}^0:\widetilde{C}^0)_{p'}/2\geq
(\widetilde{G}:\widetilde{C})_{p'}/4.$$

\textbf{Case} $\mathbf{1}$: If $\tau(s)$ is not a square in
$\mathbb{F}_q$, multiplying $s$ by a suitable scalar in
$\mathbb{F}_q$, we can assume that $\tau(s)$ is fixed. By Lemma
\ref{lemma4}, $\widetilde{C}\simeq(O^\pm_{m}(q^2)\times
\prod_{i=1}^{t}GL_{a_i}^{\alpha_i}(q^{k_i}))\cdot \mathbb{Z}_{q-1}$,
where $\alpha_i=\pm$, $\Sigma_{i=1}^tk_ia_i=n-m$. Denote the group
$\prod_{i=1}^{t}GL_{a_i}^{\alpha_i}(q^{k_i})$ by $A$ for short. Now
we will show that, when $n\geq6$,
$(\widetilde{G}:\widetilde{C})_{p'}>4q^{4n-10}$, and therefore
$\chi(1)\geq(G^\ast:C^\ast)_{p'}\geq
(\widetilde{G}:\widetilde{C})_{p'}/4>q^{4n-10}$. It is easy to see
that $|A|_{p'}\leq \prod_{i=1}^{n-m}(q^i-(-1)^i)$. Hence,
$$(\widetilde{G}:\widetilde{C})_{p'}\geq
\frac{2(q^{2}-1)\cdot\cdot\cdot (q^{2n}-1)}
{(q^n+1)(q+1)(q^2-1)\cdot\cdot\cdot
(q^{n-m}-(-1)^{n-m})|O^\pm_{m}(q^2|}=:f(m,n).$$ We have
$f(0,n)=2\prod_{i=1}^n(q^i+(-1)^i)/(q^n+1)>4q^{4n-10}$ and
$f(n,n)>4q^{4n-10}$ for every $n\geq 6$. Therefore,
$\chi(1)\geq(\widetilde{G}:\widetilde{C})_{p'}/4\geq \min\{f(0,n),
f(n,n)\}/4>q^{4n-10}$ when $n\geq 6$.

Now let us consider the case $(n,\alpha)=(5,+)$. Note that $GU_5(q)$
is not a subgroup of $O_{10}^+(q)$. Therefore, if $A\neq GL_5(q)$,
one can show that $(\widetilde{G}:\widetilde{C})_{p'}> 4q^{10}$ and
hence $\chi(1)> q^{10}$. If $A=GL_5(q)$, then the characteristic
polynomial of $s$ has the form
$P(x)=(x-\lambda)^5(x-\tau\lambda^{-1})^5$ with $\lambda\in
\mathbb{F}_q^\ast$ and $\tau:=\tau(s)$. Now we will show that
$(G^\ast:C^\ast)_{p'}= (\widetilde{G}^0:\widetilde{C}^0)_{p'}$ and
therefore
$\chi(1)\geq(\widetilde{G}:\widetilde{C})_{p'}/2=(q+1)(q^2+1)(q^3+1)(q^4+1)>q^{10}$.
Suppose that $s_1,s_2\in \widetilde{G}^0$ are
$\widetilde{G}^0$-conjugate and their images in $G^\ast$ are the
same. Then $s_1=\pm s_2$. Let $V_1:=\Ker(s-\lambda)$ and
$V_2:=\Ker(s-\tau\lambda^{-1})$. Note that $V_1$ and $V_2$ are
totally isotropic subspaces in $V$ and $V_1\cap V_2=\{0\}$. By
definition, an element in $\widetilde{G}^0$ cannot carry $V_1$ to
$V_2$. This implies that $s_1=s_2$. In other words,
$|{s^\ast}^{G^\ast}|=|s^{\widetilde{G}^0}|$ and hence
$(G^\ast:C^\ast)_{p'}= (\widetilde{G}^0:\widetilde{C}^0)_{p'}$ as
desired.

The next case is $(n,\alpha)=(5,-)$. Note that $GL_5(q)$ is not a
subgroup of $O_{10}^-(q)$. Therefore, if $A\neq GU_5(q)$, one again
can show that $\chi(1)>q^{10}$. If $A=GU_5(q)$, then the
characteristic polynomial of $s$ has the form
$P(x)=[(x-\lambda)(x-\tau\lambda^{-1})]^5$, where
$(x-\lambda)(x-\tau\lambda^{-1})$ is an irreducible polynomial over
$\mathbb{F}_q$. Note that there are $(q+1)/2$ choices for such a
pair $(\lambda, \tau\lambda^{-1})$. So there are exactly $(q+1)/2$
conjugacy classes $(s)$ in $\widetilde{G}$ such that
$\widetilde{C}\simeq GU_5(q)$. That means there are exactly $(q+1)$
such conjugacy classes $(s)$ in $\widetilde{G}^0$ since
$\widetilde{C}^0= \widetilde{C}$ in this case. Now modulo $Z$, we
have exactly $(q+1)/2$ conjugacy classes $(s^\ast)$ in $G^\ast$.
Arguing similarly as above, one has $(G^\ast:C^\ast)_{p'}=
(\widetilde{G}^0:\widetilde{C}^0)_{p'}=(q-1)(q^2+1)(q^3-1)(q^4+1)$
and therefore $C^\ast=\widetilde{C}^0/Z$. The condition
$\chi(1)=(G^\ast:C^\ast)_{p'}\psi(1)\leq q^{10}$ implies that
$\psi(1)=1$. Since $C^\ast$ has a unique unipotent character of
degree $1$, we get exactly $(q+1)/2$ characters of degree
$(q-1)(q^2+1)(q^3-1)(q^4+1)$ in this case.

\textbf{Case} $\mathbf{2}$: If $\tau(s)$ is a square in
$\mathbb{F}_q$, again, we can assume that $\tau(s)=1$. By Lemma
\ref{lemma4}, $\widetilde{C}\simeq (O^\pm_{2k}(q)\times
O^\pm_{2m-2k}(q) \times
\prod_{i=1}^{t}GL_{a_i}^{\alpha_i}(q^{k_i}))\cdot \mathbb{Z}_{q-1}$,
where $\alpha_i=\pm$, $\Sigma_{i=1}^tk_ia_i=n-m$.

It is easy to see that $|O^\pm_{2k}(q)\times
O^\pm_{2m-2k}(q)|_{p'}\leq 2\prod_{i=1}^{m}(q^{2i}-1)/(q^m-1)$ if
$m\geq 1$ and $|\prod_{i=1}^{t}GL_{a_i}^{\alpha_i}(q^{k_i}|_{p'}\leq
\prod_{i=1}^{n-m}(q^i-(-1)^i)$. Therefore, with convention that
$q^0-1=2$, we have $$(\widetilde{G}:\widetilde{C})_{p'}\geq
\frac{(q^{2(m+1)}-1)\cdot\cdot\cdot (q^{2n}-1)(q^m-1)}
{(q^n+1)(q+1)(q^2-1)\cdot\cdot\cdot (q^{n-m}-(-1)^{n-m})}=:g(m,n).
$$

$1)$ When $1\leq m\leq n-2$, $(\widetilde{G}:C)_{p'}\geq
\min\{g(1,n), g(n-2,n)\}$. Direct computation shows that
$g(1,n)=(q^n-1)\prod_{i=1}^{n-1}(q^i+(-1)^i)/(q+1)> 4q^{4n-10})$ and
$g(n-2,n)=(q^{2n-2}-1)(q^{n}-1)(q^{n-2}-1)/(q+1)(q^2-1)>4q^{4n-10}$.
Therefore, $\chi(1)>q^{4n-10}$.

$2)$ When $m=0$, then $g(0,n)=2\prod_{i=1}^{n}(q^i+(-1)^i)/(q^n+1)$.
We still have $g(0,n)>4q^{4n-10}$ if $n\geq6$. The case
$(n,\alpha)=(5,+)$ can be argued similarly as when $\tau(s)$ is not
a square in $\mathbb{F}_q$. So we only need to consider the case
$(n,\alpha)=(5,-)$. Note that $GL_5(q)$ is not a subgroup of
$O_{10}^-(q)$. Therefore, if $A\neq GU_5(q)$, one again can show
that $\chi(1)>q^{10}$. So we can assume $A=GU_5(q)$. Then the
characteristic polynomial of $s$ has the form
$P(x)=[(x-\lambda)(x-\lambda^{-1})]^5$, where
$(x-\lambda)(x-\lambda^{-1})$ is an irreducible polynomial over
$\mathbb{F}_q$. Note that there are exactly $(q-1)/2$ such a pair
$(\lambda, \lambda^{-1})$. Repeating arguments as when $\tau(s)$ is
not a square in $\mathbb{F}_q$, we get exactly $(q-1)/2$ characters
of degree $(q-1)(q^2+1)(q^3-1)(q^4+1)$ in this case.

$3)$ When $m=n-1$. If $1\leq k\leq n-2$ then
$(\widetilde{G}:\widetilde{C})_{p'}\geq
(q^{2n-2}-1)(q^{n}-1)(q^{n-2}-1)/2(q+1)^2>4q^{4n-10}$. Hence,
$\chi(1)\geq (G^\ast:C^\ast)_{p'}>q^{4n-10}$. So $k=0$ or $n-1$.
With no loss, we can assume $k=n-1$. Note that $\Spec(s)=\{1,...,1$,
$\lambda, \lambda^{-1}\}$ in this case. If
$\{1,...,1,\lambda,\lambda^{-1}\}=\{\mu,...,\mu,\mu\lambda,\mu\lambda^{-1}\}$
for some $\mu\in \mathbb{F}_q^\times$, then $\mu=1$ since $n\geq5$.
Therefore, $\widetilde{C}^0$ is the complete inverse image of
$C^\ast$ in $\widetilde{G}^0$. In other words,
$C^\ast=\widetilde{C}^0/Z$ and hence
$(G^\ast:C^\ast)_{p'}=(\widetilde{G}^0:\widetilde{C}^0)_{p'}=
(\widetilde{G}:\widetilde{C})_{p'}$. Also, since
$GL_1^{\beta}(q)\simeq \Omega_2^{\beta}(q)$, $\widetilde{C}\simeq
(O_{2n-2}^{\alpha\beta}(q)\times GL_1^{\beta}(q))\cdot
\mathbb{Z}_{q-1}$. One can show that $\widetilde{C}^0\simeq
(SO_{2n-2}^{\alpha\beta}(q)\times GL_1^{\beta}(q))\cdot
\mathbb{Z}_{q-1}$. There are two cases:
\begin{itemize}
\item $\widetilde{C}=(O^\alpha_{2n-2}(q)\times GL_1(q))\cdot \ZZ_{q-1}$
corresponding to the case $\Spec(s)=\{1,...,1$, $\lambda,
\lambda^{-1}\}$, where $\pm1\neq \lambda\in \FF^\times_q$. Note that
there are $(q-3)/2$ choices for $\lambda$, namely,
$\lambda=\gamma^k, k\in T_1$. For each such $\lambda$, there is
exactly one semi-simple conjugacy class of such elements $s$ in
$O^\alpha_{2n}(q)$ by \cite[(2.6)]{W1}. This conjugacy class is also
the conjugacy class of $s$ in $\widetilde{G}^0$. Therefore, there is
exactly one conjugacy class of semi-simple elements $(s^\ast)$ in
$G^\ast$ such that $\Spec(s)=\{1,...,1$, $\lambda, \lambda^{-1}\}$
for each $\lambda=\gamma^k, k\in T_1$.
\item $\widetilde{C}=(O^{-\alpha}_{2n-2}(q)\times GU_1(q)) \cdot \ZZ_{q-1}$
corresponding to the case $\Spec(s)=\{1,...,1$, $\lambda,
\lambda^{-1}\}$, where $\pm1\neq\lambda\in \FF^\times_{q^2}$ and
$\lambda^{-1}=\lambda^{q}$. Note that there are $(q-1)/2$ choices
for $\lambda$, namely, $\lambda=\eta^k, k\in T_2$. Similarly as
above, there is exactly one conjugacy class of semi-simple elements
$(s^\ast)$ in $G^\ast$ such that $\Spec(s)=\{1,...,1$, $\lambda,
\lambda^{-1}\}$ for each $\lambda=\eta^k, k\in T_2$.
\end{itemize}
Using Proposition \ref{DM} for the canonical homomorphism
$f:SO_{2n-2}^{\alpha\beta}(q)\times GL_1^{\beta}(q)\hookrightarrow
\widetilde{C}^0\rightarrow C^\ast$, $\beta=\pm$, we see that the
unipotent characters of $SO^{\alpha\beta}_{2n-2}(q)\times
GL_1^{\beta}(q)$ are of the form $\psi\circ f$, where $\psi$ runs
over the unipotent characters of $C^\ast$. In particular, by
Propositions 7.1, 7.2 of \cite{TZ1}, $\psi$ is trivial or
$\psi(1)\geq (q^{n-1}-\alpha\beta)(q^{n-2}+\alpha\beta q)/(q^2-1)$.
In the latter case,
$$\chi(1)\geq \frac{(q^{n}-\alpha)(q^{n-1}+\alpha\beta)}{q-\beta}\cdot
\frac{(q^{n-1}-\alpha\beta)(q^{n-2}+\alpha\beta
q)}{q^2-1}>q^{4n-10}.$$ Therefore, in this case, $\chi$ is one of
$(q-3)/2$ representations of degree
$(q^{n}-\alpha)(q^{n-1}+\alpha)/(q-1)$ or $(q-1)/2$ representations
of degree $(q^{n}-\alpha)(q^{n-1}-\alpha)/(q+1)$.

$4)$ When $m=n$. Since $(s^\ast)$ is non-trivial, we assume $1\leq
k\leq n-1$.

If $k=1$ or $n-1$, modulo $Z$, we may assume
$\Spec(s)=\{-1,-1,1,...,1\}$. Then
$\widetilde{C}=(O_2^\beta(q)\times O_{2n-2}^{\alpha\beta}(q))\cdot
\mathbb{Z}_{q-1}$ and $\widetilde{C}^0\simeq (SO_2^\beta(q)\times
SO_{2n-2}^{\alpha\beta}(q))\cdot 2\cdot \mathbb{Z}_{q-1}$. Again,
$C^\ast=\widetilde{C}^0/Z$ and
$(G^\ast:C^\ast)_{p'}=(\widetilde{G}:\widetilde{C})_{p'}$. There is
a unique possibility for $(s^\ast)$ for each $\beta=\pm$. We have
$\chi(1)=\frac{(q^{n}-\alpha)(q^{n-1}+\alpha\beta)}{2(q-\beta)}\psi(1)$
where $\psi$ is a unipotent character of $C^\ast$. Consider the
canonical homomorphism $f: (SO^\beta_2(q)\times
SO^{\alpha\beta}_{2n-2}(q))\cdot 2\hookrightarrow
\widetilde{C}^0\rightarrow C^\ast$, whose kernel is contained in the
center of $(SO^{\beta}_2(q)\times SO^{\alpha\beta}_{2n-2}(q))\cdot2$
and image contains the commutator group of $C^\ast$ since
$(SO^{\beta}_2(q)\times SO^{\alpha\beta}_{2n-2}(q))\cdot 2$ contains
the commutator group of $\widetilde{C}^0$. By Proposition \ref{DM},
the unipotent characters of $(SO^\beta_2(q)\times
SO^{\alpha\beta}_{2n-2}(q))\cdot2$ are of the form $\psi\circ f$,
where $\psi$ runs over the unipotent characters of $C^\ast$. In
particular, $\psi$ is one of two linear unipotent characters or
$\psi(1)\geq (q^{n-1}-\alpha\beta)(q^{n-2}+\alpha\beta q)/(q^2-1)$.
In the latter case,
$$\chi(1)\geq \frac{(q^{n}-\alpha)(q^{n-1}+\alpha\beta)}{2(q-\beta)}\cdot
\frac{(q^{n-1}-\alpha\beta)(q^{n-2}+\alpha\beta
q)}{q^2-1}>q^{4n-10}.$$

If $2\leq k\leq n-2$, then $$\chi(1)\geq
\frac{(\widetilde{G}:\widetilde{C})_{p'}}{4}\geq
\frac{|O_{2n}^\alpha(q)|}{4|O^\beta_4(q)|\cdot|O_{2n-4}^{\alpha\beta}(q)|}=$$$$=
\frac{(q^{2n-2}-1)(q^n-\alpha)(q^{n-2}+\alpha\beta)(q^2+\beta)}
{8(q^2-1)(q^4-1)}>q^{4n-10}.$$ \hfill$\Box$

%%% -----------------------------------------------------------------------------------

\section{Groups $Spin_{12}^\pm(3)$}

In this section, we will classify the irreducible complex characters
of $G=Spin_{12}^\alpha(3)$, $\alpha=\pm$, of degrees up to $4\cdot
3^{15}$. Note that $G^\ast=P(CO_{12}^\alpha(3)^0)$. The argument in
this section is similar to that in $\S7$. We will keep all notation
from there. Two following lemmas can be checked by direct
computation.

\begin{lemma}\label{prop-}
Suppose that $\chi\in \Irr(P(CO_{12}^-(3)^0))$ is unipotent. Then
either $\chi$ is one of the characters labeled by $0 \hspace{3pt}
6\choose-$, $1 \hspace{3pt} 5\choose-$, $2\hspace{3pt} 4\choose-$,
$0\hspace{3pt} 1 \hspace{3pt} 6\choose1$, $0 \hspace{3pt} 1
\hspace{3pt} 2\choose5$, or $\chi(1)> 4\cdot3^{15}$.
\end{lemma}

\begin{lemma}\label{prop+}
Suppose that $\chi\in \Irr(P(CO_{12}^+(3)^0))$ is unipotent. Then
either $\chi$ is one of characters labeled by $6\choose0$,
$5\choose1$, $4\choose2$, $3\choose3$ ($2$ characters), $0
\hspace{3pt}1\choose1\hspace{3pt}6$, $0\hspace{3pt}1
\hspace{3pt}2\hspace{3pt}5\choose-$, or $\chi(1)> 4\cdot 3^{15}$.
\end{lemma}

\begin{propo}\label{smallcase}
$Spin_{12}^+(3)$ has exactly $28$ irreducible complex characters of
degrees less than $4\cdot 3^{15}$ and $Spin_{12}^-(3)$ has exactly
$16$ irreducible complex characters of degrees less than $4\cdot
3^{15}$.
\end{propo}

\begin{proof}
Let $\chi$ be an irreducible complex character of $G$ of degree less
than $4\cdot 3^{15}$. Since we have already counted the number of
unipotent characters in two above lemmas, we can assume that $\chi$
is not unipotent. Suppose that $\chi$ is parameterized by
$((s^\ast),\psi)$, where $1\neq s^\ast\in
G^\ast:=P(CO_{12}^\alpha(3)^0)$ and $\psi\in\Irr(C^\ast)$,
$C^\ast=C_{G^\ast}(s^\ast)$. Recall that
$\widetilde{G}:=CO_{12}^\alpha(3)$ and
$\widetilde{G}^0:=CO_{12}^\alpha(3)^0$. Denote
$Z:=Z(\widetilde{G})=Z(\widetilde{G}^0)$ the center of
$\widetilde{G}$ as well as $\widetilde{G}^0$. Let $s$ be an inverse
image of $s^\ast$ in $\widetilde{G}^0$. Set
$\widetilde{C}=C_{\widetilde{G}}(s)$ and
$\widetilde{C}^0=C_{\widetilde{G}^0}(s)$.

\textbf{Case} $\mathbf{1}$: $\tau(s)=-1$. Then we have
$\widetilde{C}\simeq (O^\pm_{m}(9)\times
\prod_{i=1}^{t}GL_{a_i}^{\alpha_i}(3^{k_i}))\cdot \mathbb{Z}_{2}$,
where $\alpha_i=\pm$, $\Sigma_{i=1}^tk_ia_i+m=6$. Since
$\chi(1)<4\cdot3^{15}$, $(\widetilde{G}:\widetilde{C})_{3'}<16\cdot
3^{15}$. This inequality happens only when $\widetilde{C}\simeq
GL_6^{\pm}(3)\cdot \mathbb{Z}_2$. This forces $G=Spin_{12}^+(3)$.
Let us consider the case $\widetilde{C}\simeq GL_6(3)\cdot
\mathbb{Z}_2$. Note that there is a unique conjugacy class $(s)$ of
semi-simple elements in $\widetilde{G}$ so that $\widetilde{C}\simeq
GL_6(3)\cdot \mathbb{Z}_2$, which is happened when the
characteristic polynomial of $s$ is $(x^2-1)^6$. In this case,
$\widetilde{C}^0=\widetilde{C}$,
$(\widetilde{G}:\widetilde{C})=2(\widetilde{G}^0:\widetilde{C}^0)$,
and therefore there are two conjugacy classes of semi-simple
elements in $\widetilde{G}^0$, as well as in $G^\ast$, such that
$\widetilde{C}\simeq GL_6(3)\cdot \mathbb{Z}_2$. Furthermore,
$C^\ast$ is an extension by $2$ of $\widetilde{C}/Z$. Therefore,
$\chi$ is one of $4$ characters of degree
$(G^\ast:C^\ast)_{3'}=\frac{1}{2}\prod_{i=1}^5(3^i+1)$, provided
that $G=Spin_{12}^+(3)$.

Next, we consider the case $\widetilde{C}\simeq GU_6(3)\cdot
\mathbb{Z}_2$. Note that there are two semi-simple conjugacy classes
$(s)$ in $\widetilde{G}$ so that $\widetilde{C}\simeq GL_6(3)\cdot
\mathbb{Z}_2$, which is happened when the characteristic polynomial
of $s$ is $(x^2+x-1)^6$ or $(x^2-x-1)^6$. In this case,
$\widetilde{C}^0=\widetilde{C}$,
$(\widetilde{G}:\widetilde{C})=2(\widetilde{G}^0:\widetilde{C}^0)$,
and therefore there are $4$ such conjugacy classes of semi-simple
elements in $\widetilde{G}^0$. Modulo $Z$, we get exactly two
semi-simple conjugacy classes $(s^\ast)$ in $G^\ast$. Furthermore,
$C^\ast\simeq\widetilde{C}/Z$. Therefore, we get two characters of
degree $\prod_{i=1}^5(3^i+(-1)^i)$ of $G=Spin_{12}^+(3)$.

\textbf{Case} $\mathbf{2}$: $\tau(s)=1$. We have
$\widetilde{C}\simeq (O^\pm_{2k}(3)\times O^\pm_{2m-2k}(3)\times
\prod_{i=1}^{t}GL_{a_i}^{\alpha_i}(3^{k_i}))\cdot \mathbb{Z}_{2}$,
where $\alpha_i=\pm$, $\Sigma_{i=1}^tk_ia_i+m=6$. The inequality
$\chi(1)<4\cdot3^{15}$ implies that $m=0, 5$, or $6$. If $m=0$,
similarly as in case $1$, one can show that $\widetilde{C}\simeq
GU_6(3)\cdot \mathbb{Z}_2$ and $C^\ast$ is an extension by $2$ of
$\widetilde{C}/Z$. Hence, $\chi$ is one of $4$ characters of degree
$\frac{1}{2}\prod_{i=1}^5(3^i+(-1)^i)$, provided that
$G=Spin_{12}^+(3)$.

When $m=5$, it is easy to show that
$(G^\ast:C^\ast)_{3'}=(\widetilde{G}:\widetilde{C})_{3'}$. If $k=0$
or $5$, we have already shown in the proof of Theorem
\ref{orthogonal-qodd} that either $\chi$ is the unique character of
degree $(3^6-\alpha)(3^5-\alpha)/4$ or
$$\chi(1)\geq \frac{(3^6-\alpha)(3^5-\alpha)}{4}\cdot
\frac{(3^5+\alpha)(3^4-\alpha3)}{8}>4\cdot3^{15}.$$ If $1\leq
k\leq4$, $\chi(1)\geq(\widetilde{G}:\widetilde{C})_{3'}\geq
(3^{10}-1)(3^6-1)(3^4-1)/32>4\cdot3^{15}$.

Now we suppose $m=6$. First, If $k=1$ or $5$, we have
$\chi=\psi(1)(3^6-\alpha)(3^5+\alpha\beta)/2(3-\beta)$, where
$\beta=\pm$ and $\psi$ is a unipotent character of
$C^\ast=\widetilde{C}^0/Z$, $\widetilde{C}^0\simeq
(SO_2^\beta(3)\times SO_{10}^{\alpha\beta}(3))\cdot2\cdot
\mathbb{Z}_2$. When $\beta=+$, the inequality $\chi(1)<4\cdot
3^{15}$ forces $\psi$ to be linear. Therefore, $\chi$ is one of $2$
characters of degree $(3^6-\alpha)(3^5+\alpha)/4$. When $\beta=-$,
note that the unipotent characters of $(SO_2^-(3)\times
SO_{10}^{-\alpha}(3))\cdot2 $ are of the form $\psi\circ f$, where
$f$ is the canonical homomorphism $f: (SO_2^-(3)\times
SO_{10}^{-\alpha}(3))\cdot2 \rightarrow C^\ast$. Since $SO_2^-(3)$
has a unique unipotent character which is trivial, the unipotent
characters of $(SO_2^-(3)\times SO_{10}^{-\alpha}(3))\cdot2$ are
actually unipotent characters of $O_{10}^{-\alpha}(3)$. Therefore,
$\psi$ is one of $2$ unipotent characters of degree $1$; $2$
unipotent characters of degree $(3^5+\alpha)(3^4-\alpha 3)/8$; or
$\psi(1)\geq (3^{10}-3^2)/8$. Hence there are $4$ characters of
degrees less than $4\cdot 3^{15}$ in this case. Next, if $k=2$ or
$4$, we have
$$\chi(1)=\frac{(3^4+\alpha\beta)(3^{10}-1)(3^6-\alpha)}
{2(3^2-1)(3^2-\beta)}\cdot \psi(1),$$ where $\psi$ is a unipotent
character of $C^\ast=\widetilde{C}^0/Z(\widetilde{G})$,
$\widetilde{C}^0\simeq (SO_4^\beta(3)\times
SO_8^{\alpha\beta}(3))\cdot 2\cdot\mathbb{Z}_2$. Since
$\chi(1)<4\cdot 3^{15}$, $\psi$ is one of two linear unipotent
characters of $C^\ast$ for each $\beta=\pm$. This gives $4$ more
characters of degrees less than $4\cdot 3^{15}$. Finally, if $k=3$,
then $\chi(1)\geq(G^\ast:C^\ast)_{3'}\geq
(\widetilde{G}:\widetilde{C})_{3'}/2>4\cdot 3^{15}$.
\end{proof}

%%% -----------------------------------------------------------------------------------

\end{document}